\documentclass[12pt]{amsart}

\usepackage{amscd,amssymb}
\theoremstyle{plain}
\newtheorem{Thm}[equation]{Theorem}

\newtheorem{Con}[equation]{Conjecture}\newtheorem{Rmk}[equation]{Remark}
\newtheorem{Cor}[equation]{Corollary}
\newtheorem{Prop}[equation]{Proposition}

\newtheorem{Lem}[equation]{Lemma}
\newtheorem{Def}[equation]{Definition}
\numberwithin{equation}{section}

\newcommand{\mc}[1]{{}}

\newcommand{\q}{\mathbb{Q}}

\newcommand{\e}{\epsilon}
\newcommand{\z}{\mathbb{Z}}
\renewcommand{\q}{\mathbb{Q}}
\newcommand{\n}{\mathbb{N}}
\newcommand{\br}{\mathbb{R}}

\newcommand{\A}{\mathbb{A}}
\newcommand{\ba}{\backslash}

\newcommand{\G}{\Gamma}

\newcommand{\Cal}{\mathcal}

\newcommand{\bG}{\mathbf G}

\newcommand{\pr}{\rm{pr}}
\newcommand{\bM}{\mathbf M}
\newcommand{\SL}{\operatorname{SL}}\newcommand{\GL}{\operatorname{GL}}
\newcommand{\PGL}{\operatorname{PGL}}

\newcommand{\Sp}{\operatorname{Sp}}
\newcommand{\M}{\operatorname{M}}

\newcommand{\bc}{\mathbb C}
\newcommand{\bL}{\mathbf L}\newcommand{\bT}{\mathbf T}
\newcommand{\bhi}{G_{W_f}}

\newcommand{\bi}{\begin{itemize}}

\newcommand{\ei}{\end{itemize}}

\newcommand{\op}{\operatorname}

\renewcommand{\deg}{\text{deg}}
\newcommand{\vs}{\vskip 5pt}
\newcommand{\bga}{\bG(\A)}
\newcommand{\He}{\operatorname{H}}
\newcommand{\bH}{\mathbf H}\newcommand{\bN}{\mathbf N}
\begin{document}

\title[Rational points and Adelic mixing]
{Manin's and Peyre's conjectures on rational points and adelic mixing}
\author{Alex Gorodnik, Fran\c cois Maucourant, Hee Oh}
\dedicatory{Dedicated to Prof. Gregory Margulis  on the occasion of his sixtieth birthday}
\address{Mathemati 253-37, Caltech, Pasadena, CA 91125, USA \& (current:)
School of Mathematics, University of Bristol, Bristol BS8 1TW, U.K.//
IRMAR-Campus de Beaulieu, 35042 RENNES Cedex, France
//Mathemati 253-37, Caltech, Pasadena, CA 91125, USA \&
(current: )  Box 1917, 151 Thayer Street,
Brown university, Providence, RI 02912}

\abstract{
Let $X$ be the wonderful compactification of a
connected adjoint semisimple group $G$ defined over a number field $K$.
We prove Manin's conjecture on the asymptotic (as $T\to \infty$) of the number
of $K$-rational points of $X$ of height less than $T$, and give an explicit construction of a measure
on $X(\A)$, generalizing Peyre's measure,
which describes the asymptotic distribution of the rational points
$\bG(K)$ on $X(\A)$.  Our approach is based on the mixing property of 
${L}^2(\bG(K)\backslash \bG(\A))$ which
 we obtain with a rate of convergence.

\vs\vs\vs Soit $X$ la compactification merveilleuse d'un groupe
semisimple $\bG$, connexe, de type adjoint, alg\'ebrique d\'efini sur un corps de nombre $K$. Nous d\'emontrons l'asymptotique conjectur\'ee par Manin du nombre de points $K$-rationnels sur $X$ de hauteur plus petite que $T$, lorsque $T\rightarrow +\infty$, et construisons de mani\`ere explicite une mesure sur $X(\A)$, g\'en\'eralisant celle de Peyre, qui d\'ecrit la r\'epartition asymptotique des points rationnels $\bG(K)$ sur $X(\A)$. Ce travail repose sur la propri\'et\'e de m\'elange de $L^2(\bG(K)\backslash \bga)$, qui est d\'emontr\'ee avec une estim\'ee de vitesse.}\endabstract

\thanks{The first and the third authors
are partially supported by DMS-0400631, and DMS-0333397 and DMS-0629322 respectively.
The second author would like to thank Caltech for
the hospitality during his visit where the work was first conceived}


\email{a.gorodnik@bristol.ac.uk}
\email{francois.maucourant@univ-rennes1.fr}
\email{heeoh@math.brown.edu}
\maketitle
\tableofcontents
\section{Introduction}
Let $K$ be a number field and $X$ a smooth projective variety 
defined over $K$.
A fundamental problem in
modern algebro-arithmetic geometry is to describe the set $X(K)$ in terms of the geometric invariants of $X$.
One of the main conjectures in this area was made by 
Manin in the late eighties in \cite{BM}.
It formulates the asymptotic (as $T\to \infty$) of
the number of points in $X(K)$ of height less than $T$ for Fano varieties (that is, varieties with
ample anti-canonical class).

 Manin's conjecture has been proved for
flag varieties (\cite{FMT}, \cite{Pe1}), toric varieties (\cite{BT1}, \cite{BT2}), horospherical varieties \cite{ST},
equivariant compactifications of unipotent groups
 (see \cite{CT2}, \cite{ST1}, \cite{ST2}), etc. We refer to survey papers
 by Tschinkel (\cite{T}, \cite{T2}) for a more precise background on this conjecture.
Recently Shalika, Tschinkel and Takloo-Bighash proved the conjecture
for the wonderful  compactification of 
a connected semisimple adjoint group \cite{STT2}.
In this paper, we present a different proof of the conjecture, as well as
 describe the asymptotic distributions of
rational points of bounded height as conjectured by Peyre.
 Our proof relies on the computation of the volume asymptotics of height balls in \cite{STT2}.
We refer to \cite{Bi} for the comparison of these two approaches.

Although our work is highly motivated by conjectures in arithmetic geometry, our approach is almost purely (algebraic)
group theoretic. For this reason, we formulate our main results in the language of 
algebraic groups and their representations in the introduction and refer to
section \ref{manin} for the account that how these
results imply the conjectures of Manin and Peyre.

\subsection{Height function}
We begin by defining the notion of a height function
on the $K$-rational points of the projective $n$-space $\mathbb P^n$.
Intuitively speaking, the height of a rational point $x\in \mathbb P^n(K)$
measures an {{\it arithmetic} size of $x$.
In the case of $K=\q$,
it is simply given by $$\He(x)=\max_{0\le i \le n} |x_i|$$ 
where $(x_0, \cdots, x_n)$
is a primitive integral vector representing $x$.
To give its definition for a general $K$, we denote by
 $R$ the set of all normalized
absolute values $x\mapsto |x|_v$ of $K$, and
by $K_v$ the completion of $K$ with respect to $|\cdot|_v$.
For each $v\in R$, choose a norm $\He_v$ on $K_v^{n+1}$
which is simply the max norm $\He_v(x_0,
\cdots, x_n)=\max_{i=0}^n|x_i|_v$
for almost all $v$.
Then a function $\He: \mathbb P^n(K)\to \mathbb R_{>0}$ of the following form
is called a height function:
 $$\He(x):=\prod_{v\in R}\He_v(x_0, \cdots, x_n)$$ 
for $x=(x_0:\cdots : x_n)\in \mathbb P^n(K)$.
Since $\He_v(x_0, \cdots, x_n)=1$ for almost all $v\in R$,
we have $\He(x)>0$ and by the product formula,
$\He$ is well defined, i.e., independent of the choice of
representative for $x$.

It is easy to see that for any $T>0$, the number
$$N(T):=\#\{x\in \mathbb P^n(K): \He(x)<T\} $$
is finite. Schanuel \cite{Sc} computed the precise asymptotic in 1964:
 $$N(T)\sim c\cdot  T^{n+1}\quad\text{  as $T\to \infty$}$$ for some
  explicit constant $c=c(\He)>0$.


Unless mentioned otherwise,
throughout the introduction,
we let
 $\bG$ be a connected semisimple adjoint group over $K$ 
and $\iota: \bG \to \GL_N$ be a faithful representation of $\bG$ defined
over $K$ with a unique maximal weight.
 Consider the projective embedding of $\bG$ over $K$ induced by $\iota$: $$\bar
\iota: \bG \to \mathbb P(\M_N)$$ where $\M_N$
denotes the space of matrices of order $N$.
We then define a height function $\He_\iota$ on $\bG(K)$ associated to
$\iota$ by pulling back a height function on $\mathbb P(\M_N(K))$
 via $\bar \iota$.
That is, for $g\in \bG(K)$,
\begin{equation}\label{h1}\He_\iota(g):=
\prod_{v\in R} \He_v
(\iota(g)),\end{equation}where each $\He_v$ is a norm on $\M_N(K_v)$,
 which is the max norm for almost all $v\in R$. 

We note that $\He_\iota$ is not uniquely determined by $\iota$,
because of the freedom of choosing $\He_v$ locally (though only
for finitely many $v$).
\subsection{Asymptotic number of rational points}\label{in:a}
For each $T>0$, we introduce the notation for the number
of points in $\bG(K)$ of height less than $T$:
$$N(\He_\iota, T):=\# \{g\in \bG(K)\;:\; \He_\iota(g)\, < \,T\}.$$

\begin{Thm}\label{mt} There exist $a_\iota\in \q^+$, $b_\iota\in \n$ and
$c=c(\He_\iota)>0$ such that for some $\delta>0$,
$$N(\He_\iota, T)
= c\cdot T^{a_\iota} (\log T)^{b_\iota-1} \cdot (1+O((\log T)^{-\delta})). $$
\end{Thm}

 

The constants $a_\iota$ and $b_\iota$ can be defined explicitly by
combinatorial data on the root system of $\bG$ and the
unique maximal weight of $\iota$. Choose 
a maximal torus $\mathbf T$ of $\bG$ defined over $K$ containing
a maximal $K$-split torus
and a set $\Delta$ of simple roots in
the root system $\Phi(\bG,\mathbf T)$.  Denote by $2\rho$ the sum of all positive
roots in $\Phi(\bG, \mathbf T)$, and by $\lambda_\iota$ the maximal weight of
$\iota$. Define $u_\alpha, m_\alpha\in \n $, $\alpha\in \Delta$, by
 $$2\rho=\sum_{\alpha\in \Delta} u_\alpha\alpha\quad\text{and}\quad
 \lambda_\iota=\sum_{\alpha\in \Delta} m_\alpha \alpha .$$
The fact that $m_\alpha\in \n$ follows since $\bG$ is of adjoint type.
Consider the twisted action of the Galois group $\Gamma_K:=\operatorname{Gal}(\bar K/K)$ on
$\Delta$ (for instance, if the $K$-form of $\bG$ is inner,
this action is just trivial).
Then 
 \begin{equation}\label{abdef}a_\iota=\max_{\alpha\in \Delta}\frac{u_\alpha+1}{m_\alpha}\quad\text{and}\quad
b_\iota =\# \{\Gamma_K.\alpha: \frac{u_\alpha+1}{m_\alpha}=a_\iota \}.\end{equation}

Note that the exponent $a_\iota$ is independent of the field $K$,
and $b_\iota$ depends only on the quasi-split $K$-form of $\bG$.
Therefore, by passing to a finite field extension containing the splitting
field of $\bG$, $b_\iota$
also becomes independent of $K$.



\vs \noindent{\bf Remark: } 
When $\bG$ is almost $K$-simple or, more generally, when
 $\He_\iota$ is the product of height functions of the $K$-simple
factors of $\bG$, we can improve the rate of convergence in Theorem \ref{mt}: for some $\delta>0$,
\begin{align*}\label{eq:rate}
N(\He_\iota, T)= c\cdot T^{a_\iota} P(\log T)\cdot (1+O(T^{-\delta}))
\end{align*}
where $P(x)$ is a monic polynomial of degree $b_\iota-1$.

\subsection{Distribution of rational points}
For each $v\in R$, denote by $X_{\iota, v}$ the closure of
$\bar \iota (\bG(K_v))$ in $\mathbb{P}(\M_N(K_v))$, and
consider the compact space $X_\iota :=\prod_{v\in R} X_{\iota, v}$.
In section \ref{arith},
we construct a probability measure $\mu_\iota$ on $X_\iota$
which describes the asymptotic distribution of rational points in $\bG(K)$
 in $X_\iota$ with respect to the height
$\He_\iota$.
To keep the introduction concise, we give the definition of
 $\mu_\iota$ only when $\iota$ is saturated.
 A representation $\iota:\bG\to \GL_N$ is called {\it
 saturated} if the set
\begin{equation*}
\{\alpha\in \Delta: \frac{u_\alpha+1}{m_\alpha}=a_\iota
\}\end{equation*} is not contained in the root system of a proper
normal $K$-subgroup of $\bG$. 
In particular, if $\bG$ is almost $K$-simple, any representation of $\bG$ is saturated.

 Let $\tau$ denote the Haar measure on $\bga$
such that $\tau (\bG(K)\ba \bga)=1$.
Denote by $\Lambda$ the set of all automorphic characters of $\bga$
(cf. section \ref{defla}) and by $W_\iota$ the maximal compact subgroup
of the group $\bG(\A_f)$ of finite adeles, under which $\He_\iota$ is bi-invariant
(see Definition \ref{wi}).
Then the following is a positive real number (see Propositions \ref{gammas} and \ref{ri} (3), noting
$r_\iota=\gamma_{W_\iota}(e)$ in the notation therein):
\begin{equation}\label{riota} r_{\iota}:=
\sum_{\chi\in \Lambda}\lim_{s\to a_\iota^+}(s-a_\iota)^{b_\iota}
\int_{\bga} \He_{\iota}(g)^{-s} \chi(g)\; d\tau(g) .\end{equation}

For $\iota$ saturated,
the probability
measure $\mu_{\iota}$
 on $ X_\iota$ is the unique measure satisfying that
for any $\psi\in  C(X_\iota)$ invariant under a co-finite subgroup of $W_\iota$,
\begin{equation}\label{intro:muiota} \mu_\iota(\psi)
=r_{\iota}^{-1}\cdot \sum_{\chi\in \Lambda}
\lim_{s\to a_\iota^+}(s-a_\iota)^{b_\iota}
\int_{\bga} \He_{\iota}(g)^{-s} \chi(g)\;\psi(g)\; d\tau(g)
 \end{equation}
(see Theorem \ref{smuiota-g}).
We refer to \eqref{arith-measure} for the definition of $\mu_\iota$ for a general
$\iota$:
\begin{Thm}\label{e}
For any $\psi \in C(X_{\iota})$,
$$\lim_{T\to \infty} \frac{1}{N(\He_{\iota}, T)} \sum_{g\in \bG(K): \He_\iota(g)<T} \psi (g)= \int_{X_\iota} \psi \, d\mu_\iota .$$
\end{Thm}


\begin{Rmk} \rm 
\begin{enumerate}
\item For $\iota$ saturated, the measure $\mu_\iota$
coincides with the measure $\tilde\mu_\iota$
which describes the distribution of height balls in $\bga$ (see Proposition \ref{measure}).
 \item 
Although the projection $\mu_{\iota,S}$ of $\mu_{\iota}$ to
$X_{\iota, S}=\prod_{v\in S}X_{\iota, v}
$ is always equivalent to
 a Haar measure on $\bG_S=\prod_{v\in S} \bG(K_v)$ (Proposition \ref{accc}),
 it is $\bG_S$-invariant, only when the height
$\He_{\iota,S}=\prod_{v\in S}\He_{v}\circ \iota$ is $\bG_S$-invariant.
\item
The space $X_{\iota,S}$ is a compactification of $\bG_S$ which
  is an analog of the Satake compactification defined for real groups
  (see, for example, \cite{BJ}).
Theorem \ref{e} implies that the rational points $\bG(K)$
do not escape to the boundary $X_{\iota,S}-\bG_S$.
 It is interesting to compare this result with the
distribution of the integral points $\bG(\z)$ of bounded height in the Satake
compactification of $\bG(\br)$ where the limiting distribution
is supported on the boundary (see \cite{GOS}
and \cite{Mau} for more details).
\end{enumerate}
\end{Rmk}

\subsection{Counting and volume heuristic}
To explain our strategy in counting $K$-rational points of $\bG$, we
first recall the analogous results in counting integral
 points in a simple real algebraic group.
Let $G\subset \GL_N$ be a connected non-compact simple real algebraic group and
$\Gamma$ be a lattice in $G$, i.e., a discrete subgroup of
finite co-volume. Fixing a norm $\|\cdot\|$ on $\M_N(\br)$,
 set $B_T:=\{g\in G: \|g\|\le T\}$.
By Duke-Rudnick-Sarnak \cite{DRS} and Eskin-McMullen \cite{EM}
independently, it is well known  that
\begin{equation}\label{integral}\#\G\cap B_T\sim \int_{B_T} dg \quad\text{as $T\to \infty$}, \end{equation}
where $dg$ is the Haar measure on $G$ such that
$\int_{\Gamma\backslash G}dg=1$. 

Coming back to the question of counting rational points $\bG(K)$, 
we recall that $\bG(K)$ is a lattice in the adele
group $\bga$ when embedded diagonally
and that the height function $\He_\iota=\prod_{v\in R}
\He_v \circ \iota$ on $\bG(K)$ extends to $\bga$.

If we set $$B_T:=\{g\in \bga: \He_\iota(g)<T\} ,$$ 
then $B_T$ is a
relatively compact subset of $\bga$ (Lemma \ref{ff3}) and we have
the equality
 $$N(\He_\iota, T)=\# \bG(K)\cap B_T .$$
In view of (\ref{integral}), one naturally asks
whether the following holds:
\begin{equation}\label{rational}\#\bG(K) \cap B_T\sim 
\tau(B_T) \quad\text{as $T\to \infty$}. \end{equation}

It turns out that the group $\bga$ is too big for (\ref{rational})
to hold in general, due to the presence of non-trivial automorphic
 characters of $\bga$.
For a compact open subgroup $W_f$ of $\bG(\A_f)$,
denote by $\Lambda^{W_f}\subset \Lambda $ 
the set of all  
$W_f$-invariant characters in $\Lambda$.
We set 
\begin{equation*}
G_{W_f}:=\op{ker}(\Lambda^{W_f})=\cap
\{\op{ker}\chi\subset \bga: \chi\in \Lambda^{W_f}\}.
\end{equation*}
The subgroup $G_{W_f}$ is a normal subgroup of
$\bga$ with finite index (see Lemma \ref{ker}), and hence $\bG(K)$ is a lattice in $G_{W_f}$.
Denote by $\tau_{W_f}$ the Haar measure on $G_{W_f}$
normalized so that $\tau_{W_f}(\bG(K)\backslash G_{W_f})=1$.


\begin{Thm}\label{eqin} Assume that $\iota:\bG\to \GL_N$ 
is saturated.
Then for any compact open subgroup $W_f$ of $\bG(\A_f)$
under which $\He_\iota$ is bi-invariant,
$$\# \bG(K)\cap B_T\sim_T \tau_{W_f}(G_{W_f}\cap B_T).$$
\end{Thm}

 We remark that one cannot in general replace $G_{W_f}$ by
 $\bga$ (see example \ref{sec:ex}), and Theorem \ref{eqin}
does not hold for $\iota$ non-saturated.




As in the proof of Eskin-McMullen of (\ref{integral}),
our key ingredient in proving Theorem \ref{eqin}
is the mixing theorem on
$L^2(\bG(K)\backslash G_{W_f})$.

\subsection{Adelic mixing}
Let $L^2_{00}(\bG(K)\backslash \bga)$ denote the orthogonal
complement 
in $L^2(\bG(K)\backslash \bG(\A))$ to the
 direct sum of all automorphic characters.
  In the case when $\bG$
is simply connected, $L^2_{00}(\bG(K)\backslash \bga)$ coincides with the
orthogonal complement
$L^2_0(\bG(K)\backslash \bG(\A))$
to the constant functions.

 Set $\bG_\infty:=\prod_{v\in R_\infty} \bG(K_v)$
 where $R_\infty$ is the subset of $R$ of all archimedean
 valuations.
\begin{Thm}[Automorphic bound for $\bG$] \label{mmix} Let $\bG$ be
a connected absolutely almost simple $K$-group. Let $U_\infty$ be a
maximal compact subgroup of $\bG_\infty$ and $W_f$ a compact open
subgroup of $\bG(\A_f)$.
Then there exist $c_{W_f}>0$ and $r_0=r_0(\bG_\infty)>0$ such that
for any $U_\infty$-finite and $W_f$-invariant functions $\psi_1,
\psi_2\in L^2_{00}(\bG(K)\backslash \bG(\A))$,
$$|\langle \psi_1, g .\psi_2
 \rangle  |\le c_{W_f} \cdot(\op{dim}\langle U_\infty \psi_1\rangle \cdot
\op{dim}\langle U_\infty \psi_2 \rangle)^{r_0}\cdot \tilde
\xi_\bG(g)\cdot\|\psi_1\|_2\cdot \|\psi_2\|_2 \quad\text{for all $g\in
\bG(\A)$}.$$  Here, $\tilde \xi_\bG:\bga \to (0,1]$
 is an explicitly constructed proper function
which is $L^p$-integrable for some $p=p(\bG)<\infty$. (see Def.
\ref{tx}).
\end{Thm}

Using the restriction of scalars functor,
we extend this theorem to connected
almost $K$-simple adjoint (simply connected) groups (Theorem \ref{am}).
The above bounds on matrix coefficients can
also be extended to smooth functions
in certain Sobolev spaces (see Theorem \ref{smooth}).
We also mention a paper of Guilloux \cite{Gu} where an application of Theorem
\ref{mmix} was discussed in local-global principle problems.


\begin{Cor}[Adelic Mixing] \label{mix}
 Let $\bG$ be a connected absolutely simple $K$-group,
or a connected almost $K$-simple adjoint (simply connected) $K$-group.
 Then for any
 $\psi_1, \psi_2\in L^2_{00}(\bG(K)\backslash \bga) $,
$$\langle \psi_1, g .\psi_2 \rangle \to 0 $$
as $g\in \bG(\A)$ tends to infinity.
\end{Cor}
Any $W_f$-invariant function
in $ L^2(\bG(K)\backslash \bhi)$ orthogonal to constants
belongs to  $L^2_{00}(\bG(K)\backslash \bga)$ (see Lemma \ref{wif}).
Hence Corollary \ref{mix} implies
 that if $\psi_1$ and $ \psi_2$ are $W_f$-invariant functions in $
L^2(\bG(K)\backslash \bhi)$, then as $g\to \infty$,
$$\int_{\bG(K)\backslash \bhi}
\psi_1(x)\psi_2(xg)\; d\tau_{W_f}(x)
\to \int \psi_1 \, d\tau_{W_f}  \cdot \int \psi_2\, d\tau_{W_f} . $$



\vs\vs
For each $v\in R$, denote $\hat \bG_v^{\op{Aut}}\subset \hat \bG_v$
 the automorphic dual of $\bG(K_v)$, i.e., the
 subset of
unitary dual of $\bG (K_v)$ consisting of representations which are
  weakly contained in the representations appearing as $\bG(K_v)$
components of $L^2(\bG(K)\backslash \bga)^{O_f}$ for some compact
open subgroup $O_f$ of $\bG(\A_f)$.
The proof of Theorem \ref{mmix} goes roughly as follows: 
 if $\tilde
\xi_v$ is a uniform bound for the matrix coefficients of infinite
dimensional representations in $\hat \bG_v^{\op{Aut}}$, $\tilde
\xi_\bG$ is defined to be the
 product $\prod_{v\in R}\tilde {\xi_v}.$ This can be made precise using the language of
 direct integral of a representation (cf. proof of Theorem \ref{k}).
For those $v\in R$ such that the $K_v$-rank of $\bG$ is at least
$2$, the uniform bounds, say $\xi_v$, of matrix coefficients of {\it
all} infinite dimensional unitary representations of $\bG(K_v)$ were
obtained by Oh \cite{Oh1}. For these cases, one can simply take $\tilde
\xi_v=\xi_v$. In particular, if $K$-rank of $\bG$ is at least $2$
and $\bG (K_v)^+$ denotes the closed subgroup of
$\bG (K_v)$ generated by all unipotent elements in $\bG (K_v)$,
we have $\tilde \xi_\bG=\prod_{v\in R}\xi_v$ and $\tilde \xi_\bG$
works as a uniform bound for all unitary representations of $\bga$
without $\bG(K_v)^+$-invariant vectors for each $v\in R$ (see Theorem
\ref{k} for a precise statement). Moreover $\tilde \xi_\bG$ is
fairly sharp in these cases. For instance, one can show that $\tilde
\xi_\bG$ is optimal for $\bG=\SL_n$ $(n\ge 3)$, or $\Sp_{2n}$ $(n\ge
2)$ by
 \cite[5.4]{COU}.

When there is $v\in R$ with $K_v$-rank of $\bG$ one,
 finding an automorphic bound $\tilde \xi_v$ is essentially carried out by Clozel \cite{Cl1}.
In particular, several deep theorems in automorphic theory were used such as
the Gelbart-Jacquet bound \cite{GJ} toward Ramanujan conjecture,
 the results of Burger-Sarnak \cite{BS} and Clozel-Ullmo \cite{CU} on
lifting automorphic bounds,
the base changes by Rogawski \cite{Ro} and Clozel \cite{Cl2}, and Jacquet-Langlands correspondence \cite{JL}.


\subsection{Organization of the paper}
In section \ref{s-notation} we list some notations and preliminaries
which will be used  throughout the paper.
In section \ref{s-mixing}, we discuss adelic mixing and prove
 Theorem \ref{mmix}. We also extend
a theorem of Clozel-Oh-Ullmo on
the equidistribution of Hecke points \cite{COU} in this section
as an application of adelic mixing. In section \ref{s-volume},
we deduce the volume asymptotics of height balls from
the results in \cite{STT2} and
construct in subsection \ref{muiota} the probability measure $\tilde\mu_\iota$ on $X_\iota$ which
describes the asymptotic distribution of height balls in $\bga$.
The main theorem in section \ref{sec:proof} is Theorem \ref{smuiota}
on the equidistribution of rational points with respect to $\tilde\mu_\iota$
for the case when $\iota$ is saturated.
Theorems \ref{e} (saturated cases) and \ref{eqin} follow from this theorem.
In section \ref{arith}, we prove Theorem \ref{e} for general cases
(Theorem \ref{uns}) as well as Theorem \ref{mt} (Theorem \ref{t:eqgeneral}).
In section \ref{manin}, we restate our main theorems in the
context of Manin's and Peyre's conjectures. 
\subsection{Acknowledgments}
We would like to thank Wee Teck Gan, Emmanuel Peyre and Yehuda Shalom
 for helpful conversations. We thank Ramin Takloo-Bighash
for the useful comments on our preliminary version of this paper.
We are also deeply grateful to the referee for many detailed comments
 on the submitted version.

\section{Notations and Preliminaries}\label{s-notation}
We set up some notations which will be used throughout the
paper.
Let $K$ be a number field and $\bG$ a connected semisimple 
group defined over $K$. We denote by $R_K$, or simply by $R$,
the set of all normalized absolute values on $K$.
We keep the same notation $R, R_f, K_v$
as in the introduction.

\subsection{} Let $\Cal O$ denote the ring of integers of $K$
and $\Cal O_v$ the valuation ring of $K_v$.
Set $R_\infty=R-R_f$. For $v\in R_f$,
 let $q_v$ denote the order of the residue field of $\Cal O_v$.
We choose an absolute value $|\cdot|_v$  on $K_v$ normalized so that
the absolute value of a uniformizer of $\Cal O_v$ is given by
$q_v^{-1}$. 
Denote by $\A$ the adele ring over $K$ and by $\bga$ the adele group
associated to $\bG$.

Denote by $\bG(\A_f)$ (resp. $\bG_\infty$)  the subgroup of finite
(resp. infinite) adeles, i.e., $((g_v)_v)\in \bga$ with $g_v=e$ for
all $v\in R_\infty$ (resp. for all $v\in R_f$). Then
$$\bga=\bG_\infty\times \bG(\A_f) .$$

\subsection{}\label{smodel}
We fix a smooth model $\mathcal G$
 of $\bG$ over $\Cal O[k^{-1}]$ for some non-zero
$k\in \z$.
There exists a finite subset $S_0\subset R_f$ such that
for any $v\in R_f-S_0$, $\bG$ is unramified over $K_v$ and
$\mathcal G(\Cal O_v)$ is a hyperspecial compact subgroup  (cf. \cite{Ti2}).
We set $U_v=\mathcal G(\Cal O_v)$ for each $v\in R_f-S_0$.
Then for each $v\in R_f-S_0$, one has the group $A_v$ of
$K_v$-rational points of a maximal $K_v$-split torus of $\bG$ so that
the following Cartan decomposition holds:
\begin{equation}\label{eq:cartan}
\bG(K_v)=U_vA_v^+ U_v
\end{equation}
where $A_v^+$ is a closed positive Weyl chamber of $A_v$.
More precisely,
one can choose a system $\Phi_v^+$ of positive roots in
the set $\Phi_v=\Phi(\bG(K_v),
A_v)$ of all non-multipliable roots of $\bG(K_v)$ relative to $A_v$
so that
$$A_v^+=\{a\in A_v: \alpha(a)\ge 1 \;\;\text{for each $\alpha\in \Phi_v^+$}\}\quad\text{ for
$v$ archimedean} $$
$$A_v^+=\{a\in A_v: |\alpha(a)|_v\in q_v^{\n}\;\;\text{for each $\alpha\in \Phi_v^+$}\} 
\quad\text{otherwise}.$$

For $v\in S_0\cup R_\infty$, there exists
 a good maximal compact subgroup $U_v$ (cf. \cite[2.1]{Oh1} for definition)
of $\bG(K_v)$ such that
$$\bG(K_v)=U_vA_v^+\Omega_v U_v$$
where
 $\Omega_v$ is a finite subset in the centralizer of $A_v$ in $\bG(K_v)$.

In particular for any $g\in \bG(K_v)$, there exist unique $a_v\in A_v^+$ and
$d_v\in \Omega_v$ such that $g\in U_v a_v d_v U_v$.
 For $v\in R_\infty$, any maximal compact subgroup
of $\bG(K_v)$ is a good maximal compact subgroup and
$\Omega_v=\{e\}$.

\subsection{}\label{taumu} For a finite subset $S$ of $R$, let $\bG^S$ denote the subgroup of
$\bga$ consisting of $(g_v)$, with $g_v=e$ for all $v\in S$,
 and set $\bG_S:=\prod_{v\in S}\bG(K_v)$.
 Note
that $\bga=\bG_S\bG^S$. 
For each $v\in R$,
let $\tau_v$, or $dg_v$, denote a Haar measure on $\bG(K_v)$
such that $\tau_v(U_v)=1$ whenever $v\in R_f$.
Then the collection $\{\tau_v: v\in R\}$ defines a Haar measure, say $\tau$,
 on $\bga$ (cf. \cite[3.5]{PR}). We will assume
that $\tau(\bG(K)\ba \bga)=1$. This is possible by replacing $\tau_v$,
$v\in R_\infty$ with a suitable multiple of it, since
 $\bG(K)$ is a lattice in $\bga$
(cf. \cite[Theorem 5.5]{PR}).

 We denote by $\tau_S$
the product measure $\prod_{v\in S}\tau_v$ on $\bG_S$ and
by $\tau^S$ the Haar measure on
 $\bG^S$ for which the triple $(\tau, \tau_S, \tau^S)$ are
compatible with each other, i.e., $\tau=\tau_S\times \tau^S$ locally.

\subsection{}\label{defla} 
An automorphic character of $\bga$ is a continuous homomorphism
from $\bga$ to the unit circle $\{z\in \mathbb C: z\bar z=1\}$
which contains $\bG(K)$ in its kernel. 
Each automorphic character $\chi$ of $\bga$ can be considered
as a function on the quotient $\bG(K)\backslash \bG(\A)$, and
since $\tau(\bG(K)\ba \bga)=1$,
 $\chi$ belongs to $L^2(\bG(K)\backslash \bG(\A))$
and $\|\chi\|_2=1$.
Let $\Lambda$ denote the set of 
all automorphic characters of $\bga$.
Note that any two distinct elements of $\Lambda$ are orthogonal to each other.
We then have an orthogonal decomposition
$$L^2(\bG(K)\backslash \bG(\A))=L^2_{00}(\bG(K)\backslash \bga)\oplus
 \hat \bigoplus_{\chi\in \Lambda } \bc \chi .$$
where $\hat \bigoplus_{\chi\in \Lambda } \bc \chi$
is the
closure of the direct sum of $\bc \chi$'s, and
$L^2_{00}(\bG(K)\backslash \bG(\A))$ denotes its orthogonal
complement in $L^2(\bG(K)\backslash \bG(\A))$.

 If $\bG$ is simply
connected, it follows from the strong approximation property that
the only automorphic character of $\bga$ is the trivial one and hence that
$L^2_{00}(\bG(K)\backslash \bga)$ is the orthogonal
complement to the space of
constant functions (cf. Lemma \ref{zero}).

\subsection{}\label{gwf} 
For any compact open subgroup $W_f$ of $\bG(\A_f)$,
we denote by $\Lambda^{W_f}\subset \Lambda $ 
the set of all  
$W_f$-invariant characters in $\Lambda$, i.e.,
$ \Lambda^{W_f}=\{\chi\in\Lambda: \chi (w)=1\quad \text{for all $w\in W_f$}\}$.
We set 
\begin{equation}\label{eq:ghi}
G_{W_f}:=\op{ker}(\Lambda^{W_f})=\cap
\{\op{ker}\chi\subset \bga: \chi\in \Lambda^{W_f}\}.
\end{equation}
The subgroup $G_{W_f}$ is a normal subgroup of
$\bga$ with finite index (see Lemma \ref{ker}), and hence $\bG(K)$ is a lattice in $G_{W_f}$.
Denote by $\tau_{W_f}$ the Haar measure on $G_{W_f}$
normalized so that $\tau_{W_f}(\bG(K)\backslash G_{W_f})=1$.

\subsection{}\label{hed} For a group $\bG$ of adjoint type,
let $\iota:\bG\to \GL_N$ be a faithful
representation defined over $K$.
 We give a definition of a height
function $\He_\iota$ on $\bga$ associated to $\iota$ which is slightly more
general than those considered in the introduction.
 It is this class of the functions for which we prove
our main theorems.

\begin{Def}\label{hedh} A height function $\He_\iota:\bga \to \br^+$
is defined by the product $\prod_{v\in R} \He_{\iota,v}$
where $\He_{\iota, v}$ is a function on $\bG(K_v)$ for $v\in R$ satisfying the following:
\begin{enumerate}\item there exists a finite subset $S\subset R$ such that
$$\He_{\iota, v}(g)=\max_{ij} |\iota(g)_{ij}|_v\quad\text{ for all $v\in R-S$};$$
\item for $v\in S$, there exists $C>0$
such that
$$ C^{-1} \cdot \max_{ij} |\iota(g)_{ij}|_v
\le \He_{\iota,v}(g) \le C\cdot \max_{ij} |\iota(g)_{ij}|_v;$$
\item for any $v\in S\cap R_\infty$,
there exists $b>0$ such that for any small $\e>0$,
$$(1-b\cdot \e) \He_{\iota, v}(x)\le \He_{\iota, v}(gxh)\le (1+b\cdot \e) \He_{\iota, v}(x)$$
for any $x\in \bG(K_v)$ and any $g, h$ in the $\e$-neighborhood of
$e$ in $\bG(K_v)$ with respect to a Riemannian metric;
\item for any $v\in S\cap R_f$,
 $\He_{\iota, v}$ is bi-invariant under a compact open
subgroup of $\bG(K_v)$.
\end{enumerate}
\end{Def}
Note by (1) that for $(g_v)\in \bga$,
since $g_v\in \Cal G(\Cal
O_v)$ for almost all $v\in R_f$,
$\He_{\iota, v}(\iota(g_v))=1$ for almost all
$v$, and hence
 the product $\prod_{v\in R}\He_{\iota, v}(g_v)$ converges.

Note also that the class of height functions defined above does not depend on
the choice of a basis of $K^N$.

We will need the following observation on heights:
\begin{Lem}\label{l:prod} Suppose that $\iota$ has a unique maximal weight.
Let $\bG_1$ and $\bG_2$ be connected
normal algebraic $K$-subgroups of $\bG$ with $\bG=\bG_1\bG_2$
and $\bG_1\cap\bG_2=\{e\}$.
There exists $\kappa>1$ such that for any $g_1\in
\bG_1(\A)$ and $g_2\in \bG_2(\A)$,
\begin{align*}
\kappa^{-1}\cdot {\He}_\iota(g_1){\He}_\iota
(g_2)\le {\He}_\iota(g_1g_2)&\le \kappa\cdot {\He}_\iota(g_1){\He}_\iota
(g_2).
\end{align*}
\end{Lem}
\begin{proof}
Let $\lambda_\iota$ denote the highest weight of $\iota$.
Then there exists a finite subset $S\subset R$ such that for any
 $v\in R-S$,
$$\bG(K_v)=U_v A^+_v U_v\quad\text{and}\quad
\He_{v}(\iota(g))=|\lambda_\iota (a)|_v\quad\hbox{for $g=u_1 au_2\in \bG(K_v)$}
$$
where $U_v$ and $A_v^+$ are defined as in \eqref{smodel}.
 In particular, it
follows that for each $v\in R-S$, and for any $g_1\in \bG_1(K_v)$ and
$g_2\in \bG_2(K_v)$,
\begin{align*}
{\He}_{v}(\iota(g_1g_2))={\He}_{v}(\iota(g_1)){\He}_{v}(\iota(g_2)).
\end{align*}
On the other hand, for $v\in S$, $\He_{\iota,v}$ is equivalent to
$\lambda_\iota$ in the sense that there exists $\kappa_v>1$ such that
$$
\kappa_v^{-1}\cdot |\lambda_\iota(a)|_v\le \He_{\iota,v}(g) \le \kappa_v\cdot
|\lambda_\iota(a)|_v  \quad\hbox{for $g=u_1 a d u_2\in U_v A_v^+ \Omega_v U_v=
\bG(K_v)$.}
$$
This implies the lemma.
\end{proof}

For $T>0$, set
$$B_T:=\{g\in \bga:\He_\iota(g)< T\} .$$

\begin{Lem}\label{ff3} \begin{enumerate}
\item We have
\begin{equation}\label{eq_delta0}
\delta_0:=\inf _{g\in \bga}\He_\iota(g)>0.
\end{equation}
\item For each $T>0$,
  $B_T$ is a relatively compact subset of $\bga$.
In other words, the height function
$\He_\iota:\bga\to [\delta_0, \infty)$ is proper.
\end{enumerate}
\end{Lem}
\begin{proof}
By Definition \ref{hed}, there exists a finite subset $S$ such that for
all $v\in R-S$, $\He_v(\iota(g))\ge 1$ for any $g\in \bG(K_v)$. Let
$0<\delta\le 1$ be such that $\He_v(\iota(g))\ge \delta$ for $v\in
S$ and $\delta_1=\delta^{\# S}$. Then $\He_\iota(g)\ge \delta_1$
for all $g\in \bga$. Hence $\delta_0\ge \delta_1>0$.

 Note that
 $$B_T\subset \bga\cap \prod_v \{g_v\in \bG(K_v): \He_v(\iota(g_v)) 
\le \delta^{-1} T\} .$$
Since for almost all $v\in R_f$, $\He_v(\iota(g_v))\ge q_v$ whenever
$g_v\notin \mathcal G(\Cal O_v)$,
it follows that for some finite subset $S_1\subset R$, we have
 $$B_T\subset \{(g_v)_v\in\bga:
 \He_v(\iota(g_v)) \le \delta_0^{-1} T\text{ for $v\in S_1$},\;\;
g_v\in \mathcal G (\Cal O_v)\text{ otherwise}\}.$$
Since the set
$\{g_v\in \bG(K_v): \He_v(\iota(g_v))\le b\}$ is compact for any $b>0$, 
it follows that $B_T$ is a
relatively compact subset of $\bga$.
\end{proof}

\begin{Def}\label{wi} 
For a height function $\He_\iota$ of $\bga$,
define $W_{\He_\iota}$, or simply $W_\iota$, to be
$$W_\iota=\{w\in \bG(\A_f): \He_\iota(wg)=\He_\iota(gw)=\He_\iota(g)\quad
\text{ for all }g\in \bga\}.$$
\end{Def} It is easy to check that $W_\iota$ is a subgroup of
$\bG(\A_f)$, and is compact by the above lemma. Hence
 $W_\iota$ is the maximal compact open subgroup
of $\bG(\A_f)$ under which $\He_\iota$ is bi-invariant.

\section{Adelic Mixing}\label{s-mixing}
\subsection{Definition and properties of $\xi_\bG$}\label{dpxi}
 Let $\bG$ be a connected semisimple
algebraic group defined over a number field $K$.
Let $\Cal T$ denote the set of $v\in R$
such that $\bG(K_v)$ is compact, that is,
$U_v=\bG(K_v)$.
 It is well known that $\Cal T$ is a finite set.

Denote by $\Phi_v^+$ the system of positive roots in
the set of all non-multipliable roots of $\bG(K_v)$ relative to $A_v^+$ and choose  a maximal
strongly orthogonal system $\Cal S_v$ in
$\Phi_v^+$ in the sense of \cite{Oh1}
(where an explicit construction is also given).
For $v\in R-\Cal T$ and $K_v\ne \bc$, define
 the bi-$U_v$-invariant function $\xi_v=\xi_{\bG(K_v)}$ on $\bG(K_v)$ (cf. \cite{Oh1}):
  for  each $g=kadk'\in U_vA_v^+\Omega_v U_v$,
$$\xi _{v}(g)=\prod_{\alpha\in \Cal S_v} \Xi_{\PGL_2(K_v)}
 \begin{pmatrix} \alpha(a) & 0\\ 0 & 1\end{pmatrix}$$
where $\Xi_{\PGL_2(K_v)}$ is the Harish-Chandra function of $\PGL_2(K_v)$.
If $K_v=\bc$, set
$$\xi _{v}(g)=\prod_{\alpha\in \Cal S_v} \Xi_{\PGL_2(\bc)}
 \begin{pmatrix} \alpha(a) & 0\\ 0 & 1\end{pmatrix} ^{n_\alpha}$$
where $n_\alpha=1/2$ if $\alpha$ is a long root
of $\bG$, when $\bG$ is locally isomorphic to $\Sp_{2n}(\bc)$, and $n_\alpha=1$
for all  other cases.
We set $\xi_v=1$ for $v\in \Cal T$.

Since $0< \xi_v(g_v)\le 1$ for all $v\in R$ and  $\xi_v(g_v)=1$ for almost all $v$,
the following function $\xi_\bG$ is well defined:
\begin{Def}\label{xi} Define the function $\xi_\bG:\bga \to (0, 1]$ by
\begin{equation*}
\xi_\bG(g)=\prod_{v\in R}\xi_v(g_v) \quad\text{for $g=(g_v)_v\in\bga$}.\end{equation*}
\end{Def}

Set
\begin{equation}
\label{du}
U_f=\prod_{v\in R_f} U_v, \quad U_\infty=\prod_{v\in R_\infty} U_v,
\quad \text{and }U= U_f\times U_\infty .\end{equation}
Note also that $\xi_{\bG}$ is bi-$U$-invariant.

For $v\in R-\Cal T$,
 we set 
  $$\eta_v(kadk'):=\prod_{\alpha\in \Cal S_v} |\alpha(a)|_v $$
 where $kadk'\in U_vA_v^+\Omega_v U_v$ for
all $v$ with $K_v\ne \bc$.
As in the case of the definition of $\xi_v$,
if $K_v=\bc$ and for $kak'\in U_vA^+_vU_v$, we set 
$$\eta_v(kak')=
\prod_{\alpha\in \Cal S_v}|\alpha^{n_\alpha}(a)|_v$$
 with the same $n_\alpha$ defined as before.
If $v\in\Cal T$, we set $\eta_v=1$.
\begin{Lem}\label{xiin}
For any $\e>0$, there is a constant $C_{\e}>0$  such that for any $g=(g_v)_v\in \bga$,
\begin{equation}\label{eta}
 \prod_{v\in R} \eta_v(g_v)^{-1/2} \le \xi_\bG(g) \le C_{ \e} \cdot \prod_{v\in R}
\eta_v(g_v)^{-1/2+\e} . \end{equation}
In particular,
$$\xi_\bG (g)\to 0\quad\text{as $g\to \infty$ in $\bga$}.$$
\end{Lem}
\begin{proof}
For $v\in R-\Cal T$, it follows from the explicit formula for $\Xi_v$
(cf. \cite[3.8]{Oh1})
 that
for any $\e>0$, there is a constant $C_{v,\e}>0$ such that for any $g_v\in \bG(K_v)$,
\begin{equation*}
 \eta_v(g_v)^{-1/2} \le \xi_v(g_v) \le C_{v, \e} \cdot \eta_v(g_v)^{-1/2+\e} . \end{equation*}
Moreover one can take $C_{v, \e}=1$ for almost all $v$.
This implies (\ref{eta}).

To see the second claim, first note that
for any $g\in \bga$, \begin{equation}\label{inf}
\xi_\bG(g)\le\xi_v (g_{v})\le C_{v,\e} \cdot  \eta_v(g_{v})^{-1/2+\e} .\end{equation}
 Now suppose on the contrary that there exists a sequence $\{g_i\in \bga\}$
such that $g_i\to \infty$ and
$\xi_\bG(g_i)\nrightarrow 0 .$ Then by passing to a subsequence
we may assume either that there is a place $v\in R$ such that
$g_{i,v}\to \infty$ in $\bG(K_v)$ or
that there exists a sequence
 $\{v_i\in R_f-S_0\}$ such that $g_{i,{v_i}}\notin U_{v_i}$
and  $q_{v_i}\to \infty$.
If $g_{i,v}\to \infty$ as $i\to \infty$, then
 $ |\eta_v(g_{i,v})|\to \infty$ as $i\to \infty$ and hence
$\xi_\bG(g_i)\to 0$ by (\ref{eta}). Therefore the first case cannot
happen.

In the second case,
note that since $g_{i, v_i}\notin U_{v_i}$ and $\Omega_{v_i}=\{e\}$,
we have $\eta_v(g_{i,v_i})\ge q_{v_i}$ for each $i$.
Hence by (\ref{inf}) for all $i$ big enough,
$$\xi_\bG(g_i)\le  C_{v_i, \e} \cdot q_{v_i}^{-1/2+\e} \le
q_{v_i}^{-1/2+\e} .$$
This gives a contradiction since $q_{v_i}\to \infty$.
\end{proof}

\begin{Lem}\label{xhi} Let $\iota:\bG\to \op{GL_N}$ be a faithful
representation defined over $K$ with a unique maximal weight
and $\He_\iota$ be a height function on $\bga$ associated to $\iota$.
Then there exist $m\in \n-\{0\}$ and $C>0$ such that
 $$\xi_\bG (g)\le C\cdot \He_\iota^{-1/m}(g)\quad\text{for any $g\in \bga$}.$$
\end{Lem}
\begin{proof} Let $\chi$ denote the highest weight of $\iota$.
 Let $l\in \n-\{0\}$ be such
 that $\chi|_{A_v^+}\le l\cdot \log_{q_v} \eta_v$ for each $v\in R$.
 Here ${q_v}=e$ if $v\in R_\infty$.
 Without loss of generality, we may assume
$$\He_v(\iota(a_v))=q_v^{\chi( a_v)}\quad \text{ for each $a_v\in A_v^+$ and $v\in R$}.$$

Since
$\eta_v(a_v)=q_v^{\log_{q_v} |\eta_v(a_v)|}$ for $a_v\in A_v^+$,
we have for each $v\in R$,
$$\eta_v(a_v)^{-l}\le \prod_v \He_v ^{-1} (\iota(a_v))\quad\text{for
$a_v\in A_v^+$}.$$

By the continuity of $\He_v$ and the Cartan decomposition 
$\bG(K_v)=U_vA_v^+\Omega_v U_v$, there exists $r_v\ge 1$ such that
$$r_v^{-1} \He_v(\iota(a_v))\le \He_v(\iota(g_v))\le r_v \He_v(\iota(a_v))$$
for $g_v=k_v a_v d_v k_v'\in U_vA_v^+\Omega_v U_v$.
Since $\He_\iota$ is invariant under a compact open
subgroup of $\bG(\A_f)$ and hence $\He_v\circ \iota$ is invariant under
$U_v$ for almost all $v\in R_f$, we can take $r_v=1$ for almost all
 $v\in R_f$.

Therefore if $g=(g_v)\in \bga$ with $g_v=k_v a_v d_v k_v'\in U_vA_v^+\Omega_v U_v$,
\begin{align*}\prod_v \eta_v(g_v)^{-l}&=\prod_v \eta_v(a_v)^{-l} \le
 \prod_v \He_v ^{-1} (\iota(a_v)) \\& \le
r_0 \prod_v \He_v ^{-1} (\iota(g_v))
 \le r_0\cdot \He_\iota^{-1}(g) \end{align*}
where $r_0=\prod_v r_v<\infty$.

 Hence using Lemma \ref{xiin}, 
there exists $c_1>0$ such that for any $g=(g_v)\in \bga$,
 $$\xi_\bG ^{4l}(g) \le c_1\cdot \prod_v
\eta_v(g_v)^{-l}\le c_1\cdot r_0 \cdot
\He_\iota^{-1}(g). $$
This proves the claim.
\end{proof}

Theorem 7.1 in \cite{STT2} shows that for ay representation
$\iota$ of $\bG$ over $K$ with a unique maximal weight,
 the height zeta function $$\Cal Z(s):=\int_{\bga}
\He_\iota(g)^{-s}d\tau(g)$$ converges for $\Re (s)>a_\iota$
where $a_\iota$ is defined as in (\ref{abdef}).
In particular $\He_\iota^{-1}$ belongs to $L^{p}(\bga)$ for any $p>a_\iota$.
Hence as a corollary of \cite[Theorem 7.1]{STT2} using Lemma \ref{xhi}, we 
obtain the following:
\begin{Cor}\label{lp} There exists $0<p=p(\bG)<\infty $ such that
$\xi_\bG\in L^p(\bga)$.
\end{Cor}

\subsection{Restriction of scalars functor}
We prove certain funtorial properties of
$\xi_{\bG} $ for the restriction of scalars functor introduced by Weil
\cite{We}, which
we will need later to reduce the dicussion on
general almost $K$-simple groups
to that on absolutely simple $K$-groups.
We refer to \cite[I. 3.1.4]{Mar2}
and \cite[Ch. 6]{BT} for the properties of
the restrictions of scalar functor $R_{K/k}$ used in the following discussion.
Suppose that $k$ is a finite extension field of $K$ and
$\bG'$ is a connected semisimple $k$-group.
Then $\bG=R_{k/K}\bG'$ is a connected
semisimple $K$-group. 

Denote by $\{1=\sigma_{1}, \cdots, \sigma_{d}\}$ the set of
all distinct embeddings of
$k$ into the algebraic closure of $K$.
Then there is a $K$-morphism
$\mu:\bG\to \bG'$
such that the map \begin{equation}\label{rs}
\mu^\circ=({}^{\sigma_1}\mu, \cdots, {}^{\sigma_{d}} \mu):
\bG\to {}^{\sigma_1}\bG'\times \cdots\times {}^{\sigma_d} \bG'\end{equation}
is a $K$-isomorphism.
If $R^\circ_{k/K}$ denotes the inverse map to
$\mu|_{\bG(K)}$, then
 $R^\circ_{k/K}:\bG'(k)\to \bG(K)$ is a group isomorphism. 

For each $v\in R_K$,
denote by $I_v$ the set of
all valuations of $k$ extending $v$.
Then there is a natural $K_v$-isomorphism
$f_v:\bG\to \prod_{w\in I_v} R_{k_{w}/K_v} \bG '$
and the isomorphisms $f_v^{-1}\circ R_{k_w/K_v}^\circ :
\prod_{w\in I_v}\bG'(k_w)\to \bG(K_v)$,
 $v\in R_K$, $w\in I_v$, induce a topological group isomorphism, say, $j$,
of the adele group $\bG'(\A_k)$ to $\bG(\A_K)$.


\begin{Lem}\label{resc} Let $\e>0$. Then there are constants
$C_\e\ge 1 $ such that for any $g\in \bG'(\A_k)$,
$$C_\e^{-1} \cdot \xi_{\bG }(j(g))^{1+\e} \le
 \xi_{\bG'}(g)\le  C_\e \cdot \xi_{\bG }(j(g))^{1-\e}.$$
\end{Lem}
\begin{proof}

Fix $v\in R_K$.
The set $I_v$ parametrizes the set, say, of
 all distinct embeddings $\sigma$ of
$k$ into $\bar K_v$ which are non-conjugate over $K_v$,
in the way that $w\in I_v$ corresponds to
$\sigma$ with $k_w=\sigma(k)K_v$.
For each embedding $w\in I_v$,
we denote by $J_w$ the set
of all embeddings $\tau$
of $k$ into $\bar K_v$ such that $k_w=\tau(k)K_v$.
Fix $w\in I_v$.
Let $A'_{w}$ be
 the group of $k_{w}$-points of a maximal $k_{w}$-split $A'$
 torus of $\bG'$,
$\Phi(\bG'(k_{w}),A'_{w})$ the set of non-multipliable roots,
and $\mathcal S_{w}\subset \Phi(\bG'(k_{w}),A'_{w})$
be a maximal strongly orthogonal system used in the definition
of $\xi_{\bG'}$.
Then $\prod_{w\in I_v} (R_{k_{w}/K_v}A')$ is a maximal $K_v$-torus
of $\bG$ and its maximal $K_v$-split subtorus
is a maximal $K_v$-split torus of $\bG$ \cite[Ch. 6]{BT}.
For each $w\in I_v$, let $B(w)$ denote the group of $K_v$-points
of the maximal $K_v$-split torus of $R_{k_{w}/K_v}A'$ and
set $A_v=\prod_{w\in I_v} B(w)$. 
We can identify each $ B(w)$ with 
$\{({}^\tau a_w)_{\tau\in J_w}: a_w\in A'(K_v)\}$,
 and
$\Psi_w:=\Phi((R_{k_{w}/K_v}\bG')(K_v), B(w))$ with 
$\{ \alpha_w  : 
\alpha\in \Phi( \bG'(k_{w}), A'_{w})\}$
where $\alpha_w(({}^\tau a_w)_{\tau\in J_w})=
\prod_{\tau\in J_w} {}^\tau \alpha(a_w)$.
Hence
$\{ \alpha_w: \alpha\in \mathcal S_{w}\}$
is a maximal strongly orthogonal system of
$\Psi_w$.

On the other hand,
 for any $\alpha\in \mathcal S_w$,
$$|\alpha (a_w)|_w
= | \prod_{\tau\in J_w} {}^\tau \alpha(a_w) |_v .$$
Since $\prod_{w\in I_v} \mathcal S_w$ is a maximal strongly
orthogonal system for $\Phi( \bG(K_v), A_v)$,
if $a=(({}^\tau a_w)_{\tau\in J_w} )_{w\in I_v}\in A_v$,
$$\eta_{\bG(K_v)}(a)=\prod_{(\alpha_w)_{w\in I_v}
\in \prod_{w\in I_v}\mathcal S_w}
|\alpha_w( ({}^\tau a_w)_{\tau\in J_w} )|_v
=\prod_{w\in I_v}\prod_{\alpha_w\in \mathcal S_w}
 |\alpha(a_w)|_w
=\prod_{w\in I_v}\eta_{\bG'(k_w)}(a_w) .$$
Hence for all $g\in \bG'(\A_k)$,
$$\prod_{v\in R_K}\eta_{\bG(K_v)} (j(g))
=\prod_{w\in R_k}\eta_{\bG'(k_w)}(g) .$$
By Lemma \ref{xiin}, this proves the claim.
\end{proof}

\subsection{Uniform bound for matrix coefficients of  $\bG(\A)$}
Let $W_f\subset \bG(\A_f)$ be a compact open subgroup.
Write $W_v=W_f\cap \bG(K_v)$ for each $v\in R$.
Then $W_v= U_v$ for almost all $v\in R_f$.
For each $v\in R_f$, by \cite{Be},
there exists $d_{W_v} <\infty$ such that
for any irreducible unitary representation $\rho$ of $\bG(K_v)$,
the dimension of $W_v$-invariant vectors of $\rho$
is at most $d_{W_v}$.
Moreover $d_{W_v}=1$ for almost all $v\in R$ by \cite[3.3.3]{Ti2}
and \cite[Corollary 1]{Fl}.
Hence the following number is well-defined:
 $$d_{W_f}:=\prod_{v\in R_f} d_{W_v} <\infty .$$

The notation $\bG (K_v)^+$ denotes the normal subgroup
 of $\bG(K_v)$ generated by all unipotent subgroups of $\bG(K_v)$.

\begin{Thm}\label{k} Let $\bG$ be a connected absolutely
 almost simple $K$-group
with $K$-rank at least $2$. Let $W_f$ be a compact open subgroup
of $\bG(\A_f)$.
Let $\pi$ be any unitary representation of $\bga$ without
no non-trivial $\bG (K_v)^+$-invariant vector for every $v\in R$.
Then for any $U_\infty$-finite and $W_f$-invariant unit vectors $x$ and $y$,
\begin{equation}\label{q}
|\langle \pi(g)x, y\rangle| \le  d_0\cdot c_{W_f}
\cdot (\op{dim}\langle U_\infty x\rangle \cdot \op{dim}\langle U_\infty
 y\rangle)^{(r+1)/2}  \cdot \xi_\bG(g)\quad
\text{ for all $g\in \bga$}\end{equation}
 where $c_{W_f}:=  d_{W_f} \cdot \prod_v[U_v:U_v\cap W_v]\cdot (\max_{d\in \Omega_v}
[U_v: dU_vd^{-1}])$ and $d_0, r \ge 1$ depend only on
 $\bG$.
 Moreover if $\bG(K_v)\ncong\op{Sp}_{2n}(\mathbb C)$ locally for any $v\in R_\infty$, $d_0=1$ and $r=1$.

If $\bG$
is a connected almost $K$-simple adjoint (simply connected) $K$-group
with $K$-rank at least $2$, and $\pi$
has no non-trivial $\bL(K_v)^+$-invariant vectors for any connected $K_v$-normal
subgroup $\bL$ of $\bG$,
then
\eqref{q} holds with $\xi_{ \bG}$ replaced by
$ \xi_{\bG}^{1-\e}$ for any $\e>0$.
\end{Thm}

As a corollary,
we obtain the
 adelic version of Howe-Moore theorem \cite{HM} on the 
vanishing of matrix coefficients:
\begin{Cor}\label{kk} Let $\bG$ and $\pi$ be as in Theorem \ref{k}.
Then for any vectors $x$ and $y$,
$$\langle \pi(g)x, y\rangle \to 0\quad \text{as $g\to \infty$ in $\bga$}.$$
\end{Cor}
\begin{proof}
It is easy to see that it suffices to prove the claim for
a dense subset of vectors in the Hilbert space $V$ associated to
$\pi$. Considering the restriction $\tilde\pi:=\pi|_U$ to $U$,
 $V$ decomposes into a direct sum of irreducible unitary representations of
the compact group $U$ each of which is finite dimensional  
by Peter-Weyl theorem. Hence
the set $V_0$ of $U$-finite vectors is dense in $V$.
Now if $x, y\in V_0$, then $x$ and $y$ are invariant under
a finite index subgroup $W_f$ of $U_f$.
Hence by applying the above theorem,
we obtain that for some constant $c_0>0$,
$$|\langle \pi(g)x, y\rangle| \le c_0 \cdot \xi_{\bG}(g)\quad\text{for all $g\in
\bga$} .$$
Since $\xi_{\bG}(g)\to 0$ as $g\to\infty$,
this implies the claim.
\end{proof}

The proof of Theorem \ref{k} is based on theorems in \cite{Oh1}. More precisely, recall:
\begin{Thm} \cite[Theorem 1.1-2]{Oh1}\label{o}
Suppose that
the $K_v$-rank of $\bG$ is at least $2$.
Let $\pi_v$ be a unitary representation of $\bG(K_v)$
without any non-trivial $\bG(K_v)^+$-invariant vectors.
Then for any $U_v$-finite unit vectors $x$ and $y$,
$$|\langle \pi_v(g)x, y\rangle|\le d_v \cdot c_v \cdot
 (\op{dim}\langle U_v x\rangle \cdot \op{dim}\langle U_v y\rangle)^{r_v/2} \cdot \xi_v(g)
\quad\text{ for any $g\in \bG(K_v)$}$$
 where $c_v=\max_{d\in \Omega_v} [U_v: dU_vd^{-1}]$ and $d_v, r_v \ge 1 $ depend only on
 $\bG(K_v)$. Moreover whenever $\bG(K_v)\ncong\op{Sp}_{2n}(\mathbb C)$ locally, $d_v=1$ and $r_v=1$.
\end{Thm}
In the case when $\bG(K_v)\cong\op{Sp}_{2n}(\mathbb C)$ locally, the
above theorem was stated only for $U_v$-invariant vectors in \cite{Oh1}.
However if we replace Proposition 2.7 in \cite{Oh1}
by the remark following it, the same proof works for the above claim.

\begin{proof}[\bf Proof of Theorem \ref{k}]
We first assume that $\bG$ is absolutely almost simple.
For $g=(g_v)_v\in \bga$, choose a finite subset $S_g$ of places
containing $$\{v\in R_f: g_v\notin U_v \}\cup R_\infty  .$$

Note that for $v\in R-S_g$,
we have $g_v\in U_v$ and hence $\xi_v(g_v)=1$.
Therefore for $g=(g_v)_v\in \bga$,
$$\xi(g)=\prod_{v\in S_g}\xi_v(g_v) .$$
Let $G_g=\prod_{v\in S_g}\bG(K_v)$ and $W_g=\prod_{v\in S_g\cap
R_f}W_v$. As a $G_g$ representation,
$\pi$ has a Hilbert integral decomposition:
$$\pi=\int_{z\in Z_g} \oplus ^{m_z}\rho_z \;d\nu(z)$$
where $Z_g$ is the unitary dual of $G_g$ and $\rho_z$
is irreducible, $m_z$ is a multiplicity for each $z\in Z_g$ and $\nu$ is a measure on $Z_g$
(see \cite{Di} or \cite[Section 2.3]{Zi}).
We may assume that for all $z$, $\rho_z$ has no $\bG(K_v)^+$-invariant vector (see
\cite[Prop. 2.3.2]{Zi}).

If we write $\Cal L_z=\oplus ^{m_z} \rho_z$, $x=\int x_z d\nu(z)$
and $y=\int y_z d\nu(z)$ with $$x_z=\sum_{i=1}^{m_z} x_{zi}\quad
\text{and} \quad   y_z=\sum_{i=1}^{m_z} y_{zi}\in \Cal L_z,$$ we have
$$\langle x, y \rangle=\int_{Z_g} \sum_{i=1}^{m_z}
 \langle x_{zi} ,y_{zi} \rangle\, d\nu(z). $$
It follows from the definition of
a Hilbert direct integral that
 $$\op{dim}\langle U_{\infty}x_{zi}\rangle \le \op{dim}\langle U_\infty x_{z}\rangle \le \op{dim}\langle U_\infty x\rangle ,$$
$x_{zi}$ is $W_g$-invariant
for almost all $z$ and all $i$, and similarly for $y$.
Without loss of generality, we assume the above holds for all $z$.
We claim that
\begin{align}\label{zi}|\langle \rho_{z}(g) x_{zi}, y_{zi}\rangle  | \le
 c_{W_f}  \cdot d_0 \cdot   \xi_{\bG}(g)\cdot
(\op{dim}\langle U_\infty x\rangle \cdot \op{dim}\langle U_\infty y
\rangle)^{(r+1)/2}\cdot \|x_{zi}\|\cdot \|y_{zi}\|
\end{align}
 where $r=\max_v {r_v}$ and $c_0=
d_{W_f}\prod_v( c_v \cdot [U_v:U_v\cap W_v])<\infty, d_0=\prod_v d_v
<\infty $ with $c_v,d_v, r_v$
  as in Theorem \ref{o}.
By \cite{Be}, we may write
$\rho_{z}=\otimes_{v\in S_g}\rho_{z(v)}$ where $\rho_{z(v)}$ is an
irreducible representation of $\bG(K_v)$ without no non-trivial
$\bG(K_v)^+$-invariant vectors. Since the finite linear combinations
of pure tensor vectors are dense, it suffices to prove (\ref{zi})
assuming $x_{zi}$ and $y_{zi}$ are finite sums of pure tensors. Hence
we can write
$$ x_{zi}=\sum_{j}
\bigotimes_{v\in S_g} x_{z{ij}(v)}\;;\; y_{zi}=\sum_{k}\bigotimes_{v\in S_g} y_{z{ik}(v)}$$
where for each $v\in S_g$, $x_{zij(v)}$ (resp. $y_{zik(v)}$) are
mutually orthogonal and the number of summands for $x_{zi}$ (resp.
$y_{zi}$) is at most $\op{dim}\langle U_\infty x\rangle\cdot d
_{W_f}$ (resp. $\op{dim}\langle U_\infty y\rangle\cdot d_{W_f}$).
Hence by Cauchy-Schwarz inequality, for $x_{zij}= \prod_{v\in S_g}
x_{zij(v)}$ and $y_{zij}= \prod_{v\in S_g} y_{zij(v)}$
$$\sum_j \|x_{zij}\|
\le (\op{dim}\langle U_\infty x\rangle\cdot d_{W_f})^{1/2}\|x_{zi}\|;
\;\text{and}\;\;\sum_k \|y_{zik}\|
 \le (\op{dim}\langle U_\infty y \rangle\cdot d_{W_f})^{1/2}\|y_{zi}\| .$$

Since for $v\in R_f$
$$\op{dim} \langle U_v x\rangle \le [U_v:W_v\cap U_v]\quad\text{and}\quad
\op{dim} \langle U_v y\rangle \le [U_v:W_v\cap U_v],$$
 by Theorem \ref{o}, we have for $c_0=\prod_v c_v$,
\begin{align}\label{iine}
&|\langle \rho_{z}(g) x_{zi}, y_{zi}\rangle  |
\le \sum_{j, k}  \prod_{v\in S_g} |\langle \rho_{z(v)}(g_v)x_{z{ij}(v)}, y_{z{ik}(v)}\rangle | \\
 &
 \le c_0 \cdot d_0 \cdot \prod_{v\in S_g}  \xi_v(g_v) \cdot
(\op{dim}\langle U_\infty x\rangle \cdot \op{dim} \langle U_\infty
y\rangle)^{r/2} (\prod_{v\in R_f}[U_v:W_v\cap U_v])\cdot
 \left(\sum_{j, k}\|x_{z{ij}}\|\cdot \|y_{z{ik}}\|\right)\notag \\
 &\le c_0 \cdot  d_0 \cdot\xi_{\bG}(g) \cdot
(\op{dim}\langle U_\infty x\rangle \cdot \op{dim}\langle U_\infty
y \rangle)^{(r+1)/2}\cdot
 \prod_{v\in R_f}[U_v:W_v\cap U_v])\cdot d_{W_f} \left(\|x_{zi}\|\cdot \|y_{zi}\|\right)\notag\\
&= c_{W_f} \cdot  d_0 \cdot\xi_{\bG}(g) \cdot (\op{dim}\langle
U_\infty x\rangle \cdot \op{dim}\langle U_\infty y
\rangle)^{(r+1)/2}\cdot \left(\|x_{zi}\|\cdot \|y_{zi}\|\right)
 \notag\end{align}
 proving
 (\ref{zi}).
Therefore again by Cauchy-Schwarz inequality,
\begin{align}\label{ine}
|\langle (\oplus^{m_z}\rho_z)(g)(x_z), y_z\rangle|&\le
\sum_i|\langle \rho_{z}(g) x_{zi}, y_{zi}\rangle  |
\\&\le c_{W_f} \cdot  d_0 \cdot \xi_\bG(g)\cdot (\op{dim}\langle U_\infty x\rangle \cdot \op{dim}
\langle U_\infty y\rangle)^{(r+1)/2} \cdot  \|x_z\|\cdot
\|y_z\|\notag .
 \end{align}
By integrating over $Z_g$,
we obtain (\ref{q}).

Since $\bG$ is adjoint (resp. simply connected),
there exists a finite separable extension $k$ of $K$ and a connected
adjoint (resp. simply connected) absolutely almost simple $k$-group $\bG'$
such that
$\bG=R_{k/K}\bG'$ by \cite[3.1.2]{Ti1}. 
Then the topological isomorphism $j^{-1}: \bG(\A_K)\to
\bG'(\A_k)$ described prior to Lemma \ref{resc} 
 maps $\bG(K_v)$ and $\bG(K_v)^+$ to 
$\prod_{w\in I_v}\bG'(k_w)$ and $\prod_{w\in I_v}\bG'(k_w)^+$
respectively.
Let $\pi$ be a representation on $\bga$ satisfying
the hypothesis. Then 
the representation, say, $\pi'$, on $\bG'(\A_k)$ induced by
$\pi$ has no $\bG'(k_w)^+$-invariant vectors for each $w\in R_k$.
Since the $k$-rank of $\bG'$ is equal to the $K$-rank
of $\bG$ and hence is at least $2$,
by the assertion already proved for
$\bG'(\A_k)$, we
deduce that
for any $U_\infty$-finite and $W_f$-invariant $x, y$,
\begin{align*}
|\langle \pi(g) x,  y \rangle|&=|\langle \pi'(j^{-1}(g)) x,  y \rangle|
\\ & \le d_0\cdot
 c_{W_f}
 \cdot (\op{dim}\langle U_\infty x\rangle \cdot \op{dim}\langle U_\infty y\rangle)^{(r+1)/2}
 \cdot \xi_{\bG '}(j^{-1}(g)) \quad\text{for all $g\in \bG(\A)$} .\end{align*}
By \eqref{eqtil}, we may replace $\xi_{\bG'}(j^{-1}(g))$ by
$\xi_{\bG}(g)^{1-\e}$, finishing the proof.

\end{proof}

\vs\vs

\subsection{Automorphic bound for $\bG(\A)$}
If $\bG$ has $K$-rank at most one, the analogue of Theorem \ref{k}
does not hold in general. However if we look at those infinite
dimensional representations occurring in
 $L^2(\bG(K)\backslash \bG(\A))$,
we still obtain a similar upper bound.


We first state the following conjecture:
\begin{Con}\label{ram} Let $\bG$ be a connected absolutely almost
simple $K$-group. Let $W_f$ be a compact open subgroup
of $\bG(\A_f)$.
Then for any $U_\infty\times W_f$-invariant unit vectors
$f, h\in L^2_{00}(\bG(K)\backslash \bga)$,
$$|\langle f, g. h\rangle |\le c_{W_f}\cdot  \xi_\bG (g)
\quad \text{for all $g\in \bga$}$$
where $c_{W_f}>0$ is a constant depending only on $\bG$ and $W_f$.
\end{Con}
The above holds for groups of $K$-rank at least $2$ by Theorem \ref{k}.
For $\bG=\PGL_2$, Conjecture \ref{ram}
 is essentially equivalent to the Ramanujan conjecture.
We will prove a weaker statement of Conjecture \ref{ram} where the function $\xi_\bG$ is replaced by
a function $\tilde \xi_\bG$ with slower decay such that $\xi_\bG \le \tilde\xi_\bG\le \xi_\bG ^{1/2}$.

\begin{Def}\label{tx}
Let $\bG$ be a connected almost $K$-simple group.
For each $v\in R$, write $\bG$ as an almost direct product 
$\bG_{v}^1\bG_{v}^2$ where $\bG_{v}^1$ is
the maximal semisimple normal $K_v$-subgroup of $\bG$
such that every simple normal $K_v$-subgroup of $\bG_{v}^1$ has $K_v$-rank one.
 Note that $\bG_{v}^2$ is then the maximal semisimple
normal $K_v$-subgroup of $\bG$ without any $K_v$-normal subgroup
of rank zero or one.
We define a function $\tilde \xi_{\bG}:\bga\to (0,1]$ by
 $$\tilde\xi_{\bG}:=\prod_{v\in R}
\left(\xi_{\bG_{v}^1(K_v)}^{1/2}\cdot \xi_{\bG_{v}^2(K_v)}\right).$$
 \end{Def}
If $\bG$ is absolutely almost simple
and $R_1:=\{v\in R:\op{rank}_{K_v}(\bG)=1\}$.
 then
$$\tilde\xi_{\bG}=
\prod_{v\in R_1}\xi_{\bG(K_v)}^{1/2}\cdot
\prod_{v\in R-R_1}\xi_{\bG(K_v)} .$$
If $\bG=R_{K/k}\bG'$ for some finite extension field
$k$ and for a connected absolutely almost simple $k$-group $\bG'$,
then, for any $w\in R_k$ extending $v\in R_K$,
the $k_{w}$-rank of
a connected simple $k_{w}$-subgroup $\bH'$ of $\bG'$
is equal to the $K_{v}$-rank of the $K_v$-subgroup
$R_{k_{w}/K_v}\bH'$ of $\bG$ (cf. \cite[Ch I. 3.1]{Mar2}).
Using this,
the proof of Lemma \ref{resc} also shows that
there is $C_\e>1$ such that for any $g\in \bG'(\A_k)$
\begin{equation}\label{eqtil}
C_\e^{-1} \cdot \tilde \xi_{\bG }(j(g))^{1+\e} \le
 \tilde \xi_{\bG'}(g)\le  C_\e \cdot \tilde \xi_{\bG }(j(g))^{1-\e}.
\end{equation}
where $j$ denotes the topological group isomorphism
 from $\bG'(\A_k)$ to $\bga$ described prior to Lemma \ref{resc}.

\begin{Thm}[Automorphic bounds]\label{am}
Let $\bG$ be a connected absolutely almost simple $K$-group,
For a compact open subgroup
$W_f$ of $\bG(\A_f)$, there exist $r=r(\bG)\ge 1$ and $c_{W_f}>0$
such that for any $U_\infty$-finite
and $W_f$-invariant unit vectors $x, y\in L^2_{00}(\bG(K)\backslash \bga)$,
$$|\langle x, g .y \rangle|\le c_{W_f}
 \cdot (\op{dim}\langle U_\infty x\rangle \cdot \op{dim}\langle U_\infty y\rangle)^{(r+1)/2}
 \cdot \tilde \xi_\bG(g) \quad\text{for all $g\in \bG(\A)$} .$$
One can take
 $r=1$ provided for any $v\in R$, $\bG(K_v)$
has no subgroup locally isomorphic
to $\op{Sp}_{2n}(\mathbb C)$ $(n\ge 2)$ locally.

If $\bG$
is a connected almost $K$-simple adjoint (simply connected) $K$-group,
then
the above inequality holds with $\tilde\xi_{ \bG}$ replaced by
$\tilde \xi_{\bG}^{1-\e}$ for any $\e>0$.

\end{Thm}

Recall that for unitary representations $\rho_1$ and $\rho_2$ of $\bG(K_v)$,
$\rho_1$ is said to be weakly contained in $\rho_2$ if every diagonal matrix coefficients of $\rho_1$ can be approximated
uniformly on compact subsets by convex combinations of diagonal matrix coefficients of $\rho_2$.
For each $v\in R$, denote by $\hat \bG_v$
 the unitary dual of $\bG(K_v)$ and by $\hat \bG_v^{\op{Aut}}\subset \hat \bG_v$
 the automorphic dual of $\bG(K_v)$ as defined in the introduction.
The following theorem was first obtained by Burger and Sarnak
for $v$ archimedean
\cite[Theorem 1.1]{BS} and generalized by Clozel and Ullmo to all $v$
\cite[Theorem 1.4]{CU}.
\begin{Thm}\label{cu} Let $\bG$ be a connected absolutely
almost simple $K$-group.
Let $\bH\subset \bG$ be a connected semisimple $K$-subgroup.
Then for any $v\in R$ and for any $\rho_v \in \hat\bG_v^{\op{Aut}}$,
any irreducible representation of $\bH(K_v)$ weakly contained in 
$\rho_v|_{\bH(K_v)}$is contained in $ \hat \bH_v^{\op{Aut}}$.
\end{Thm}

\begin{Lem}\label{zero} For any $v\in R$ such that $\bG(K_v)$ is non-compact,
$L^2_{00}(\bG(K)\backslash \bga)$ has no non-zero
 $\bG(K_v)^+$-invariant function.
\end{Lem}
\begin{proof} 
Let $\Cal L_v$ denote the set of $f\in L^2_{00}(\bG(K)\backslash \bga)$ fixed by $\bG(K_v)^+$.
We need to show that $\Cal L_v=\{0\}$.
Let $\bG^{\{v\}}$ denote the subgroup
of $\bga$ consisting of elements whose $v$-component is trivial.
Consider the family of continuous functions $f\in C_c(\bG^{\{v\}})$ of the form
$f=\prod_{w\in R-\{v\}} f_w $ where each $f_w$ is a continuous function of $\bG(K_w)$
  such that $f_w|_{U_w}=1$ for almost all $w$.
By considering the convolutions with these functions,
we obtain a dense family of
 the continuous functions belonging to $\Cal L_v$.
Hence
it suffices to show that any continuous function $f\in \Cal L_v$ is trivial.
Let $f\in \Cal L_v$ be continuous.
Let $\tilde {\bG}$ be the simply connected cover of
$\bG$ and denote by $\op{pr}:\tilde \bG\to \bG$ the covering map.
Consider
the projection map $$\tilde \bG(K)\backslash \tilde{\bG}(\A)\to
 \bG(K)\backslash \bga .$$
Let $\tilde f$ be the pull back of $f$.
Since the image of $\tilde \bG(K_v)$ 
under the map $\op{pr}$ is
 $\bG(K_v)^+$,
the function $\tilde f$ is left $\tilde \bG (K)$-invariant
 and right $ \tilde \bG (K_v)$-invariant. On the other hand,
 the strong approximation property implies that $\tilde \bG (K)\tilde \bG (K_v)$ is dense
in $\tilde{\bG}(\A)$ (cf. \cite[Theorem 7.12]{PR}).
Therefore $\tilde f$ is constant, and hence $f$ factors through 
the image of $\tilde{\bG} (\A)$ in $\bga$. Since $\bga/\tilde \bG (\A)$ is abelian,
$L^2( \bG(K) \tilde \bG (\A) \ba \bga)$ is a sum of automorphic characters of $\bga$.
Since $f$ is orthogonal to $\Lambda$,
$f=0$.
\end{proof}

\begin{proof}[\bf Proof of Theorem \ref{am}]
We first treat the case when $\bG$ is
absolutely almost simple.
  The case when $K$-rank is at least $2$ follows from Theorem \ref{k} and Lemma \ref{zero}.
Suppose first that $\bG$ has $K$-rank one.
By \cite[Theorem 3.4]{COU}, for $v\in R_1$,
any infinite dimensional  $\rho_v\in \hat {\mathbf G}_v^{\op{Aut}}$,
and $U_v$-finite vectors $x_v,y_v$, 
\begin{equation}\label{ccu}
|\langle \rho_v(g)(x_v), y_v\rangle|
\le  c_v\cdot \xi_v (g)^{1/2}\cdot (\op{dim}\langle U_vx_v \rangle \cdot \op{dim}
\langle U_vy_v\rangle)^{1/2} \end{equation}
for any $g\in \bG(K_v)$.
Combining this with Theorem \ref{o},
we can derive the desired bound by the same argument as in
the proof of Theorem \ref{k}.

Now suppose $\bG$ is $K$-anisotropic. 
For $v\in R_1$, we claim \eqref{ccu} holds.
In \cite{Cl1}, it is analyzed what kind of
$\rho_v$ occurs in this situation, and this
is the main case which was not known before Clozel's work.
We give a brief summary. If $R_1\ne \emptyset$, it
follows from the classification theorem by Tits \cite{Ti1} that $\bG$ is
of Dynkin type $\Cal A$. \cite[Theorem 1.1]{Cl1} says that
there exists a $K$-embedding of $K$-subgroup $\bH$ of type $\Cal A$
such that $\bH$ has $K_v$-rank one whenever $v\in R_1$. Let $v\in
R_1$. Then up to isogeny, one has either that $\bH=\PGL_1(D)$ for a
quaternion algebra $D$ over $K$ and $\bH= \PGL_2$  over $K_v$, or
$\bH=\operatorname{PGU}(D, *)$ for a division algebra $D$ of prime
degree $d$ over a
quadratic extension $k$ of $K$ with a second kind involution $*$,
and $\bH= \op{PGU}(n-1, 1)$ over $K_v$ (with $n\ge 3$).
 In the former case one uses the
Jacquet-Langlands correspondence \cite{JL} to transfer the
Gelbart-Jacquet automorphic bound of $\op{PGL_2}$ to $\bH(K_v)$ via
Theorem \ref{cu}. In the second case which is hardest,
 by the base changes obtained by
Rogawski \cite{Ro} and Clozel \cite{Cl2},
 we can use the bound of $\PGL_n(F_w)$ given by
Theorem \ref{o} to get a bound
for $\bH(K_v)$ where $w$ is a place of $k$ lying above $v$
and $k_w$ is a quadratic extension of $K_v$.
This proves the claim.
 Combining with Theorem \ref{o}
for those places $v\in R-(R_1\cup \Cal T)$ as in the proof of
Theorem \ref{k}, we obtain the desired bound.

Now for the case when $\bG$ is adjoint (resp. simply connected)
almost $K$-simple,
 the same argument
used in the proof of Theorem \ref{k} applies, since (in the notation therein)
 the topological isomorphism $j: \bG(\A_K)\to
\bG'(\A_k)$ 
 induces
 an equivariant isometry between $L^2_{00}(\bG(K)\ba \bga)$ and
$L^2_{00}(\bG'(k)\ba \bG'(\A_k))$.
\end{proof}

\subsection{From $U_\infty$-finite vectors to smooth vectors}
In Theorem \ref{am}, we can relax $U_\infty$-finite conditions
to smooth conditions provided we replace the $L^2$-norms by 
$L^2$-Sobolev norms. For a precise formulation,
let $X_1, \cdots, X_{m}$ be an orthonormal
basis of the Lie algebra $\op{Lie}(U_\infty)$ with respect to an $\op{Ad}$-invariant scalar product.
Then the elliptic operator
\begin{equation}\label{elliptic}
\Cal D:=1-\sum_{i=1}^m X_i^2
\end{equation} lies in the center of the universal enveloping algebra of
 $\op{Lie}(U_\infty)$.
We say a function $f$ on $\bG(K)\backslash \bG(\A)$ is smooth
if $f$ is invariant under some
compact open subgroup of $\bG(\A_f)$ and smooth for the action of $\bG_\infty$.
\begin{Thm} \label{smooth} Let $\bG$ be a connected almost absolutely
simple $K$-group. 
and $W_f$ be a compact open subgroup of
$\bG(\A_f)$. Then there exist an explicit $l\in \n$ and
$c_{W_f}>0$ such that
for any $W_f$-invariant smooth functions 
$\varphi, \psi\in L^2_{00}(\bG(K)\backslash \bga)$ with
$\|\Cal D^l (\varphi)\|<\infty$
and $\|\Cal D^l (\psi)\|<\infty$,
$$|\langle \varphi, g. \psi\rangle |\le c_{W_f}\cdot  \tilde \xi_\bG (g) \cdot
\|\Cal D^l (\varphi)\| \cdot\| \Cal D^l(\psi)\|\quad\text{for all $g\in \bga$}.$$
If $\bG$
is a connected almost $K$-simple adjoint (simply connected) $K$-group,
then
the above inequality holds with $\tilde\xi_{ \bG}$ replaced by
$\tilde \xi_{\bG}^{1-\e}$ for any $\e>0$.
\end{Thm}
\begin{proof} Deducing this from Theorem \ref{am} is quite standard in view of
the results of Harish-Chandra explained in \cite[Ch 4]{Wa}.
We give a sketch of the proof.
Denote by $\pi$ the representation $ L^2_{00}(\bG(K)\backslash \bga)$.
Then $\pi=\oplus_{\nu\in \hat{U}_\infty} \pi_\nu$ where $\pi_\nu$ is
the $\nu$-isotypic component of $\pi$ and $\Cal D$ acts as a scalar, say, $c_\nu$
on each $\pi_\nu$.
We write $\varphi=\sum_{\nu\in \hat{U}_\infty}\varphi_\nu$ and
$\psi=\sum_{\nu\in \hat{U}_\infty}\psi_\nu$.
One has $\|\varphi_\nu\|=c_\nu^{-l}\|\Cal D^l \varphi_\nu\|$ and similarly for $\psi$.
Then $$|\langle \varphi, g. \psi\rangle |\le \sum_{(\nu_1, \nu_2)\in \hat{ U}_\infty\times \hat{U}_\infty}
|\langle \varphi_{\nu_1}, g.\psi_{\nu_2}\rangle |$$
Using Theorem \ref{am}, we then obtain
 \begin{align*}
 &|\langle \varphi, g. \psi \rangle |\\& \le c_{W_f}\cdot \tilde{\xi}_\bG (g)\left(
 \sum_{\nu\in \hat U_\infty}\|\varphi_{\nu} \|\op{dim}\langle
 U_\infty \varphi_{\nu} \rangle ^{(r+1)/2} \right)\left(\sum_{\nu\in \hat U_\infty}\| \psi_{\nu}\|\op{dim}\langle
 U_\infty \psi_{\nu} \rangle ^{(r+1)/2} \right) \\
 & \le
 c_{W_f}\cdot
\tilde{\xi}_{\bG}(g)\cdot \|\Cal D^l (\varphi)\|\cdot \|\Cal D^l( \psi)\|
 \cdot
 \sum_{\nu\in \hat U_\infty} {c_{\nu}} ^{-2l}
 \dim ( \nu )^{r+1} \end{align*}
Now if $l\in \n$ is sufficiently large, then $\sum_{\nu}c_{\nu} ^{-2l}
 \dim ( \nu )^{r+1} <\infty$ \cite[Lemma 4.4.2.3]{Wa}.
 This proves the claim.
\end{proof}

\subsection{From $K$-simple groups to semisimple groups}
 If $\bG$ is a connected semisimple $K$-group, we say that a
sequence $\{g_i\in \bga\}$ tends to infinity strongly if for any
non-trivial connected simple normal $K$-subgroup $\bH$ of $\bG$,
$p_\bH (g_i)$ tends to $\infty$ as $i\to \infty$, where
$p_\bH:\bG(\A)\to\bG(\A)/\bH(\A)$ denotes the canonical projection.

\begin{Thm}[Mixing for $L^2(\bG(K)\backslash \bG(\A))$] \label{smix}
 Let $\bG$ be a product of connected almost $K$-simple $K$-groups.
 Then for any
 $\varphi, \psi \in L^2_{00}(\bG(K)\backslash \bga) $,
$$\langle \varphi, g .\psi \rangle \to 0 $$
as $g\in \bG(\A)$ tends to infinity strongly.
\end{Thm}
\begin{proof} Write $\bG=\bG_1\times \cdots \times \bG_m$ where
each $\bG_i$ is a connected absolutely almost simple $K$-group.
By Theorem \ref{am} and Peter-Weyl theorem (cf. Corollary \ref{kk}),
for each $1\le i\le m$,
and  for any $\varphi_i, \psi_i\in L^2_{00}(\bG_i(K)\backslash \bG_i(\A)) $,
\begin{equation}\label{peter}\langle \varphi_i, g_i .\psi_i\rangle\to 0\end{equation}
as $g_i\to \infty$ in $\bG_i(\A)$.

Consider $\otimes_{i=1}^m L^2(\bG_i(K)\backslash \bG_i(\A))$ as a
subset of  $L^2(\bG(K)\backslash \bG(\A))$. The finite sums of
the functions of the form $\psi=\otimes_{i=1}^m  \psi_i\in
L^2(\bG(K)\backslash \bga)$, $\psi_j\in L^2(\bG_j(K)\backslash
\bG_j(\A))$, such that for at least one $j$, $\psi_j\in
L^2_{00}(\bG_j(K)\backslash \bG_j(\A))$ form a dense subset of the
space $L^2_{00}(\bG(K)\backslash \bG(\A))$. Hence it suffices to
prove the claim for $\varphi= \otimes_{i=1}^m \varphi_i$ and $\psi=\otimes_{i=1}^m
\psi_i$ of such type.  Suppose $\psi_j\in L^2_{00}(\bG_j(K)\backslash
\bG_j(\A))$ for some $1\le j\le m$.  If $g=(g_1, \cdots, g_m)$ with
$g_i\in \bG_i(\A)$, then
$$ |\langle \varphi, g .\psi \rangle |=\prod_{i=1}^m |
\langle \varphi_i, g_i .\psi_i \rangle|\le
|\langle \varphi_j, g_j .\psi_j\rangle |\cdot
( \prod_{i\ne j} \|\varphi_i\|\cdot \|\psi_i\|). $$
If $\varphi_j'$ denotes the projection of $\varphi_j$ to $L^2_{00}(\bG_j(K)\backslash \bG_j(\A))$,
then
$$\langle \varphi_j, g_j .\psi_j\rangle= \langle \varphi_j', g_j .\psi_j\rangle .$$
Since $g\to \infty $ strongly and hence $g_j\to\infty$, we obtain
$\langle \varphi_j', g_j .\psi_j \rangle \to 0$ by (\ref{peter}). This proves
the claim.
\end{proof}

By the following proposition, the above theorem applies
to connected semisimple adjoint (simply connected) $K$-groups.
\begin{Prop}\cite[3.1.2]{Ti1}\label{decom} Any connected semisimple adjoint (resp. simply connected) $K$-group
decomposes into a direct product of adjoint (resp. simply connected)
almost $K$-simple $K$-groups.
\end{Prop}

\subsection{Equidistribution of Hecke points}\label{hecke}
In this subsection which is not needed in the rest of the paper,
 we explain applications of the adelic mixing in
the equidistribution problems of Hecke points considered in \cite{COU}.
Let $K=\q$. Let $S$ be a
finite set of primes including the archimedean prime $\infty$. If
$\Gamma$ is an $S$-arithmetic subgroup of $\bG_S$
 (here $\q_\infty=\br$) and $a\in\bG(\q)$, then the Hecke
operator $\op{T}_a$ on $L^2(\G\backslash \bG_S)$ is defined by
$$\op{T}_a(\psi)(g)=\frac{1}{\deg(a)}\sum_{x\in\G\backslash \G a\G}\psi(xg)$$
 where $\deg(a)=\#\G\backslash \G a\G$.
Theorem \ref{mmix} extends the main result in \cite{COU} where some cases
of $\q$-anisotropic groups were excluded (see \cite{EO}). In fact, the following
corollary immediately follows from Theorem \ref{am} and
Proposition 2.6 in \cite{COU}:
\begin{Cor} Let $\bG$ be a connected simply connected
almost $\q$-simple
 $\q$-group and $S$ a finite set of primes
including $\infty$. Suppose that $\bG_S$ is non-compact. Let $\Gamma
\subset \bG(\q)$ be an $S$-congruence subgroup of $\bG_S$.
 For any $\e>0$, there exists a constant $c=c(\Gamma, \e)>0$ such that
$$\|\op{T}_a\|\le c\cdot \tilde{\xi}_{\bG}^{1-\e}(a)\quad
\text{for any $a\in \bG(\q)$}. $$
\end{Cor}
This corollary in particular implies that for any sequence $a_i\in
\bG(\q)$ with $\deg(a_i)\to \infty$, and for any $\psi\in C_c(\bG_S)$,
$$\lim_{i\to \infty} \frac{1}{\deg(a_i)}\sum_{x\in \G a_i\G}
\psi(x)=\frac{1}{\tau_S(\Gamma\ba \bG_S)}
\int_{\bG_S} \psi(g)\,d\tau_S .$$ 
It is interesting to note that unlike the rational points $\bG(\q)$
of bounded height (Theorem \ref{e}), the Hecke points are
equidistributed  in $\bG_S$ with respect to the invariant measure.

The following corollary presents a 
stronger version of property $(\tau)$ of $\bG$ proved by Clozel \cite{Cl1}:

\begin{Cor}
 Let $\bG$ be a connected simply connected almost $K$-simple
 $K$-group. Let $\pi$ denote the quasi-regular representation of
$\bga$ on $L^2_0(\bG(K)\backslash \bG(\A))$. Let $W$ be a maximal
compact subgroup of $\bga$.
Then there exist an explicit $p=p(\bG)<\infty $ such that
any $W$-finite matrix coefficient of $\pi$
is $L^p(\bga)$-integrable. 

\end{Cor}

\section{Volume asymptotics and 
Construction of $\tilde\mu_{\iota}$}\label{s-volume}
\subsection{Analytic properties of height zeta functions}
Let $\bG$ be a connected adjoint semisimple algebraic group
over $K$.
Let $\iota:\bG\to \GL_N$ be a faithful
representation defined over $K$ which has a unique maximal weight.
Recall the constants $a_\iota$ and $b_\iota$ defined in the
introduction (\ref{abdef}): Choosing a maximal
 torus $\mathbf {T}$ of $\bG$ defined over $K$ containing
a maximal $K$-split torus and a set of simple roots $\Delta$
of the root system $\Phi(\bG, \mathbf T)$,
 denote by $2\rho$ the sum of all positive roots and $\lambda_\iota$
the unique maximal weight of $\iota$.

If $$2\rho=\sum_{\alpha\in \Delta} u_\alpha\alpha\quad\text{and}\quad
 \lambda_\iota=\sum_{\alpha\in \Delta} m_\alpha \alpha $$
then
 \begin{equation}\label{abdeff}a_\iota=\max_{\alpha\in \Delta}\frac{u_\alpha+1}{m_\alpha}\quad\text{and}\quad
b_\iota =\# \{\Gamma_K.\alpha: \frac{u_\alpha+1}{m_\alpha}=a_\iota \}\end{equation}
where $\Gamma_K$ is the absolute Galois group over $K$.

We fix a height function $\He_\iota=\prod_{v\in R}\He_{\iota,v}$
 as defined in \eqref{hed} 
in the rest of this section.
Recall the notation $$B_T:=\{g\in \bga:\He_\iota(g)\le T\}.$$

Given an automorphic character $\chi$, we consider the following
functions:
$$\Cal Z_S(s,\chi):=\int_{\bG_S}\He_\iota(g)^{-s}\chi(g)
\,d\tau_S(g);\quad\quad \Cal Z^S(s,\chi):=\int_{\bG^S}\He_\iota(g)^{-s}\chi(g)
\,d\tau^S(g).$$

\begin{Lem}\label{mv}  There exists $\e>0$ such that 
 the following hold for any finite $S\subset R$:
\begin{enumerate}
\item $\mathcal Z_{S}(s,\chi)$ absolutely converges for $\Re(s)\ge a_\iota -\e$.
 \item
 $\tau_S(B_T\cap \bG_S)=O(T^{a_\iota-\e})$
 where the implied constant depends on $S$.
\end{enumerate}\end{Lem}
\begin{proof}
For (1), Recall the Cartan decomposition for each $v$:
$\bG(K_v)=U_v A_v^+\Omega_vU_v$ \eqref{smodel}.
Since the definition of $a_\iota$ does not depend on a particular
choice of $\mathbf T$, we may assume $A_v\subset \mathbf T(K_v)$.
Choose the set of simple roots $\Delta_v=\{\alpha_1, \cdots, \alpha_r\}$ 
in the root system $\Phi(\bG(K_v), A_v)$
so that the restriction of $\Delta$ to $A_v$ is contained in
$\Delta_v\cup\{0\}$.

If $2\rho_v$ denotes the sum of all positive roots in
$\Phi(\bG(K_v), A_v)$ and
 $u'_{i}$ denotes the sum
of all $u_\alpha$'s for those $\alpha$ such that $\alpha|_{A_v}=1$,
i.e., $u'_i
=\sum\{u_\alpha: \alpha|_{A_v}=\alpha_i\}$,
 then $2\rho_v=2\rho|_{A_v}=
\sum_{i=1}^r 
u'_i\alpha_i .$
Similarly, if $\lambda_{\iota, v}:=\lambda_\iota|_{A_v}$
and $m'_{i}=\sum\{m_\alpha: \alpha|_{A_v}=\alpha_i\}$,
we have $\lambda_{\iota, v}=
\sum_{i=1}^r 
m'_{i} \alpha_i .$

Observe that $$a_{v}:=\max_{1\le i\le r}\frac{u_{i}'}
{m_{i}'} \le 
\max_{\alpha\in \Delta}\frac{u_\alpha}{m_\alpha}<a_{\iota} .$$

Since $\lambda_\iota$ is the unique maximal weight of $\iota$,
we may assume, without loss of generality, that
$$\He_{\iota, v}(k ad  k')=q_v^{\log_{q_v}|\lambda_\iota (a)|}$$
where $k, k'\in U_v, a\in A_v^+, d\in \Omega_v$ and
$q_v=e$ if $v\in R_\infty$.

For $v\in R_\infty$, it is well known (cf. Prop. 5.28 in \cite{Kn}) that
 $dg_v=\delta(X)dk_1 dX dk_2$ where
 for any $\e>0$,
there exists $C_\e>0$ such that
$$\delta(X)<C_\e \exp((1+\e)2\rho_v(X))\quad\text{for all $X\in \log(A_v^+)$ }.$$

Hence if $\sigma>0$,
\begin{align*}
&\int_{\bG(K_v)} \He_{\iota,v} (g_v)^{-\sigma} dg_v \le
C_\e \int_{\log A_v^+}\exp(-\sigma(\lambda_{\iota,v}(X)+(1+\e)2\rho_v(X)))dX
\\
&\le C_\e \prod_{i=1}^r
 \int_{x_{i}=0}^\infty\exp(-x_{i}
(\sigma m'_{i} -u'_{i} (1+\e)))dx_{i}
 .\end{align*}
Hence the above converges for any $\sigma 
> a_v$, 
proving the claim for $v\in R_\infty$.

Let $v\in R_f$. Without loss of generality, we may assume that
$\He_{\iota, v}$ is bi-$U_v$-invariant, and hence
$$\int_{\bG(K_v)} \He_{\iota,v} (g_v)^{-\sigma} dg_v =
 \sum_{ad\in A_v^+\Omega_v}\He_{\iota,v} (ad)^{-\sigma}
\tau_v(U_v ad U_v) .$$
By \cite[Lemma 4.1.1]{Sil},  there exists $c_1>0$
such that $\tau_v(U_v ad U_v)\le c_1\cdot q_v^{2\rho_v(a)}$ for all
 $ad\in A_v^+\Omega_v$.
Hence for some constant $c>0$,
$$\int_{\bG(K_v)} \He_{\iota,v} (g_v)^{-\sigma} dg_v \le
c \cdot \sum_{a\in A_v^+} q_v^{-\sigma \lambda_{\iota, v}(a) +2\rho_v(a)}
= c \prod_{i=1}^r\sum_{j=0}^\infty
 q_v^{-(\sigma m'_i -u'_i) j},$$
where the last term converges for any $\sigma 
> a_v$.
Put $$\e=\frac12 (a_\iota-
\max_{\alpha\in \Delta}\frac{u_\alpha}{m_\alpha}) $$
so that $a_\iota-\e >a_v$ for all $v\in R$.
The above argument proves the claim (1) for this choice of $\epsilon$. 
Also,
if $v(t):=\tau_S(\{g\in \bG_S: \He_{\iota, S}(g)\le t\})$, and
 $\sigma>a_\iota-\e$,
$$
\int_0^\infty t^{-\sigma}dv(t)=\int_{\bG_S}
 \He_{\iota,S} (g)^{-\sigma} d\tau_S(g)
<\infty
$$
Now the second claim follows from the properties of Laplace/Mellin
transform (see, for example, \cite[Ch.II,\S 2]{wi}). 
\end{proof}


One of the main contribution of the paper
by Shalika, Takloo-Bighash and Tschinkel \cite{STT2} is the regularization of
$\Cal Z^S(s, \chi)$ via the Hecke $L$-functions.
Their result stated as \cite[Theorem 7.1]{STT2},
together with the results in Tate's thesis on the meromorphic continuation of Hecke $L$-functions
and their boundedness on vertical strips (cf. \cite[Prop. 3.16]{Bu}),
implies the following:
\begin{Thm} \label{ab} Let $S$ be a finite subset of $R$ and $a_\iota$, $b_\iota$ as in (\ref{abdeff}).
Then $\Cal Z^S(s,\chi)$ converges absolutely when
$\Re(s)>a_\iota$, and  there exists $\e>0$ such that
$ \Cal Z^S(s,\chi)$
 has a meromorphic continuation to $\Re(s)>a_\iota-\e$ with a
unique pole at $s=a_\iota$
of order at most $b_\iota$. The order of the pole is exactly $b_\iota$ for $\chi=1$. Moreover,
for some constants $\kappa \in \br$ and $ C>0$,
$$\left|\frac{(s-a_\iota)^{b_\iota}\Cal Z^S(s,\chi)}{s^{b_\iota}}\right|\le
C\cdot   (1+|\op{Im}(s)|)^\kappa$$
for $\Re(s)>a_\iota-\e$.
\end{Thm}
In \cite{STT2}, it is assumed that $\He_{\iota, v}$ is smooth
for $v\in (R-S)\cap R_\infty$ which is stronger than
the condition (3) in Definition \ref{hed}. 
This implies Theorem \ref{ab} for any $S$ including $R_\infty$.
On the other hand,
 by Lemma \ref{mv}, for any finite $S_1\subset R$, there exists $\epsilon>0$ such that
 the product $\mathcal Z_{S_1}(s,\chi):=\int_{\bG_{S_1}}\He_\iota(g)^{-s}\chi(g)
\,d\tau_{S_1}(g)$ absolutely converges for all $\Re(s)>a_\iota-\epsilon$.
Therefore, for any $S_2\subset S$,
 the product $\mathcal Z^{S_2}=\mathcal Z_{S-S_2}\mathcal Z^S$ satisfies
the properties listed in Theorem \ref{ab}, provided $\mathcal Z^S$ does.
Therefore Theorem \ref{ab} holds for any finite $S\subset R$.

\vs We use the following version of Ikehara Tauberian theorem to
deduce the volume asymptotics 
 from Theorem \ref{ab}.
\begin{Thm}\label{ik} Fix $a>0$ and $\delta_0>0$.
Let $\alpha(t)$ be a non-negative non-decreasing function on
$(\delta, \infty)$ such that
$$f(s):=\int _{\delta_0}^\infty  t^{-s} \, d\alpha $$
converges for $\Re(s)>a$. Suppose that for some $\e>0$,
\begin{itemize}
\item  $f(s)$ has a meromorphic continuation to the half plane
$\Re(s)>a-\e>0$ and has a unique pole at $s=a$ with order
$b$;
\item For some $\kappa \in \br$ and $C>0$,
$$\left|\frac{f(s)(s-a)^b}{s^b}\right|\le C \cdot (1+|\op{Im}(s)|)^\kappa$$
for $\Re(s)>a-\e$.
\end{itemize}
Then for some $\delta>0$,
$$\int_{\delta}^T d\alpha=\alpha(T)-\alpha(\delta)
=\frac{c}{a(b-1)!}\cdot T^aP(\log T)+O(T^{a-\delta})\quad\text{as $T\to \infty$}$$
where $c=\lim_{s\to a} (s-a)^b f(s)$ and $P(x)$ is a monic polynomial
of degree $b-1$.
\end{Thm}
\begin{proof}
This can be proved by repeating the same argument as in the appendix
of \cite{CT1} simply replacing the sum $\sum_n  n^{-s}\alpha_n$
 by the integral $\int_{\delta_0} ^\infty t^{-s} \, d\alpha(t)$.
\end{proof}
\vs


\subsection{Definition of $\gamma_{W_f}^S$}
Recall from \eqref{wi} that $W_\iota$ denotes the maximal compact subgroup
of $\bG(\A_f)$ under which $\He_\iota$ is bi-invariant.
For any co-finite subgroup $W_f$ of $W_\iota$,
recall from (\ref{gwf}) the definition
$G_{W_f}:=\ker(\Lambda^{W_f}) $
where
$\Lambda^{W_f}$ is the subset
of $W_f$-invariant characters in $\Lambda$.
We will deduce the asymptotic volume of the intersection
$B_T\cap G_{W_f}$, more generally,
that of $B_T\cap g G_{W_f}\cap \bG^S$ for any $g\in \bga$ and any finite $S\subset R$.
In this subsection, we will define a function
$\gamma_{W_f}^S:\bga\to \br_{>0}$ which appears in the main asymptotic of 
these volumes.
\begin{Def}\label{eq:delta_S}
For a finite $S\subset R$ and a co-finite subgroup
$W_f$ of $W_\iota$, 
 define a function $\gamma^S_{W_f}:\bga\to \br_{>0}$ by
\begin{equation*}
\gamma^S_{W_f}(g):=\sum_{\chi\in\Lambda^{W_f}}
c^S_{\chi} \cdot \chi(g)
\quad\text{with}\quad c^S_{\chi}:= \lim_{s\to a_\iota} (s-a_\iota)^{b_\iota} 
\Cal Z^S(s,\chi).
\end{equation*}
For simplicity, when $S=\emptyset$,
we set $\gamma_{W_f}=\gamma^\emptyset_{W_f}$ and $c_{\chi}=
c^\emptyset_{\chi}$.
\end{Def}

By Theorem \ref{ab},
the limits appearing in the definition of $\gamma^S_{W_f}$ 
exist. To show $\gamma_{W_f}^S$ is well-defined,  it remains to show 
 the following:

\begin{Prop}\label{gammas} 
For any $S\subset R$ and
any co-finite subgroup $W_f$ of $W_\iota$,
$\gamma^S_{W_f}(g)>0$ for any $g\in \bga$.
\end{Prop}

 We need some preliminaries to prove this proposition.
\begin{Lem}\label{ker} The following statements hold for any compact open subgroup
$W_f$ of $\bG(\A_f)$:
\begin{itemize}
\item[(1)] If $\bG_\infty^\circ$ denotes 
the identity component of $\bG_\infty$, 
 $$\bG(K) \bG_\infty^\circ W_f\subset G_{W_f} .$$

\item[(2)] $\#\Lambda^{W_f}=  [\bga:G_{W_f}]<\infty .$
\item[(3)] For any $g\in \bga$,
$$\sum_{\chi\in \Lambda^{W_f}} \chi(g)=\begin{cases} \# \Lambda^{W_f}&\text{if $g\in G_{W_f}$} \\
0 &\text{otherwise}\end{cases} .$$
\end{itemize}
\end{Lem}
\begin{proof}
 Since $\bG_\infty^\circ$ is a connected semisimple group,
 $\bG_\infty^\circ\subset \ker (\chi)$ for any $\chi\in \Lambda$.
On the other hand, $\chi (\bG(K))=1$ for any $\chi\in \Lambda$, by the definition
of an automorphic character. Hence (1) follows.
Since $\bG_\infty^\circ$ has a finite index in $\bG_\infty$, it
follows from \cite[Theorem 5.1]{PR}
 that there exist finitely many $u_1, \cdots, u_h\in \bga$ such that
$$\bga=\cup_{i=1}^h \bG(K) u_i \bG_\infty^\circ W_f.$$
It follows by (1) that $[\bga:G_{W_f}]<\infty$. Now the
quotient $G_{W_f}\ba \bga $ is a finite abelian group.
In particular, $G_{W_f}$ is an open subgroup
of $\bga$ and hence any character of the group
$G_{W_f}\ba \bga $ can be considered as a continuous character of
$\bga$ which is trivial on $G_{W_f}$, that is, an element
of $\Lambda^{W_f}$. Conversely, any
element of $\Lambda^{W_f}$ defines a character of
$G_{W_f}\ba \bga $.

Now consider the scalar product on
the space functions of $G_{W_f}\ba \bga $ given by
$$\langle \psi_1 |\psi_2 \rangle:=\frac{1}{[\bga: G_{W_f}]}\sum_{x\in
G_{W_f}\ba \bga }\psi_1(x)\overline{ \psi_2(x)} .$$
We claim that $\Lambda^{W_f}$ forms an orthonormal
set with respect to this scalar product.
For $\chi, \chi'\in \Lambda^{W_f}$,
observe that
\begin{align*} &\sum_{x\in
G_{W_f}\ba \bga }\chi(x)\overline{ \chi'(x)}
\\&=\tau (\bG(K)\ba G_{W_f})^{-1}
\int_{g\in \bG(K)
\ba G_{W_f} }
\left( \sum_{x\in
G_{W_f} \ba \bga }\chi(gx)\overline{ \chi'(gx)} \right) d\tau (g)\\
&=\tau (\bG(K)\ba G_{W_f})^{-1}
\int_{\bG(K)
\ba \bga}\chi(x)\overline{ \chi'(x)} \; d\tau(x) .\end{align*}

Therefore
 we have
\begin{align*} \langle \chi |\chi' \rangle
&=
\frac{1}{[\bga: G_{W_f}]} \tau (\bG(K)\ba G_{W_f})^{-1} \int_{\bG(K)
\ba \bga}\chi(x)\overline{ \chi'(x)} \; d\tau(x)\\
&
=\tau(\bG(K)\ba \bga)^{-1}\langle \chi, \chi'\rangle_{L^2(\bG(K)\ba\bga)}
\\
&=\langle \chi, \chi'\rangle_{L^2(\bG(K)\ba\bga)}\end{align*}
since $\tau(\bG(K)\ba \bga)=1$.
Since $\Lambda$ is an orthonormal subset of 
$L^2(\bG(K)\ba\bga)$, the claim follows.
 Now (2) and (3)
 follow from the duality of finite groups.
\end{proof}

\begin{Lem}\label{positive} Fix a finite subset $S\subset R$.
 Let $U\subset \bG^S$ be an open subset
such that $\bG^S=F U$ for some finite subset $F$ of $ \bG^S$.
Then for $(\sigma\in \br)$
$$\liminf_{\sigma\to a_\iota} \;
(\sigma-a_\iota)^{b_\iota} \int_U \bH_\iota (h)^{-\sigma}
d\tau ^S(h) >0 .$$
\end{Lem}
\begin{proof}
 By Theorem \ref{ab},
$$c_0:=\lim_{s\to a_\iota}(s-a_\iota)^{b_\iota} \int_{\bG^S} \bH_\iota (h)^{-s}
d\tau ^S(h) $$ exists and is non-zero. It is then clear that $c_0>0$ since
$\bH_\iota$ is a positive function on $\bG^S$.
For any $f\in F$,
we can find $ c_f\ge 1$ such that for all $h\in \bG^S$,
 $c_f^{-1}\bH_\iota(h) \le \bH_\iota(fh)\le c_f\bH_\iota(h)$.
Without loss of generality, we may assume $\bH_\iota(h)\ge 1$ for all $h\in \bG^S$.
Now for $ \sigma<a_\iota+1$,
we have
$$\int_{\bG^S} \bH_\iota (h)^{-\sigma}
d\tau ^S(h) \le (\max_{f\in F} c_f )^{a_\iota+1} \cdot  \int_{U} \bH_\iota (h)^{-\sigma}
d\tau ^S(h).
$$
Hence the claim follows.
\end{proof}

\begin{Lem}\label{eq:haar} For any finite $S\subset R$ and a co-finite subgroup $W_f$ of $W_\iota$,
 there is a map
$g\mapsto s_g: \bG_S\to \bG^S$
 which factors through $(\bG_S\cap G_{W_f})\ba \bG_S$ for which the following hold:
\begin{enumerate}
 \item 
\begin{equation*}
G_{W_f}=\bigcup_{g\in (G_{W_f}\cap \bG_S)\ba \bG_S}
(G_{W_f}\cap\bG_S)g \, s_g(G_{W_f}\cap\bG^S);
\end{equation*}
\item for any $\varphi\in C_c(\bhi)$, \begin{equation*}
\int_{\bhi} \varphi \, d\tau_{W_f} =\int_{g\in\bG_S}\int_{h\in
\bG^S\cap \bhi} \varphi (gs_gh)\,\, d\tau^S(h)\, d\tau_S (g)
\end{equation*}
if $\tau^S$ and $\tau_S$ are normalized so that
$\tau_{W_f}=\tau^S\times \tau_S$ locally.

\end{enumerate} 
\end{Lem}
\begin{proof}
Let $\op{pr}$ denote the restriction of
 the projection map $\bga\to \bG_S$ to $G_{W_f}$.
Since $\bG(K)$ is dense in $\bG_S$
by the weak approximation
and the image $\op{pr}(G_{W_f})$ is an
 open subgroup
containing $\bG(K)$, the map $\op{pr}$ is surjective.

Note that
 $(G_{W_f}\cap \bG^S)(G_{W_f}\cap \bG_S)$ is a normal subgroup of $G_{W_f}$, and that 
the map $\pr$
induces an isomorphism, say $\tilde \pr$,
between $(G_{W_f}\cap \bG_S)(G_{W_f}\cap \bG^S)\ba G_{W_f}$ and $(G_{W_f}
\cap \bG_S)\ba \bG_S$.
For each $g\in (G_{W_f}
\cap \bG_S)\ba \bG_S$,
choose $s_g\in \bG^S\cap \tilde \pr^{-1}(g)$.
This yields the decomposition
$$G_{W_f}=\cup_{g \in (G_{W_f}\cap \bG_S)\ba \bG_S}  (G_{W_f}\cap
\bG^S)(G_{W_f}\cap \bG_S) s_g\; g.$$
Since $(G_{W_f}\cap \bG_S) s_g \, g =g\, s_g (G_{W_f}\cap \bG_S)$, (1) follows.
It is easy to deduce (2) from (1).
\end{proof}

\noindent{\bf Proof of Proposition \ref{gammas}: }
For $g\in \bG^S$, consider
\begin{align}\label{ff2}
\sum_{\chi\in \Lambda^{W_f}}\Cal
Z^S(s,\chi)\chi(g)&=
\sum_{\chi\in \Lambda^{W_f}} \int_{\bG^S} \He_\iota(h)^{-s}\chi(gh)\,
d\tau^S(h)\\&=(\#\Lambda^{W_f}) \int_{g^{-1}G_{W_f}\cap \bG^S} \He_\iota(h)^{-s}\,
d\tau^S(h)\notag
\end{align}
where the second equality holds by Lemma \ref{ker}(3).

Hence $$ (\#\Lambda^{W_f})^{-1} \cdot \gamma^S_{W_f}(g)=
\lim_{s\to a_\iota} (s-a_\iota)^{b_\iota} 
 \int_{g^{-1}G_{W_f}\cap \bG^S} \He_\iota(h)^{-s}\,
d\tau^S(h), $$
which is equal to, along $\sigma\in \br$,
$$\liminf_{\sigma
\to a_\iota}\; (\sigma-a_\iota)^{b_\iota}\int_{g^{-1} G_{W_f}\cap \bG^S} \He_\iota(h)^{-\sigma}\,
d\tau^S(h).$$
Therefore by Lemma \ref{positive},
$\gamma^S_{W_f}(g)>0$ for $g\in \bG^S$.
Since $\bga=\bG_S \bG^S$, it suffices to prove
$\gamma^S_{W_f}(hx)>0$ for any $h\in \bG_S$ and $x\in \bG^S$.
Let $s_h\in \bG^S$ be as defined in Lemma \ref{eq:haar}. By (1)
of the same lemma,
$$\chi(h)=\chi(s_h^{-1})$$
 for any $\chi\in \Lambda^{W_f}$.
This implies that for any $x\in \bG^S$,
$\gamma^S_{W_f}(hx)= \gamma^S_{W_f}(s_h^{-1}x)$.
Since $s_h^{-1}x\in \bG^S$, by Theorem \ref{vol},
$\gamma^S_{W_f}(s_h^{-1}x)>0$ and hence
$\gamma^S_{W_f}(hx)>0$.
This finishes the proof.

The functions $\gamma_{W_f}^S$ are related for different $S$'s by the following:
\begin{Prop}\label{ri}
Let $S\subset S'$ be finite subsets of $R$
and $W_f$  a co-finite subgroup of $W_\iota$,
\begin{enumerate}
 \item for any $g\in \bG^S$,
$$\gamma_{W_f}^{S}(g)
 =\int_{h\in \bG_{S_0}}  \He_{\iota,S_0} ^{-a_\iota}(h)
\, \gamma^{S'}_{W_f} (hg) \; d\tau_{S_0},
\quad\text{for $S_0=S'-S$}.
$$ 
\item In particular, for any finite $S\subset R$,
$$\gamma_{W_f}(e)
 =\int_{ \bG_{S}}  \He_{\iota,S} ^{-a_\iota}
\, \gamma^{S}_{W_f}  \; d\tau_{S} .$$
\item $\gamma_{W_f}=\sum_{\chi\in \Lambda} c_\chi \chi$; in particular,
 $\gamma_{W_f}=\gamma_{W_\iota}$.
\end{enumerate} 
\end{Prop}
\begin{proof}
Note that $\bG^S=\bG_{S_0}\bG^{S'}$ and $\tau^{S}=\tau_{S_0}\times \tau^{S'}$.
Since $\int_{\bG_{S_0}}  \He_{\iota,S_0} ^{-a_\iota}\, \chi\; d\tau_{S_0}$
exists for any $\chi\in \Lambda$ by Lemma \ref{mv}, we deduce 
\begin{align*} \gamma_{W_f}^{S}(g) &= \sum_{\chi\in \Lambda^{W_f}}
\left(\int_{\bG_{S_0}}
 \He_{\iota, S_0}^{-a_\iota}  \chi \, d\tau_S \right)  
 \cdot \left(\lim_{s\to a_\iota^+}(s-a_\iota)^{b_\iota} \mathcal Z^{S'}(s, \chi)
\right)\cdot \chi(g)
\\
&= 
\int_{h\in \bG_{S_0}}
 \He_{\iota, S_0}(h)^{-a_\iota}  \; 
\left(\sum_{\chi\in \Lambda^{W_f}} c_\chi^{S'}\; \chi(gh) \right) \; d\tau_{S_0}(h) 
\\&=
\int_{h\in \bG_{S_0}}
\He_{\iota, S_0}(h)^{-a_\iota} \gamma^{S'}_{W_f} (gh)\; d\tau_{S_0}(h).
 \end{align*}
Hence (1) follows. By putting $S=\emptyset$, (2) follows.

Since $\He_\iota$ is $W_f$-invariant,
it follows from Lemma \ref{wif} below $c_\chi=0$ for any $\chi\in \Lambda-\Lambda^{W_f}$. Hence
the claim holds.
\end{proof}

\begin{Lem}\label{wif} Let $Y=\bga$ or $\bG(K)\ba \bga$.
 Let $W_f$ be a co-finite subgroup of $W_\iota$ and let $\chi\in \Lambda-\Lambda^{W_f}$.
Then for any $W_f$-invariant function $\psi$ on $Y$,
we have
$$
\int_{Y}\chi(g)\psi(g)\, d\tau (g)
=0,
$$
if the integral exists.
In particular, if the support of $\psi$ is
contained in $\bG(K)\ba G_{W_f}$ and $\int_{\bG(K)\ba G_{W_f}}
\psi \; d\tau =0$,
then for any $\chi\in \Lambda$,
$$\int_{\bG(K)\ba G_{W_f}} \chi \cdot \psi \; d\tau =0 .$$

\end{Lem}
\begin{proof}
Since $\chi\in \Lambda -\Lambda^{W_f}$
there exists $w\in W_f$ such that $\chi(w)\ne 1$.
Since $\psi$ is $W_f$-invariant and $\tau$
is invariant,
\begin{align*}
&\int_{\bga}\chi(g)\psi(g)\, d\tau_0 (g)\\
&=
\int_{\bga} \chi(wg)\psi(wg)\, d\tau (g)
\\&=\chi(w) \int_{\bga} \chi(g)\psi(g)\, d\tau (g).
\end{align*}
This equality implies the first claim immediately.
For the second claim, it suffices to note that for $\chi\in \Lambda^{W_f}$,
$\chi=1$ on $G_{W_f}$.
\end{proof}

\subsection{Volume asymptotic}
\begin{Thm}\label{vol} Let $a_\iota\in \q^+$ and $b_\iota\in \n$ be as in (\ref{abdeff}). Then
for any finite subset $S\subset R$, any co-finite
subgroup $W_f$ of $W_\iota$  and  $g\in \bG^S$,
there exists
a monic polynomial $P(x)$ of degree $b_\iota-1$ and a positive
real number $\delta$ such that
\begin{equation}\label{eq:lead}
\tau^S(B_T\cap gG_{W_f}\cap \bG^S)
= \frac{\gamma^S_{W_f}(g^{-1})}
 {\#\Lambda^{W_f}\cdot a_\iota(b_\iota-1)!}\cdot
T^{a_\iota} P(\log T)+O(T^{a_\iota-\delta}).
\end{equation}
 
In particular,
as $T\to \infty$,
\begin{equation}\label{eq:lead2}
\tau^S(B_T\cap gG_{W_f}\cap \bG^S)
\sim  \frac{\gamma^S_{W_f}(g^{-1})}
 {\#\Lambda^{W_f}\cdot a_\iota(b_\iota-1)!}\cdot
T^{a_\iota} (\log T)^{b_\iota-1}.\end{equation}
\end{Thm}

\begin{proof}
By Lemma \ref{ff3}, $B_T$ is a relatively compact subset of $\bga$
and hence $\tau^S(B_T\cap \bG^S)<\infty$ for each $T\ge 1$ and for
any finite $S$. By the same lemma,
$\delta_0:=\inf _{g\in \bga}\He_\iota(g)>0$. Define
$$\alpha(t)=\tau^S(B_t\cap g G_{W_f}\cap
\bG^S)\quad\text{for }t\in [\delta_0, \infty) ,$$ and $$f(s)=
\int_{\delta_0}^\infty t^{-s}\, d\alpha.$$
Then by \eqref{ff2},
\begin{align*}
f(s)=\int_{gG_{W_f}\cap \bG^S} \He_\iota(h)^{-s}\,
d\tau^S(h)= (\#\Lambda^{W_f})^{-1}\sum_{\chi\in \Lambda^{W_f}}\Cal
Z^S(s,\chi)\chi(g^{-1}).
\end{align*}

Since
$$\lim_{s\to a_\iota} (s-a_\iota)^{b_\iota} f(s) =
(\#\Lambda^{W_f})^{-1}\gamma^S_{W_f}(g) >0$$ by Proposition \ref{gammas},
 the claim follows from Theorems \ref{ab} and \ref{ik}.
\end{proof}

Lastly in this subsection, we show that
the volume asymptotic for $
\tau_{W_f}(B_T\cap G_{W_f})$ is independent of $W_f$'s.

\begin{Cor}\label{tauf} For $g\in\bga$
 and any co-finite subgroup $W_f$ of $W_\iota$,
 $$\tau (B_T\cap g G_{W_\iota})
\sim_T [ G_{W_\iota}:
G_{W_f}]\cdot
\tau (B_T\cap gG_{W_f} ). $$
In particular,
 $$\lim_{T\to\infty}\frac{ \tau_{W_f}(B_T\cap G_{W_f})}{
  \tau_{W_\iota}(B_T\cap G_{W_\iota})}=1. $$
\end{Cor}
\begin{proof}
Since $
\gamma_{ W_f}
(g) = \gamma_{ W_\iota}
(g) $ by Proposition \ref{ri},
 the first claim follows from Theorem \ref{vol} and Lemma \ref{ker}(2).
Since the restriction of $\tau$ to
$G_{W_f}$ is equal to
$[ \bga:G_{W_f}]\cdot \tau_{W_f}$,
the second claim follows from the first one.
\end{proof}

\subsection{Construction of $\tilde \mu_{\iota}$}\label{muiota}
Let $\bar \iota:\bG\to\mathbb P({\op{M_N}})$
denote the projective embedding obtained by the composition of $\iota$ with
the canonical projection from $\GL_N \to  \mathbb P({\op{M_N}})$.
For each $v\in R$, denote by $X_{\iota, v}$ the closure of
$\bar\iota(\bG(K_v))$ in $\mathbb P(\M_N(K_v))$, and set
$$X_{\iota}=\prod_{v\in R}X_{\iota, v} .$$
Fix a height function $\He_\iota=\prod_{v\in R} \He_{\iota,v}$ on the
associated
 adele group $\bG(\A)$ relative to $\iota$ as in Definition \ref{hed}.
For finite subset $S\subset R$ and
set $\He_{\iota,S}= \prod_{v\in S} \He_{\iota, v}$
and $X_{\iota, S}=\prod_{v\in S}X_{\iota, v}$.
 Without loss of generality, we may consider $\bG(K_v)$ as a subset of $X_{\iota,v}$.
We will first construct a family measures $\{\mu_{\iota, W_f, S}\}$
 on $X_{\iota, S}$ for all finite $S$ and for all co-finite subgroups $W_f$ of $W_\iota$,
 and put them together to define the measure $\tilde \mu_\iota$ on $X_{\iota}$. 


Recall the definition of the function $\gamma^S_{W_f}$ on $\bga$
from \eqref{eq:delta_S}, and the notation $\gamma_{W_\iota}(e)
=\gamma_{W_\iota}^\emptyset (e)$.
By Lemma \ref{mv} and Propositions \ref{gammas}, \ref{ri} (2), the following 
 is a well-defined  probability measure  on $\bG_S$ (which will be in fact
considered as a measure on $X_{\iota, S}$):
\begin{align*}
d\mu_{\iota,W_f, S}(g)&:=
\gamma_{W_\iota}(e)^{-1}\cdot
 { \He_{\iota,S}(g)^{-a_\iota} \cdot \gamma^S_{W_f}(g) \, d\tau_S(g)}.
\end{align*}

\begin{Rmk}\label{ac}\rm
Let $\bG_S'$ denote the derived subgroup of $\bG_S$. Then $[\bG_S:\bG_S']<\infty$
\cite[Proposition 3.17]{PR}.
Then for any $\psi\in C(X_{\iota, S})$,
since the projection of $\gamma^S_{W_f}$ to $\bG_S$
 factors through $\bG_S'$, 
we deduce that $$\mu_{\iota, W_f, S}(\psi)=\gamma_{W_\iota}(e)^{-1}
\sum_{u\in \bG_S/\bG_S'} \gamma^S_{W_f}(u)\cdot
\int_{u\bG_S'}\He_{\iota, S}(g)^{-a_\iota}\psi(g) \;d\tau_S(g).$$
Since $\gamma^S_{W_f}(u)>0$ for each $u$,
it follows that the measure $\mu_{\iota, W_f, S}$ is 
equivalent to a Haar measure on $\bG_S$, considered as a measure
on $X_{\iota, S}$.
\end{Rmk}

For $W_f<W_\iota$,
denote by $C(X_\iota)^{W_f}$ the closed subspace of $C(X_\iota)$
consisting of functions which are (right)-invariant under $W_f$.

\begin{Thm}\label{smuiota-g}
 There exists a unique probability measure
  $\tilde \mu_\iota$ on $X_{\iota}$ such that
  for any $\psi\in \bigcup_{W_f<W_\iota\;\text{co-finite}}
C(X_\iota)^{W_f}$,
\begin{equation}\label{t:miws} \tilde \mu_{\iota}(\psi)
=\gamma_{W_\iota}(e)^{-1}\cdot \sum_{\chi\in \Lambda}
\lim_{s\to a_\iota^+}(s-a_\iota)^{b_\iota}
\int_{\bga} \He_{\iota}(g)^{-s} \chi(g)\;\psi(g) d\tau(g). \end{equation}
 \end{Thm}
\begin{proof}
Define a linear functional $\mu_{\iota, W_f}$ on $C(X_\iota)^{W_f}$ by
\begin{equation}\label{miws} \mu_{\iota, W_f}(\psi)
=\gamma_{W_\iota}(e)^{-1}\cdot \sum_{\chi\in \Lambda^{W_f}}
\lim_{s\to a_\iota^+}(s-a_\iota)^{b_\iota}
\int_{\bga} \He_{\iota}(g)^{-s} \chi(g)\;\psi(g) d\tau(g). \end{equation}

We first claim that 
 $\mu_{\iota, W_f}$ 
is well-defined, positive and bounded by $1$ and $\mu_{\iota, W_f}(1)=1$. 
For each finite set of places $S$,
let $C(W_f, S)$ denote the subset of
$C(X_\iota)^{W_f}$ consisting of functions which factor through
$X_{\iota, S}$.
The restriction of $\psi\in C(W_f, S)$
 to $X_{\iota, S}$ will also be denoted by
$\psi$ by abuse of notation.
By Proposition \ref{ri} (1),
the measures $\mu_{\iota, W_f, S}$ 
are compatible in the sense that
for $S\subset S'$, the restriction
of $\mu_{\iota, W_f, S'}$ to $C(W_f, S)$
coincides with
$\mu_{\iota, W_f, S}$.
Observe that for $\psi\in C(W_f,S)$,
\begin{align} \label{same}\mu_{\iota, W_f}(\psi)
  =\mu_{\iota, W_f,S}(\psi). \end{align}
Hence the limit exists in \eqref{miws} for all $\psi\in
C(W_f, S)$ for each finite $S$.
Since $\bigcup_S C( W_f, S)$ is dense in $C(X_\iota)^{W_f}$
where $S$ ranges over all finite subsets of $R$
and $\mu_{\iota, W_f}$
is a linear functional
with 
$$|\mu_{\iota, W_f}(\psi)|\le\|\psi\|_\infty \quad\text{
 for any $\psi\in C(X_\iota)^{W_f}$}$$
it follows that 
the limit exists in \eqref{miws} for any $\psi\in C_c(X_\iota)^{W_f}$
and hence $\mu_{\iota, W_f}$ is well-defined.
The other claims on $\mu_{\iota, W_f}$ are now clear.

By applying Lemma \ref{wif} for $\psi=\He_\iota^{-s}$,
 the family $\mu_{\iota, W_f}$ of linear functionals
on $C(X_\iota)^{W_f}$ is compatible, in the sense that
if $V_f\subset W_f$ are co-finite subgroups of $W_\iota$, then
$$\mu_{\iota, V_f}|_{C(X_\iota)^{W_f}}=
\mu_{\iota, W_f} .$$


Hence \eqref{t:miws} is well-defined on
$\mathcal C_0:=\bigcup_{W_f<W_\iota\;\text{co-finite}}
C(X_\iota)^{W_f}$.

Since $\mathcal C_0$ is dense
in $C(X_\iota)$, 
there exists a unique positive linear functional $\tilde \mu_\iota$
on $C(X_\iota)$ satisfying \eqref{t:miws}. \end{proof}

\begin{Prop}\label{accc} For any finite $S\subset R$,
 the projection $\tilde{\mu}_{\iota,S}$ of $\tilde{\mu}_{\iota}$ on $X_{\iota,S}$ is equivalent to 
 a Haar measure of $\bG_S$.
\end{Prop}

\begin{proof}
 Recall that $\bG_S'$ denotes the derived subgroup of $\bG_S$.
We first claim that if $V_f$, $W_f$ are co-finite subgroups 
of $(W_\iota \cap \bG'_S) \times (W_\iota \cap \bG^S)$, then $\mu_{\iota,W_f,S}=\mu_{\iota,V_f,S}$. Indeed, by definition of
 $\mu_{\iota,W_f,S}$ and $\gamma^S_{W_f}$, it is sufficient to show that for all $\chi \notin \Lambda^{W_f}$, 
$c_\chi^S=0$, and so by symmetry:
  $$\gamma^S_{V_f}(g)=\sum_{\chi \in \Lambda^{V_f}\cap \Lambda^{W_f}} c_\chi^S \cdot  \chi(g) =\gamma^S_{W_f}(g).$$
  Let $\chi \notin \Lambda^{W_f}$,
so that $\chi(w)\neq 1$ for some $w \in W_f$. 
Write $w=w_S w^S$, where $w_S \in \bG_S'$ and $ w^S \in \bG^S$.
 Then $\chi(w^S)=\chi(w) \neq 1$ since $w_S \in \bG_S'$,
and hence
  $$c_\chi^S= \lim_{s \rightarrow a_\iota} (s-a_\iota)^{b_\iota} \int_{\bG^S} 
\He_\iota(g)^{-s}\chi(g)d\tau^S(g).$$
  Since $W_f\subset W_\iota$, and  $w^S \in W_\iota$,
  $$\int_{\bG^S} \He_\iota(g)^{-s}\chi(g)d\tau^S(g)=\int_{\bG^S} \He_\iota(gw^S)^{-s}\chi(gw^S)d\tau^S(g)$$
  $$= \chi(w^S) \int_{\bG^S} \He_\iota(g)^{-s}\chi(g)d\tau^S(g),$$
  which
 proves that this integral is zero for all $s$, and so $c_{s,\chi}=0$.
   
  Define $\tilde{\mu}_{\iota,S}=\mu_{\iota,W_f,S}$ for any $W_f \subset (W_\iota \cap \bG'_S) \times (W_\iota \cap \bG^S)$. 
This measure is absolutely continuous, by the remark \ref{ac} on
   $\mu_{\iota,W_f,S}$. We claim that $\tilde{\mu}_{\iota,S}$
is precisely
 the projection of $\tilde{\mu}_\iota$.
 For $\psi \in C(W_f,S)$ for some
 $W_f<(W_\iota \cap \bG'_S) \times (W_\iota \cap \bG^S)$, we have
  $$\tilde{\mu}_\iota(\psi)=\mu_{\iota,W_f}(\psi)= 
\mu_{\iota,W_f,S}(\psi)=\tilde{\mu}_{\iota,S}(\psi).$$
  Since the union $\cup_{W_f} C(W_f,S)$,
where $W_f$ ranges over co-finite subgroups of $\bG(\A_f)$
contained in $(W_\iota \cap \bG'_S) (W_\iota \cap \bG^S)$,
 is dense in $C(X_{\iota,S})$, this finishes the proof.
 \end{proof}

Let $\iota:\bG\to \hbox{GL}_N$ be an absolutely irreducible representation defined over $K$
with the highest weight $\lambda_\iota$. Set
$$
\Delta_\iota=\{\alpha\in\Delta:\, \frac{u_\alpha+1}{m_\alpha}=a_\iota\}.
$$
For $\alpha\in\Delta$, we denote by $\check{\alpha}$ the corresponding coroot.
It follows from Theorem 7.1 in \cite{STT2} that if for a finite subset $S\subset R$
and an automorphic character $\chi$,
$$
c_{\chi}^S=\lim_{s\to a_\iota^+} (s-a_\iota)^{b_\iota} \int_{\bG^S}
\He_\iota(g)^{-s} \chi(g)\,dg\ne 0,
$$
then
\begin{equation}\label{eq:tshn1}
\chi(\check{\alpha})=1\quad\hbox{for all $\alpha\in \Delta_\iota$},
\end{equation}
and conversely if (\ref{eq:tshn1}) holds, then $c_{\chi}^S\ne 0$ for all
sufficiently large $S\subset R$.

\vs
\begin{Rmk} \label{sec:ex}\end{Rmk}
 We discuss some examples to illuminate the
properties of the measure $\tilde \mu_\iota$.
\begin{enumerate}
\item Suppose that $\lambda_\iota$ is a multiple of
  $2\rho+\sum_{\alpha\in\Delta} \alpha$. In particular, this holds
  for $\lambda_\iota$ corresponding to the anticanonical
class and for all rank 1 groups.

Then $\Delta_\iota=\Delta$.
If a character $\chi\in\Lambda$ satisfies (\ref{eq:tshn1}) then it follows from
  the Cartan decomposition (\ref{eq:cartan}) that $\chi(\bG^S)=1$ for sufficiently large
  $S$ and by the weak
  approximation, $\chi=1$.
This shows that $c_{\chi}^S=0$ for every finite $S\subset R$ and every
  $\chi\in \Lambda^{W_f}-\{1\}$, so $\gamma^S_{W_f}$
is equal to
 $\lim_{s\to a_\iota} (s-a_\iota)^{b_\iota} \mathcal Z^S(s, 1)$.
Hence by Theorem \ref{vol},
\begin{equation}\label{eq:num}
\#\bG(K)\cap B_T \sim_T \tau(B_T);
\end{equation}
and
\begin{equation}\label{eq:m}
\tilde \mu_\iota=\prod_{v\in R} \frac{\He_{\iota,v}(g_v)^{-a_\iota}\, d\tau_v(g_v)}{\int_{\bG(K_v)}
  \He_{\iota,v}(g_v)^{-a_\iota}\, d\tau_v(g_v)}.
\end{equation}

\item Suppose that
 $K$ has class number one, $\bG$ is $K$-split
 and $W_\iota=\prod_{v\in R_f} \bG(\mathbb{\mathcal O}_v)$ 
(with respect to the canonical model
      over the ring $\mathcal{O}$ of integers). 

According to Remark in Section 2 in
      \cite{GaO},
$$\bG(\A)=\bG({K}) \bG_\infty^\circ W_\iota.
$$
Hence, $\Lambda^{W_\iota}=\{1\}$, 
and consequently, \eqref{eq:num} holds and
$\mu_{\iota, W_\iota}$ is given by \eqref{eq:m}.

\item
 (cf. Example 8.10, \cite{STT2}) Let $\bG=\hbox{PGL}_{4}$
and $\iota$ be the adjoint representation.
By \cite[\S8.2]{PR}, there exists a lattice $L\subset \mathfrak{pgl}_4(K)$
(i.e., an $\mathcal O$-module of full rank) such that
$\bG$ has class number 2 with respect to $L$.
We take the height function $\He=\prod_{v\in R} \He_v$ 
where $\He_v$ is the maximum norm with respect to $L$ for $v\in R_f$.
The group $W_\iota$ is given by 
 $\prod_{v\in R_f}\text{Stab}_{\bG(K_v)}(L\otimes \mathcal O_v)$.
By \cite[\S8.2]{PR}, $\bG(K)\bG_\infty W_\iota$ is a normal subgroup of
index 2 in $\bga$. If we additionally assume that the number field $K$ is
totally complex, then $\bG_\infty$ is connected and, hence,
$\Lambda^{W_\iota}=\{1,\chi\}$ for some automorphic character $\chi$ of order 2.
Every automorphic character of $\bG(\A)$ is of the form $\eta\circ
\det$ where $\eta$ is a Hecke character such that $\eta^4=1$. Since
the map $\det: \hbox{PGL}_4(K_v)\to K^\times_v/(K^\times_v)^4$ is
surjective for every $v\in R$,
it follows that $\chi=\eta\circ \det$ with $\eta^2=1$.
In this case, the roots and coroots are given by
\begin{align*}
\alpha_i(\hbox{diag}(a_1,\ldots,a_4))=a_ia_{i+1}^{-1},\quad
\check{\alpha}_i(t)=\hbox{diag}(\underbrace{t,\ldots, t}_i,1,\ldots,1)
\end{align*}
for $i=1,2,3$, and
$$
\lambda_\iota=\alpha_1+\alpha_2+\alpha_3,\quad 2\rho=3\alpha_1+4\alpha_2+3\alpha_3.
$$
Hence, $a_\iota=5$, $b_\iota=1$, $\Delta_\iota=\{\alpha_2\}$.
Then (\ref{eq:tshn1}) is equivalent to $\eta^2=1$,
and we deduce that $c_{\chi}^S\ne 0$ for sufficiently large finite
$S\subset R$. Since the function $\gamma^S_{W_\iota}=c_{1}^S+c_{\chi}^S\chi$
restricted to $\bG_S$ is not constant for sufficiently large $S\subset
R$, we conclude that
$$
\mu_{\iota,W_\iota, S}\ne \prod_{v\in S} \frac{\He_{\iota,v}(g_v)^{-a_\iota}\, 
d\tau_v(g_v)}{\int_{\bG(K_v)}
  \He_{\iota,v}(g_v)^{-a_\iota}\, d\tau_v(g_v)}.
$$

We also note that in this case, Theorem \ref{vol} implies that for 
an automorphic character $\chi$
such that $c_{\chi}\ne 0$, we have 
$$
\lim_{T\to\infty} \frac{\tau(B_T\cap
  G_{W_\iota})}{\tau(B_T)}=c_{1}^{-1}\cdot
\frac{1}{2}(c_{1}+c_{\chi})\ne \frac{1}{2}.
$$
In particular, it may happen that  in Theorem \ref{eqin},
$\tau(B_T)$ is not asymptotic to
$[\bga:G_{W_\iota}]\cdot\tau(B_T\cap G_{W_\iota})$ as $T\to\infty$.

\end{enumerate}

\subsection{Equidistribution of height balls $B_T\cap G_{W_f}$
with respect to $\tilde \mu_\iota$}
\begin{Prop}\label{measure} Let $W_f<W_\iota$ be a co-finite subgroup. 
Then for any $\psi\in C(X_\iota)^{W_f}$,
$$\lim_{T\to\infty} \frac{1}{\tau_{W_f} (B_T\cap G_{W_f})}
 \int_{B_T\cap G_{W_f}} \psi \,d\tau_{W_f}
= \tilde \mu_{\iota} (\psi). $$
 \end{Prop}
\begin{proof}
To prove the proposition, we may assume
$\psi\in C(W_f, S)$ for some finite $S$, since these
functions form a dense subset of $C(X_\iota)^{W_f}$. 
Let $\chi_{B_T}$ denote the characteristic function of
the set $B_T$.
By Lemma \ref{eq:haar} (2), we have a map $g\in \bG_S\mapsto s_g\in \bG^S$ such that 
$$ \int_{B_T\cap \bhi} \psi \, d\tau_{W_f}=
\int_{g\in\bG_S} \psi (g) \int_{h\in
\bG^S\cap \bhi}  \chi_{B_T}(gs_gh) \,\, d\tau^S(h)\, d\tau_S (g).
$$

Since $\He_{\iota}(gs_g h)=\He_{\iota, S}(g)\He_\iota(s_g h)$,
\begin{align*}  \int_{h\in
\bG^S\cap \bhi}  \chi_{B_T}(gs_gh) \,\, d\tau^S(h)
&=\tau^S\{ h\in  
\bG^S\cap \bhi : \He_{\iota}(s_g h) < T \He_{\iota,S}(g)^{-1} \}
\\ &=
\tau^S(s_g^{-1} B_{T\cdot \He_{\iota,S}(g)^{-1}}\cap G_{W_f}
\cap \bG^S)\\
& =\tau^S(B_{T\cdot \He_{\iota,S}(g)^{-1}}\cap s_gG_{W_f}
\cap \bG^S)
\end{align*}

Hence 
\begin{equation}\label{eq:dom} \int_{B_T\cap \bhi} \psi \, d\tau_{W_f}=
\int_{g\in\bG_S} \psi(g)  \tau^S(B_{T\cdot \He_{\iota,S}(g)^{-1}}\cap s_gG_{W_f}
\cap \bG^S)
d\tau_S (g).\end{equation}


Setting
$$y_T(g):=\frac{\tau^S (B_{T \cdot \He_{\iota,S}(g)^{-1}}\cap s_gG_{W_f}\cap
 \bG^S)}{\tau^S (B_{T}\cap G_{W_f}\cap \bG^S)},$$
we claim that for some constant $C>0$,
\begin{equation}\label{yt}
y_T(g)\le C \cdot
 \He_{\iota,S}(g)^{-a_\iota}\quad\text{ for any $g\in \bG_S$.}
\end{equation}

By \eqref{eq:lead2},
for any $h\in \bG^S$,
there exists a constant $c_h \ge 1$ such that for any
$T>2$,
\begin{equation}\label{ab-re} c_h^{-1} \cdot T^{a_\iota}
(\log T )^{b_\iota-1}\le
 \tau^S (B_{ T }\cap hG_{W_f}\cap \bG^S)
 \le c_h \cdot T^{a_\iota}
(\log T )^{b_\iota-1}.\end{equation}
Since $\{s_g\in \bG^S: g\in \bG_S\}$ is a finite set,
$c:= \max\{ c_{s_g} : {g\in \bG_S}\} <\infty$. 
It follows by \eqref{ab-re} that for any $g\in \bG_S$ with $\He_{\iota,S}(g)\le  T/2$,
$$ \tau^S (B_{ T \cdot \He_{\iota,S}(g)^{-1}}\cap s_gG_{W_f}\cap \bG^S)
 \le c\cdot \He_{\iota,S}(g)^{-a_\iota}T^{a_\iota}
(\log T \He_{\iota,S}(g)^{-1})^{b_\iota-1}.$$
On the other hand
 for any $g\in \bG_S$ satisfying $ T/2\le  \He_{\iota,S}(g) \le T \delta_0^{-1}$
where $\delta_0:=\inf_{g\in \bga}
 \He_\iota (g)>0 $ (see Lemma \ref{ff3}),
$$ \tau^S (B_{ T \cdot \He_{\iota,S}(g)^{-1}}\cap s_gG_{W_f}\cap \bG^S) 
\le \tau^S (B_{2}\cap  \bG^S)\le
d\cdot \He_{\iota,S}(g)^{-a_\iota} T^{a_\iota}.
$$
where $d={\delta_0^{-a_\iota}}{\tau^S (B_{2}\cap  \bG^S)}$.
Also, for $\He_{\iota,S}(g) > T \delta_0^{-1}$, $y_T(g)=0$.

Hence by applying \eqref{ab-re} once more now to $h=e$,
 we obtain the inequality (\ref{yt}). Since $\He_{\iota,S}^{-a_\iota}\in L^1(\bG_S)$
by Lemma \ref{mv},
 it follows that $y_T$ belongs to $L^1(\bG_S)$.
Since by Theorem \ref{vol},
$$
y_T(g)\to \gamma_{ W_f, S}(s_g^{-1})\He_{\iota,S}(g)^{-a_\iota}\gamma^S_{W_f}(e)^{-1}\quad\hbox{as $T\to\infty$,}
$$
 we apply the
dominated convergence theorem to (\ref{eq:dom}) and
deduce that
\begin{equation*}
\lim_{T\to \infty}
\frac{\int_{B_T\cap \bhi} \psi
 \, d\tau_{W_f}}{\tau^S(B_T\cap G_{W_f}\cap \bG^S)}
= \gamma^S_{W_f}(e)^{-1}\;\int_{\bG_S} \psi(g)\gamma^S_{W_f}(s_g^{-1})
\He_{\iota, S} (g)^{-a_\iota}  \;d\tau_S(g).
\end{equation*}

Using $\gamma^S_{W_f}(s_g^{-1})=\gamma^S_{W_f}(g)$ and the definition of
the measure $\mu_{\iota, W_f, S}$,
we have
\begin{equation}\label{fc}
\lim_{T\to \infty}
\frac{\int_{B_T\cap \bhi} \psi \, d\tau_{W_f}}{\tau^S(B_T\cap G_{W_f}\cap \bG^S)}
= \tilde \mu_{\iota}(\psi)\cdot\gamma_{W_f}^S (e)^{-1}
\end{equation}

 Taking $\psi=1$, we also get
\begin{equation}\label{fsc}
\lim_{T\to\infty}
\frac{\tau_{W_f}(B_T\cap \bhi)}{
 \tau^S(B_{T}\cap G_{W_f}\cap 
\bG^S)}=\gamma_{ W_f}^S(e)^{-1} .\end{equation}
Therefore combining (\ref{fc}) and (\ref{fsc}), we obtain
\begin{equation*}
\lim_{T\to\infty}\frac{\int_{B_T\cap \bhi} \psi \, d\tau_{W_f}}{
   \tau_{W_f}(B_T\cap
\bhi)}=\tilde \mu_\iota(\psi).
\end{equation*}

\end{proof}

\section{Equidistribution for saturated cases}\label{sec:proof}
Let $\bG$ be a connected adjoint semisimple
group defined over $K$ and $\iota$ be a faithful
 representation of $\bG$ to $\GL_N$ defined over $K$ which has
a unique maximal weight.
We recall that $\iota:\bG\to \GL_N$ is called {\it
 saturated} if the set
\begin{equation}\label{rs2}
\{\alpha\in \Delta: \frac{u_\alpha+1}{m_\alpha}=a_\iota
\}\end{equation} is not contained in the root system of a proper
normal subgroup of $\bG$. 

 Fix a height function $\He_\iota$ on
$\bga$ associated to $\iota$, and let $\tilde\mu_\iota$ be the probability measure
on $X_\iota$ constructed
in Theorem \ref{smuiota-g} (associated to $\He_\iota$).

\begin{Thm}\label{smuiota}
Suppose that $\iota$ is saturated. 
For $\psi \in C(X_\iota)$,
$$\lim_{T\to \infty}\frac{1}{\tau_{W_f} (B_T\cap G_{W_f})}
 \sum_{g\in\bG(K): \He_\iota(g)<T} \psi (g) = \int_{X_\iota} \psi \,
 d\tilde \mu_{\iota} $$
where $W_f$ is any-cofinite subgroup of $W_\iota$.
 \end{Thm}

Note that by Corollary \ref{tauf}, the above equality is independent of $W_f$.
To derive the asymptotic formula for the number of $K$-rational
points, it suffices to take $\psi=1$ in Theorem \ref{smuiota}, and
hence we obtain Theorem \ref{eqin}.
\begin{Cor} If $\iota$ is saturated,
 we have, as $T\to\infty$,
$$\# \{g\in \bG(K):\bH_\iota(g)<T\}\sim \tau_{W_f}(B_T\cap G_{W_f})
$$ 
for any co-finite subgroup $W_f$ of $W_\iota$.
\end{Cor}

Combining Theorem \ref{smuiota} with Theorem \ref{vol}, we
deduce Theorems \ref{mt} (without an error term)
and \ref{e} for the
saturated case. The rate of convergence as well as nonsaturated 
case are discussed in the next
section.

\begin{Rmk}\label{acc} \rm Note that
for any finite $S$,
the projection of $\tilde \mu_\iota$ to $X_{\iota, S}$
is $\tilde \mu_{\iota, S}$, which is equivalent to a Haar measure on $\bG_S$ by
Proposition \ref{accc}.
 Any open subset $X_\iota$ contains a subset of the form
$V_SX_{R-S}$ where $V_S$ is an open subset of $X_{\iota, S}$,
which again contains $(gW_f\cap \bG_S)\bG^S$ for some $g\in \bG_S$
and some co-finite subgroup $W_f$ of $W_\iota$.
Now
$$\tilde \mu_{\iota}((gW_f\cap \bG_S )\bG^S)=\tilde \mu_{\iota,S}(gW_f\cap \bG_S)>0.$$
This shows that $\tilde \mu_\iota$ has full support on $X_\iota$.
\end{Rmk}


The rest of section is devoted to a proof of Theorem \ref{smuiota}.
We recall the following facts:
\begin{Lem}\label{mn} Let $\bG$ be a connected semisimple adjoint $K$-group,
 $\iota:\bG\to \GL (V)$ a $K$-rational representation of $\bG$, 
and let $\bM$ be a connected
normal $K$-subgroup of $\bG$.
 \begin{enumerate}
\item  There
exists a connected normal $K$-subgroup $\bN$ of $\bG$ so that
$\bG=\bM \bN$ and $\bM\cap \bN=\{e\}$, and $\bM$ is semisimple adjoint.

\item For each $x\in \bN(\A)$,
the function $g \mapsto \He_\iota(gx)$ defines a height function on $\bM(\A)$
with respect to the restriction $\iota|_\bM$.
\item If $\iota$ has a unique maximal weight
$\lambda_\iota$, then
the restriction $\iota|_\bM$ has a  unique maximal weight.

\end{enumerate}
\end{Lem}
\begin{proof}
 (1) follows directly from Proposition \ref{decom}. Let $x\in \bN(\A)$.
Whenever $x\in \bN(K_v)\cap W_\iota$, which is the case
for almost all $v$, we have
$\bH_{\iota,v}(gx)=\bH_{\iota, v}(g)$ for all $g\in \bM(K_v)$.
Using this, it is easy to verify that
the function $g\mapsto \He_\iota(gx)$ is a height function as 
defined in Definition \ref{hedh}.

For (3), let $\bf T$ be a maximal torus of $\bG$  defined over $K$,
 $\Pi$ denote the set of all weights of $\iota$ with respect to $\bf T$,
and $\Pi'$ denote the set obtained by restricting elements of $\Pi$ to $\bM\cap \bT$.
If $V=\oplus_{\lambda\in \Pi} V_\lambda$ is the weight space decomposition for $\iota$,
the weight space decomposition for $\iota|_\bM$ is of the form
$V=\oplus_{\beta\in\Pi'} W_\beta$ where $W_\beta:=\oplus
\{  V_{\lambda} :{\lambda|_{\bT\cap \bM}=\beta}\}$. In particular, any weight of $\iota|_\bM$ is the
restriction of a weight of $\iota$ to $\bT\cap \bM$.
Hence if $\beta_\iota$ is the restriction of $\lambda_\iota$ to $\bT\cap \bM$,
and $\beta\in \Pi'$, not equal to $\beta_\iota$
then $\beta_\iota-\beta$ is a non-zero sum of positive roots of $\bM$ with respect to
 $\bM\cap \bT$. Therefore $\beta_\iota$ is the unique maximal weight of $\iota|_\bM$.
\end{proof}

\begin{Lem}\label{tfae}
The following are equivalent.
\begin{enumerate}
\item $\iota$ is saturated.
\item For any proper connected
 normal $K$-subgroup $\bM$ of $\bG$,
 \begin{equation*}
{\tau_{\bM}(B_T\cap \bM({\A}))}=O( T^{a_\iota}(\log T)^{b_{\iota}-2})
\end{equation*}
where $\tau_{\bM}$ is a Haar measure
 on $\bM(\A)$.\end{enumerate}
\end{Lem}
\begin{proof} Assume (1).
Since $\iota$ has a unique maximal weight,
 the restriction $\iota|_{\bM}$ of $\iota$ has a unique maximal weight
as well. There
exists a connected normal $K$-subgroup $\bN$ of $\bG$ so that
$\bG=\bM \bN$ and $\bM\cap \bN=\{e\}$. 
By Lemma \ref{l:prod},
without loss of generality, we may assume that
$\He_\iota $ is the product of height functions
$\He_{\iota|_{\bM}}$ and $\He_{\iota|_{\bN}}$ and $\He_{\iota}(e)=1$.
Hence $B_T\cap \bM({\A})=\{x\in \bM(\A):\He_{\iota|_{\bM}}(x) \le T\}$.
Hence
by Theorem \ref{vol}, 
$$\tau_{\bM} (B_T\cap \bM(\A))
\sim c \cdot T^{a}(\log T)^{b-1}$$
where $a$ and $b$ are defined as in
(\ref{abdeff}) for $\iota|_{\bM}$.
Now the saturated condition means that 
$a=a_\iota$ and
$b \le b_\iota-1$. Hence (2) holds.
To show the other direction,
suppose that $\iota$ is not saturated. Then there exists
 connected proper normal $K$-subgroups $\bM$ and $\bN$
 of $\bG$ as above 
such  that
 $\{\alpha\in \Delta: u_{\alpha}+1=a_\iota \cdot m_\alpha\}$
is contained in the root system of $\bM$.
By Theorem \ref{vol}, there exists $c>0$
such that for all large $T$, 
$$\tau_{\bM} (B_T\cap \bM(\A)) \ge c \cdot T^{a_\iota}(\log T)^{b_\iota-1}.$$
Hence (2) does not hold.
\end{proof}

The following lemma is  the main reason why we need
 the assumption of
$\iota$ being saturated for the proof of Theorem \ref{smuiota}.

By Lemma \ref{mn}, $\bG$ is a product of connected $K$-simple
adjoint subgroups.
\begin{Lem}\label{satur}
Suppose that $\iota$ is saturated.
Write $\bG=\bG_1\times \cdots\times \bG_m$ as a product of
connected $K$-simple subgroups.
Then for any fixed $C>0$,
 \begin{equation*}
\lim_{T\to \infty}\frac{\tau_{W_f}((B_T-B^C)\cap \bhi)}{\tau_{W_f}(B_T\cap \bhi)}=0
\end{equation*}
 where
 $ B^C:=\{(g_1,\cdots , g_m)\in \bG(\A): \He_\iota(g_i)>C \hbox{ for each
  $i=1,\cdots ,m$}\}.$
\end{Lem}
\begin{proof} Since $\bhi$ is non-compact,
$\tau_{W_f}(B_T\cap \bhi)\to \infty$ as $T\to \infty$.
If $m=1$, the claim follows immediately from this.
 Suppose $m\ge 2$.
Without loss of generality, we may assume that
$\He_\iota (g_1, \cdots, g_m)=\prod_{i=1}^m \He_\iota (g_i)$.
For each $i$, let $B_i^C$ denote the subset of $B_T$
consisting of $g=(g_1, \cdots, g_m)$ with $\He_\iota(g_i)\le C$.
If we denote by $\tau_i$ and $\tau^i$ Haar measures on 
$\bG_i(\A) $ and $\prod_{j\ne i} \bG_j(\A)$
such that $\tau_{W_f}=(\tau_i\times \tau^i)|_{\bhi}$, then
$$ \tau_{W_f}(B_i^C) \le C_0\cdot \tau^i(
\{g'\in \prod_{j\ne i} \bG_j(\A): \He_\iota(g')<\delta_i^{-1} \cdot T\}) $$
where $\delta_i=\inf_{g\in \bG_i(\A)} \He_\iota(g)$ and
$C_0=\tau_i(\{g\in \bG_i(\A): \He_\iota(g)\le C\})$.

By the previous lemma,  for some $c_i>0$,
$$\tau^i
(\{g'\in \prod_{j\ne i} \bG_j(\A): \He_\iota(g')<\delta_i^{-1} \cdot T\})
=O( T^{a_\iota}(\log T)^{b_\iota-2}).$$

Hence for each $i$,
$$\lim_{T\to \infty}\frac{\tau_{W_f}(B_i^C\cap \bhi)}
{\tau_{W_f}(B_T\cap \bhi)}=0.$$
 Since $B_T-B^C\subset \cup_{i=1}^m B_i^C$, this proves the claim.
\end{proof}

In the following, we 
fix a co-finite subgroup $W_f$
of $W_\iota$.
Set $$Y_{W_f}=\bG(K)\backslash G_{W_f} .$$

For a fixed $\psi \in C(X_\iota)^{W_f}$, we
define a function $F_T$ on $G_{W_f}\times G_{W_f}$ by
$$ F_T(g, h)=\sum_{\gamma\in \bG(K)} \psi(g^{-1}\gamma h)\cdot \chi_{B_T} (g^{-1}\gamma h). $$
Since $B_T$ is a compact subset of $\bga$, the above sum is finite and
since $F_T$ is $\bG(K)\times \bG(K)$-invariant,
we may consider $F_T$ as a function on $Y_{W_f}\times Y_{W_f}$

Note that 
$$F_T(e, e)=\sum_{\gamma\in \bG(K): \He_\iota(\gamma)\le T} \psi(\gamma).$$

\begin{Prop}[Weak-convergence]\label{wc}
Suppose that $\iota$ is saturated.
 For $i=1,2$, let
 $\alpha_i\in C(Y_{W_f})$ be a $W_f$-invariant function and $ \int_{Y_{W_f}}
 \alpha_i \, d\tau_{W_f} =1$.
If $\alpha(x, y):=\alpha_1(x)\alpha_2(y)$, then
  $$ \lim_{T\to \infty}\frac{1}{\tau_{W_f}(B_T\cap G_{W_f})} \int_{Y_{W_f}\times Y_{W_f}}F_T\cdot \alpha \,\, d(\tau_{W_f}\times
  \tau_{W_f}) =\int_{X_\iota} \psi \,
 d \tilde\mu_\iota   .$$
\end{Prop}

\begin{proof}
Observe that $\alpha$ is $W_f\times W_f$-invariant and 
\begin{align}\label{ft}
& \langle F_T, \alpha\rangle_{Y_{W_f}\times Y_{W_f}} \\
&=\int_{x\in Y_{W_f}}
\int_{y\in Y_{W_f}}
\left(\sum_{\gamma\in \bG(K)} \psi(x^{-1} \gamma y) \chi_{B_T}(x^{-1}\gamma y)\right)\,\alpha_1(x)\alpha_2(y) \, d\tau_{W_f} (y)
d\tau_{W_f}(x) \notag\\
&=
\int_{x\in Y_{W_f}}\int_{h\in \bhi}
\psi(x^{-1} h) \chi_{B_T}(x^{-1}h)\alpha_1(x)\alpha_2(h) \, d\tau_{W_f} (h)d\tau_{W_f}(x)\notag\\
&=
\int_{g\in \bhi }\psi(g) \chi_{B_T}(g) \left(\int_{x \in Y_{W_f}} \alpha_1(x) \alpha_2(xg)\, d\tau_{W_f} (x)\right)
d\tau_{W_f}(g)
\notag\\
&=\int_{g\in B_T\cap\bhi}\psi(g) \langle \alpha_1, g\cdot\alpha_2 \rangle \, d\tau_{W_f} (g) .
\notag\end{align}

As in the above lemma \ref{satur},
we write $\bG=\bG_1\times \cdots\times \bG_m$.
Now for a sequence $\{g\in \bga\}$, that $g\to \infty$ strongly
means precisely the $\bG_i(\A)$-component $g_i$ of $g$ tends to $\infty$
for each $1\le i\le m$.
Since the height function $\He_\iota$ is proper, this is
again equivalent to saying that  $\He_\iota(g_i)\to
\infty$ for each $i=1,\ldots,m$.

For $i=1,2$, we claim that $\alpha_i-1\in L^2_{00}(Y_{W_f})$,
that is,
for any automorphic character $\chi$ of $\bga$,
$$\int_{Y_{W_f}} \alpha \cdot \chi \;d \tau_{W_f} =\int_{Y_{W_f}}\chi  \;d \tau_{W_f}.$$
If $\chi$ is $W_f$-invariant, then $\chi=1$ on $G_{W_f}$.
Since $\int_{Y_{W_f}} \alpha \;d\tau_{W_f}=1$,  the claim is clear.
If $\chi\in \Lambda-\Lambda^{W_f}$,
we only need to apply Lemma \ref{wif} for $\psi=\alpha -1_{Y_f}$,
where $1_{Y_f}$ is the characteristic function of $Y_{W_f}$,
since we may consider $\psi$ as a function on $\bG(K)\ba \bga$ which is $0$ outside $Y_{W_f}$.

Hence, by Theorem \ref{smix}, for any given $\e>0$, there exists $C>0$ such
that for all
$g\in B^C$,
\begin{equation}
\label{alrate}\left| \langle \alpha_1, g\cdot\alpha_2 \rangle - 1
\right|=\left| \langle \alpha_1-1, g\cdot(\alpha_2 -1) \rangle\right| <\e \quad
\end{equation}
Hence
\begin{align*}
&\left| \int_{g\in B_T\cap\bhi} \psi(g) \langle  \alpha_1, g.\alpha_2 \rangle \, d\tau_{W_f} (g) -
\int_{g\in B_T\cap\bhi} \psi(g) \,d\tau_{W_f}(g) \right|\\
< &\sup |\psi|\cdot (\|\alpha_1\|\cdot \|\alpha_2\| +1)\cdot
\tau_{W_f}((B_T-B^C)\cap \bhi) +\e\cdot \sup |\psi| \cdot
 \tau_{W_f}(B_T\cap B^C \cap
\bhi)
\end{align*}
where $\|\alpha_i\|$ is the $L^2$-norm of $\alpha_i\in L^2(Y_{W_f})$ for each $i$.
By Lemma \ref{satur}, it follows that
$$\limsup_{T\to\infty}
\frac{1}{\tau_{W_f}(B_T\cap \bhi)} \left| \int_{g\in B_T\cap\bhi} \psi(g) (\langle  \alpha_1, g.\alpha_2 \rangle - 1) \,d\tau_{W_f}(g) \right|
\le \e \cdot \sup |\psi| .$$
 Since $\e>0$ is arbitrary, by (\ref{ft}), this proves 
  $$ \lim_{T\to \infty}\frac{1}{\tau_{W_f}(B_T\cap G_{W_f})} \int_{Y_{W_f}\times Y_{W_f}}F_T\cdot \alpha \,\, d(\tau_{W_f}\times
  \tau_{W_f}) =\lim_{T\to \infty}\frac{ \int_{B_T\cap\bhi} \psi\,d \tau_{W_f} }
{\tau_{W_f}(B_T\cap \bhi)}
.$$
By Lemma \ref{measure}, this proves Proposition \ref{wc}.
\end{proof}

\begin{proof}[\bf Proof of Theorem \ref{smuiota}]
 It suffices to prove our theorem for
nonnegative functions $\psi\in C(X_\iota)^{W_f}$ for each co-finite subgroup
$W_f$ of $W_\iota$.
 Fix $\e>0$. Let
$W_\infty$ be a symmetric neighborhood of $e$ in $\bG_\infty^\circ$
 such that
$$W_\infty B_T W_\infty\subset B_{(1+\e)T}\quad\text{and}\quad
B_{(1-\e)T}\subset \bigcap _{g, h\in W_\infty}g B_T h \quad\text{for
all $T>1$.}$$

By the uniform continuity of $\psi$,  replacing
$W_\infty$ by a smaller one if necessary, we may assume that
\begin{equation}\label{ne}
 \psi (g^{-1} x h)-\e\le \psi(x)\le \psi(g^{-1} x h) +\e\quad\text{
for all $x\in X$ and $g, h\in W_\infty$.} \end{equation}

Define
$$ F_T^{\pm}(g, h)=\sum_{\gamma\in \bG(K)} (\psi(g^{-1}\gamma h)\pm \e) \cdot \chi_{B_T} (g^{-1}\gamma h) .$$

We claim
that for any $T>1$ and for any $g, h\in W:=W_\infty\times W_f$,
\begin{equation}\label{i}
F_{(1-\e)T}^{-} (g, h)\le  F_T(e,e) \le F_{(1+\e)T}^{+} (g,h).\end{equation}
To see this, observe that
if $\gamma\in \bG(K)$ with $\He_\iota(\gamma)<T$, and $g, h\in W$
then  $$\He_{\iota}(g^{-1}\gamma h)\le (1+\e)T
\quad\text{and}\quad f(\gamma)\le \psi(g^{-1}\gamma h)+\e .$$
Hence $$F_T(e,e)=\sum_{\gamma\in \bG(K), \He_\iota(\gamma)<T}\psi(\gamma)
\le\sum_{\gamma\in \bG(K), \He_\iota (g^{-1}\gamma h) <(1+\e) T} (\psi(g^{-1}\gamma h)+\e)
=F_{(1+\e)T}^+(g, h),$$
proving the right inequality in (\ref{i}). The other inequality can be proved
similarly.
 Now let $\phi\in
C_c(Y_{W_f})$ be a non-negative $W_f$-invariant function such that
$\text{supp} (\phi)\subset \bG(K)\backslash \bG(K) W$ and
$\int_{Y_{W_f}}\phi \, d\tau_{W_f}=1$. By integrating (\ref{i}) over
$Y_{W_f}\times Y_{W_f}$ against the function  $\alpha(x, y)=\phi(x)\cdot
\phi(y)$, we obtain $$\langle F_{(1-\e)T}^{-}, \alpha \rangle \le
F_T(e,e)\le \langle F_{(1+\e)T}^{+}, \alpha \rangle .$$

Note that Theorem \ref{vol} implies the following:
there exist sequences $\{a_\e\ge 1\} $ and $\{b_\e\le 1\}$
such that  $a_\e \to 1$ and $b_\e\to 1$ as $\e\to 0$
 for all sufficiently small $\e>0$,
\begin{equation}\label{wr}
b_\e\le \liminf_T\frac{\tau_{W_f}(B_{(1-\e)T}\cap G_{W_f}))}
{\tau_{W_f}(B_T\cap G_{W_f})}
\le \limsup_T
\frac{\tau_{W_f}(B_{(1+\e)T}\cap G_{W_f})}
{\tau_{W_f}(B_T\cap G_{W_f})}\le a_\e.\end{equation}

Hence by applying Proposition \ref{wc},
\begin{align*}
&\limsup_T \frac{F_T(e,e)} {\tau_{W_f}(B_T\cap G_{W_f}
))}\le \limsup_T \frac{\langle F_{(1+\e)T} ^{+},\alpha\rangle
}{\tau_{W_f} (B_T\cap G_{W_f}) }
\\&\le  \limsup_T \frac{\langle F_{(1+\e)T} ^{+},\alpha\rangle }
{\tau_{W_f} (B_{(1+\e) T}\cap G_{W_f})} \cdot
\limsup_T
\frac{\tau_{W_f}(B_{(1+\e)T}\cap G_{W_f})}{\tau_{W_f}
(B_{T}\cap G_{W_f})} \\
&\le a_\e \cdot {\int_{X_\iota}
(\psi+\e)d\mu_{\iota,W_f} } \le a_\e \cdot
 (\int_{X_\iota}
\psi d\mu_{\iota,W_f} +\e)
\end{align*}
 and similarly,
$$
 b_\e \cdot  (\int_{X_\iota}
\psi d\tilde\mu_\iota -\e)
 \le \liminf_T
\frac{F_T(e,e)}{\tau_{W_f}(B_T\cap G_{W_f})}. $$
Taking $\e\to 0$, 
$$\lim _T \frac{F_T(e,e)}{\tau_{W_f} (B_T\cap G_{W_f})}=
 \int_{X_\iota}
\psi \; d\mu_{\iota,W_f} .$$
 This finishes the proof of Theorem \ref{smuiota}.
\end{proof}

 In fact, for the case $\psi=1$,
the computation in the proof of Theorem \ref{smuiota}
can be simplified significantly, and it applies to general
families of balls $B_T$, which we presently introduce.

For an increasing sequence $\{B_T\}$ of relatively compact subsets
of $\bga$ and a compact open subgroup $W_f\subset \bG(\A_f)$,
we call $\{B_T\}$ {\it $W_f$-well rounded}
if the following holds:
\begin{enumerate}
\item $W_f B_T W_f=B_T$ for any $T>1$;
\item for any small $\e>0$,
there exists a neighborhood $W_\e\subset \bG_\infty ^\circ$ of $e$ such that
$$W_\e B_TW_\e\subset B_{(1+\e)T} \quad\text{and}\quad B_{(1-\e)T}\subset
\bigcap_{g, h\in W_\e}gB_T h$$
for all $T>1$;
 \item
 $\tau_{W_f}(B_T\cap G_{W_f})\to \infty$ as $T\to \infty$ and
there exist constants $a_\e\ge 1 $ and $b_\e\le 1$
tending to $1$ as $\e\to 0$
such that for all sufficiently small $\e>0$,
\begin{equation*}
b_\e\le \liminf_T\frac{\tau_{W_f}(B_{(1-\e)T}\cap G_{W_f})}
{\tau_{W_f}(B_T\cap G_{W_f})}
\le \limsup_T
\frac{\tau_{W_f}(B_{(1+\e)T}\cap G_{W_f})}
{\tau_{W_f}(B_T\cap G_{W_f})}\le a_\e.\end{equation*}
\end{enumerate}
The proof of Theorem \ref{smuiota} gives
\begin{Prop}\label{p:asympt}
Let $\bG$ be a connected absolutely almost simple $K$-group, and
 let $W_f$ be a compact open subgroup
of $\bG(\A_f)$.
Then for any $W_f$-well rounded sequence $\{B_T\}$ of relatively compact subsets
of $\bga$,
  $$ \# \bG(K)\cap B_T  \sim_{T\to \infty} \tau_{W_f} (B_T\cap G_{W_f}) .$$
\end{Prop}

\section{Arithmetic fibrations and Construction of $\mu_\iota$}\label{arith}
In this section we prove the main theorems for a
general case, that is, without the saturation assumption on $\iota$.
We let $\bG$ be a connected semisimple
adjoint group defined over $K$ and $\iota$ a faithful
 representation of $\bG$ to $\GL_N$ defined over $K$ which has
a unique maximal weight $\lambda_\iota$.
Fix a height function $\He_\iota$ on $\bga$ associated to $\iota$. 

Let $\bT$, $\Phi(\bG, \bT)$, $\Delta$, $2\rho=\sum_{\alpha\in \Delta} u_\alpha \alpha$,
 $\lambda_\iota=\sum_{\alpha\in \Delta} m_\alpha \alpha$ and $a_\iota, b_\iota$ be as defined in Section \ref{in:a}.
 Let $\bM$ be the
smallest connected normal $K$-subgroup of $\bG$ whose root system
contains the set
$$
\{\alpha\in \Delta: \frac{u_\alpha+1}{m_\alpha}=a_\iota \}.
$$

Let $\bN$ be a connected normal $K$-subgroup of $\bG$ such that
$\bG=\bM \bN$ and $\bM\cap \bN=\{e\}$. Let $\pi: \bG\to \bN$ be
the canonical projection. Note that any element of $\bga$ can be
uniquely written as $g_1g_2$ with $g_1\in \bM(\A)$ and $g_2\in
\bN(\A)$. We denote by $\tau_{\bM}$ the Haar measure on $\bM(\A)$ such that
$\tau_{\bM}(\bM(\A)/\bM(K))=1$.

For each $x\in \bN(\A)$, the function $\He_\iota^x(g):= \He_\iota(gx)$ defines a
height function on $\bM(\A)$ with respect to $\iota':=\iota|_\bM$,
and $\iota'$ has a unique maximal weight by Lemma \ref{mn}.
Also the definition of $\bM$ implies that
$\iota'$
 is saturated, $a_{\iota'}=a_\iota$ and $b_{\iota'}=b_\iota$.
Set $$V_\iota:= \bM(\A_f)\cap W_\iota.$$
 Then
Theorems \ref{mt} and \ref{smuiota} (for the saturated cases) imply that for
 each $x\in \bN(K)$,
\begin{enumerate}
 \item there exists $c_x>0$ such that

\begin{align}\label{rates}
N_{\pi^{-1}(x)}( \He_\iota,T)=\#\{g\in \bM(K): \He_\iota(gx)< T\}\sim
c_x\cdot  T^{a_\iota}(\log T)^{b_\iota-1};
\end{align}
\item the following number $$
r_{x,\iota}=\sum_{\chi\in \Lambda}
\lim_{s\to a_\iota^+}(s-a_\iota)^{b_\iota}
\int_{\bM(\mathbb{A})} \He_{\iota}(gx)^{-s} \chi(g)\; d\tau_{\bM}(g)
$$
is a positive real number;
\item  if $X_\bM$ denotes the closed subspace
$\prod_{v\in R} \overline{\iota'(\bM(K_v))}$ of $X_\iota$, there exists a
probability measure $\tilde\mu_{x,\iota}$ on $X_{\bM}$
 such that for any $\psi\in C(X_\bM)$
which is invariant under a co-finite subgroup of $V_{\iota}$,
$$
\tilde\mu_{x,\iota}(\psi)=r_{x,\iota}^{-1}\cdot \sum_{\chi\in \Lambda}
\lim_{s\to a_\iota^+}(s-a_\iota)^{b_\iota}
\int_{\bM(\mathbb{A})} \He_{\iota}(gx)^{-s} \chi(g)\;\psi(g) d\tau_{\bM}(g).
$$
\end{enumerate}


Noting that $N(\He_\iota, T)=\sum_{x\in\bN(K)} N_{\pi^{-1}(x)} (\He_\iota, T)$,
we restate Theorem \ref{mt} in the introduction.
 \begin{Thm}\label{uns} We have 
\begin{enumerate} \item $c_{\He_\iota}:=\sum_{x\in\bN(K)} c_x <\infty$;
 \item 
for some $\delta>0$,
\begin{equation}\label{eq:naive}
N(\He_\iota,T)= c_{\He_\iota}\cdot T^{a_\iota}(\log T)^{b_\iota-1} (1+
O((\log T)^{-\delta})) .\end{equation} 
\end{enumerate}
\end{Thm}


Set for $T>0$ and $x\in \bN(\A)$,
$$B_T^x=\{g\in \bM(\A):\, \He_\iota^x(g) < T\} .$$
Since $x$ commutes with $\bM(\A)$, each height function ${\He_\iota^x}$ 
on $\bM(\A)$ is
invariant under $V_\iota$.
 Let
$Y_\bM=\bM(K)\backslash M_{V_\iota}$ and $\tau$ be the invariant
probability measure on $Y_\bM$. For each $x\in \bN(K)$, set
$$
F_T^x(g,h):=\sum_{\gamma\in \bM(K)} \chi_{B_T^x} (g^{-1}\gamma
h),\quad g,h\in \bM(\A).
$$
We may consider $F_T^x$ as a function on $Y_\bM\times Y_\bM$. Write
$\bM=\bM_1\cdots \bM_r$ as a product of connected $K$-simple
$K$-groups
and denote by $\tau_i$ the invariant probability measure
on $\bM_i(K)\backslash
\bM_i(\A)\cap M_{V_\iota}$.
 For a collection of smooth $(W_f\cap
\bM_i(\A))$-invariant functions $\psi_i\in C_c(\bM_i(K)\backslash
\bM_i(\A)\cap M_{V_\iota})$,
 such that $\int \psi_i\, d\tau_i =1$ for each 
 $1\le i\le r$, define
 $\psi\in C_c(Y_\bM)$ and $\alpha\in C_c(Y_\bM\times Y_\bM)$ by
  $$\psi(z_1, \cdots, z_r):=\prod_{i=1}^r \psi_i(z_i) \quad\text{and}\quad \alpha(y_1, y_2):= \psi(y_1)\psi(y_2) .$$
 
\begin{Lem}
\begin{enumerate}
 \item  There is a constant $c>0$ such that
for any $x\in \bN(K)$, 
$$c_x\le c\cdot\He_\iota(x)^{-a_\iota}.$$
\item
There exist $l\in \n$ and $\delta>0$, independent of
 $x$, such that 
 for any $x\in \bN(K)$ and $T\gg \He_\iota(x)$, $$ \langle
F_T^x,\alpha\rangle_{Y_\bM\times Y_\bM}=c_x\cdot T^{a_\iota}(\log
T)^{b_\iota-1}+ O( d_x \cdot C_\psi ' \cdot T^{a_\iota}(\log
T)^{b_\iota-1-\delta})
$$
where $d_x= \He_\iota(x)^{-a_\iota}
(\log \He_\iota(x))^{b_\iota-1}$ and $C_\psi'=\max (1, \max _i \|\Cal D^l \psi_i\|^{2r} )$
 and $\Cal D$ is the elliptic operator defined in
(\ref{elliptic}).
\end{enumerate}
\end{Lem}

\begin{proof} 

As in the proof of Proposition \ref{wc}, we derive that
\begin{align*}\left<F_T^x,\alpha \right>&=
\int_{g\in B_T^x\cap M_{V_\iota}}\left<\psi,g.\psi\right> \,
d\tau(g)\end{align*}

Note that \begin{align*} | \left<\psi,g.\psi\right>-1|&=
|\prod_{i=1}^r\left<\psi_i,g_i.\psi_i\right>-1|\\
 &=|\sum_{i=1}^r
(\prod_{j=1}^{i-1} \left<\psi_j,g_j.\psi_j\right>) (\langle \psi_i,
g_i .\psi_i \rangle -1)| \\
&\le r\cdot  C_{\psi} \cdot \max_i| \langle \psi_i, g_i .\psi_i
\rangle -1|\\
&=r\cdot  C_{\psi} \cdot \max_i |\langle \psi_i-1,
 g_i.(\psi_i-1)\rangle|
\end{align*}
where $C_{\psi}=\max ( 1, \max_i \|\psi_i\|^{2r-2})$.
 Since
$\psi_i-1\in L_{00}^2(\bM_i(K)\backslash \bM_i(\A)\cap
M_{V_\iota})$ for each $i$, we deduce from Theorem \ref{smooth}
 that
\begin{align}\label{ftx1}
\left|\left<F_T^x,\alpha \right>-\tau (B_T^x\cap
  M_{V_\iota})\right|
\le 2r\cdot (\prod_i c_{W_f\cap \bM_i(\A)}) \cdot C_\psi' \cdot
\int_{g=g_1\cdots g_r\in B_T^x\cap
  M_{V_\iota}} ( \max_{i}\; \tilde
\xi_{\bM_i}(g_i) ^{1/2}) \, d\tau(g)
\end{align}
where $C_\psi'=\max (1, \max _i \|\Cal D^l \psi_i\|^{2r} )$ for some
large $l$.

Since $\tilde \xi_{\bM_i}\le {\xi_{\bM_i}}^{1/2}$,  it follows from
Lemma \ref{xhi} that there exist $m\in \n$ and $C_1>0$ such that for
any $1\le i\le r$,
$$
\tilde \xi_{\bM_i}^{1/2}(g_i) < C_1\cdot  \He_\iota(g_i)^{-1/m}
\quad\text{for any $g_i\in \bM_i(\A)$.}
$$

Define a function on $\M (\A)$ by
$$\tilde \He(g_1\cdots g_r):=\min_i \He_\iota(g_i),\quad g_i\in\bM_i(\A) .$$

Let $\kappa$ be as in Lemma \ref{l:prod} for $\bG_1=\bM$ and
$\bG_2=\bN$ so that
 $B_T^x\subset B_{\kappa T \cdot \He_\iota(x)^{-1}}$.
 It then follows from (\ref{ftx1}) that
\begin{align}\label{ftx2}
\left|\left<F_T^x,\alpha \right>-\tau(B_T^x\cap
  M_{V_\iota})\right|
& < C_2 \cdot C_\psi' \cdot \int_{B_{\kappa T\cdot \He_\iota(x)
^{-1}}\cap
  M_{V_\iota}}  \tilde \He(g)^{-1/m}\, d\tau(g)
\end{align}
for a constant $C_2>0$ independent of $x$.

 Since $\iota'$ is
saturated, by Lemma \ref{tfae},
for every proper normal $K$-subgroup $\bf L$ of $\bM$,
$$
\tau_{\bf L}(B_T\cap {\bf L}(\A))\ll (\log T)^{-1}
\tau(B_T\cap M_{V_\iota})
$$
where $\tau_{\bf L}$ is a Haar measure on $\bf{L}(\A)$.

For each $C>1$, set
$$B^C=\{g\in \bM(\A):\, \tilde \He(g)>C\} .$$

Note that
$$(B_T-B^C)\cap M_{V_\iota}
\subset \cup_{i=1}^r \Omega_i$$
where $\Omega_i=\{g=g_1\cdots g_r\in M_{V_\iota}:
 \He_\iota(g_i)\le C, \;\;\He_\iota(g)<T \}$.
Now denoting by $\bL^{(i)}$ the subgroup of $\bM$ generated by
$\bM_1, \cdots, \bM_{i-1}, \bM_{i+1},\cdots, \bM_r$, let
$\kappa_i>1$ be a constant as in Lemma \ref{l:prod} for
$\bG_1=\bM_i$ and $\bG_2=\bL^{(i)}$. 
Let $\delta_0:=\inf_{g\in \bga}\He_\iota(g)>0$ (Lemma \ref{ff3}).
Then for any $C\gg 1$,
 \begin{align*} \tau(\Omega_i) &\le
\int_{\He_\iota(g_i)<C} \tau_{\bL_i(\A)}(B_{\kappa_i \delta_0^{-1} T}\cap
\bL^{(i)}_{V_\iota\cap \bL(\A_f)}) \,d\tau_{\bM_i} (g_i)\\
&\ll
 C^{a_{\iota}}(\log
C)^{b_\iota-1}(\log T)^{-1} \tau(B_{\kappa_0 T}\cap M_{V_\iota})
\end{align*}
where $\kappa_0=\max_i(\kappa_i \delta_0^{-1})$.

Hence for any $C\gg 1$ and $T\gg C$,
$$\tau ((B_T-B^C)\cap M_{V_\iota})
\ll C^{a_{\iota}}(\log
C)^{b_\iota-1}(\log T)^{-1} \tau(B_{\kappa_0 T}\cap M_{V_\iota}).$$
Therefore
\begin{align}\label{ftx3}
 \int_{B_{T}\cap  M_{V_\iota} } \tilde \He ^{-1/m}\, d\tau &=
 \int_{B_{T}\cap B^C\cap
  M_{V_\iota}} \tilde \He^{-1/m}\, d\tau+ \int_{(B_{T}- B^C)\cap
  M_{V_\iota}}  \tilde \He^{-1/m}\, d\tau\\
&\ll (C^{-1/m}+ \delta_0^{-1/m} \cdot C^{a_{\iota}}(\log
C)^{b_\iota-1}(\log T)^{-1})\cdot
\tau(B_{\kappa_0 T}\cap M_{V_\iota}) \notag \\
& \ll
 (\log T)^{-\delta}\cdot \tau(B_{\kappa_0 T} \cap
M_{V_\iota}) \quad\text{ for $C=(\log T)^{1/(2 a_\iota)}$ }
\notag\end{align} for some $\delta>0$.
 We now deduce from
(\ref{ftx2}) and (\ref{ftx3}) that
\begin{equation}\label{ftx33}
\left<F_T^x,\alpha \right>=\tau(B_T^x\cap M_{V_\iota})+O\left(
C_\psi'\cdot  (\log T)^{-\delta}\cdot \tau(B_{\kappa_0\kappa T\cdot
\He_\iota(x)^{-1}}\cap M_{V_\iota})\right)
\end{equation}
for some $\delta>0$.
Let  $S\subset R$ be as in the proof of Lemma \ref{l:prod}, that is,
for any
 $v\in R-S$,
$$\bG(K_v)=U_v A^+_v U_v\quad\text{and}\quad
\He_{v}(\iota(g))=\chi(a)\quad\hbox{for $g=u_1 au_2\in \bG(K_v)$}.
$$

Denote by $\tau_S$ and $ \tau^S$
 Haar measures on $\bM_S$ and $\bM^S$ respectively
such that $\tau=\tau_S\times \tau^S$ locally.

Recall from \eqref{eq:haar}
 the map $\bM_S \to (M_{V_\iota}\cap \bM^S )\ba \bM^S$ by $g\mapsto [s_g]$
 and as in \eqref{eq:dom} we have
\begin{align}\label{ftx4}
\tau(B^x_T\cap  M_{V_\iota}) &=
 \int_{g\in \bM_{S}} \tau^S(B_{\kappa T \cdot
\He_{\iota}^{-1}(gx)}\cap s_g M_{V_\iota}\cap \bM^S)\,d\tau_{S}(g).
\\
&= \int_{g\in B_{\delta_0^{-1} \kappa
T\cdot \He_\iota^{-1}(x)}\cap \bM_{S}} \tau^S(B_{\kappa T \cdot
\He_{\iota}^{-1}(gx)}\cap s_g M_{V_\iota}\cap \bM^S)\,d\tau_{S}(g).
\notag\end{align}
By Theorem \ref{vol}, there is $c_0>0$ such that
for all $g\in \bM_S$, $$
\tau^S(B_{T}\cap  s_g M_{V_\iota}\cap \bM^S) =c_0\cdot
 \gamma_{V_\iota,S}(s_g^{-1})\cdot 
T^{a_\iota}(\log T)^{b_\iota-1} +O(T^{a_\iota}(\log T)^{b_\iota-2})
.$$
Here the implied constant can be taken uniformly for all
$g\in \bG_S$, since there are only finitely many cosets
$s_gM_{V_\iota } \cap \bM^S$.

 Note that
$\gamma_{V_\iota}^S(s_g^{-1})=\gamma_{V_\iota}^S(g)=\gamma_{V_\iota}^S(gx)$
since $x\in \bN(K)$ and it is bounded.
We deduce that when $\He_\iota(gx)\ll T/\delta_0$,
\begin{align*}
&\tau^S(B_{\kappa T\He_{\iota}^{-1}(gx)}\cap  M_{V_\iota}\cap \bM^S)\\
=&c\cdot \gamma_{V_\iota}^S(g) (T\cdot \He_\iota^{-1}(gx))^{a_\iota}(\log T)^{b_\iota-1}
+O((T\cdot \He_\iota^{-1} (gx))^{a_\iota}(\log
\He_\iota(gx))^{b_\iota-1}(\log T)^{b_\iota-2}).
\end{align*}
for $c=c(S, V_\iota, \kappa)>0$.
To estimate the integral over the domain $\He_\iota(gx)\gg
T/\delta_0$, it suffices to note that by Lemma \ref{mv},
$$
\tau_S(B_{ T\cdot \He_\iota^{-1}(x)}\cap \bM_{S})\ll (T \He_\iota^{-1}
(x))^{a_\iota-\epsilon}.
$$
Since by Lemmas \ref{l:prod} and \ref{mv},
\begin{align*}
\int_{g\in \bM_{S}} \gamma_{V_\iota}^S(g) \He_\iota(gx)^{-a_\iota} (\log
\He_\iota(gx))^{b_\iota-1}\, d\tau_{S}(g)\ll \He_\iota(x)^{-a_\iota}
(\log
\He_\iota(x))^{b_\iota-1},
\end{align*}
it follows from the above estimates that for $T\gg \He_\iota(x)$,
$$\tau(B^x_T\cap M_{V_\iota})= c_x T^{a_\iota}(\log T)^{b_\iota-1}+
O( d_x T^{a_\iota}(\log T)^{b_\iota-2}),
$$
where
\begin{align}\label{cx}
c_x=c \cdot\int_{g\in \bM_{S}}\gamma_{V_\iota}^S(gx)
\He_\iota(gx)^{-a_\iota}\,d\tau_{S}(g)\ll \He_\iota(x)^{-a_\iota} .
\end{align}
Hence combining (\ref{ftx33}) and (\ref{ftx4}), we have for $T\gg \He_\iota(x)$,
$$ \left<F_T^x,\alpha\right>=c_x \cdot T^{a_\iota}(\log T)^{b_\iota-1}+
O( d_x \cdot C_\psi' \cdot  T^{a_\iota}(\log T)^{b_\iota-1-\delta} ). $$

\end{proof}

A key ingredient in deducing Theorem \ref{uns}
is the following stronger version of (\ref{rates}):

\begin{Prop}\label{p:rate}
There exists $\delta>0$ such that for each $x\in \bN(K)$
and for any $T \gg \He_\iota(x)$,
\begin{equation}\label{eq:errr}
N_{\pi^{-1}(x)}(\He_\iota,T)=c_x\cdot  T^{a_\iota}(\log T)^{b_\iota-1}
+O( d_x\cdot T^{a_\iota}(\log T)^{b_\iota-1-\delta})
\end{equation}
where $d_x= \He_\iota(x)^{-a_\iota}
(\log \He_\iota(x))^{b_\iota-1}$ and the implied constant is independent of $x$. 
\end{Prop}
\begin{proof}
 Let $\phi_\epsilon$ be a smooth symmetric nonnegative
function on $\bM_{\infty}$, which is a product $\prod_{i=1}^r
\phi_{i, \e}$ of smooth functions on the simple factors of
$\bM_\infty$, $\int_{\bM_{\infty}} \phi_\epsilon\,  d\tau_\infty=1$
and $\hbox{supp}(\phi_\epsilon)$ is contained in the Riemannian ball
at $e$ in $\bM_{\infty}$ of radius $\epsilon$, and for some
$\rho>0$, $\max_i \|\Cal D^{l}\phi_{i, \e}\|^{2r} \ll \e^{-\rho}$
(see, for example, Lemma 4.4 in \cite{GaO}). By the definition of $\He_\iota$ in
\eqref{hed}, there exists $b>0$ such that
$$
\hbox{supp}(\phi_\epsilon)\cdot B_T^x\cdot\hbox{supp}(\phi_\epsilon)\subset B_{(1+b\epsilon)T}^x
$$
for every $T>1$ and $x\in \bN(K)$.

Define
$$
\psi_\epsilon(g)=\frac{1}{\tau^{R_f}(V_\iota)}
\sum_{\gamma\in \bM(K)} \phi_\epsilon(\gamma g_\infty)
\cdot \chi_{V_\iota}(\gamma g_f),\quad g=g_\infty g_f\in \bM_\infty\bM(\A_f).
$$

Define $\alpha_\e(y_1,
y_2)=\psi_\e(y_1)\psi_\e(y_2)$ for $(y_1, y_2)\in Y_\bM\times Y_\bM$. Then
\begin{align*}
&N_{\pi^{-1}(x)}(\He_\iota,T)\le
\left<F_{(1+b\epsilon)T}^x,\alpha_\epsilon\right>\\
=&\,c_x T^{a_\iota}(\log T)^{b_\iota-1}+ O(c_x \cdot \epsilon\cdot
T^{a_\iota}(\log T)^{b_\iota-1}+  d_x\cdot \epsilon^{-
\rho}T^{a_\iota}(\log T)^{b_\iota-1-\delta}).
\end{align*}
Setting $\epsilon=(\log T)^{-\delta/(\rho+1)}$, we derive the upper
estimate for $N_{\pi^{-1}(x)}(\He_\iota,T)$. The lower estimate is
proved similarly.
\end{proof}

\begin{proof}[\bf Proof of Theorem \ref{uns}]
According to the choice of $\bN$, for any simple root $\alpha\in \Delta$
whose restriction to $\bN$ is a root, we have
\begin{equation}\label{nsm}
\frac{u_\alpha+1}{m_\alpha} <a_\iota .\end{equation}
Since $\bN(K)$ is a discrete subgroup of $\bN(\A)$,
we can find an open neighborhood $U:=U_\infty \times U_f$ of the identity
in $\bN(\A)$ such that $\gamma U\cap \gamma' U=\emptyset $ for all $\gamma\ne \gamma' \in \bN(K)$.
We may assume $U_f\subset W_\iota$ and $ B_TU_\infty \subset B_{2 T}$ for all $T\gg 1$.
Since $\tau_\bN(\gamma U)=\tau_\bN(U)$ by the invariance of $\tau_\bN$,
we deduce
\begin{equation*}
N_\bN(\He_\iota, T) =\tau_\bN(U)^{-1} \cdot \tau_\bN \left(\bigcup_{\gamma\in \bN(K): \He_\iota(x) \le T} \gamma U\right)
\le \tau_\bN(U)^{-1} \cdot \tau_\bN(B_{2T}\cap \bN(\A)).
\end{equation*}
Therefore Theorem \ref{vol}, applied to $\tau_\bN(B_{2T}\cap \bN(\A))$, 
 together with \eqref{nsm} yields that there exists $\epsilon>0$ such that $$
N_\bN(\He_\iota,T)=O(T^{a_\iota-\epsilon}).$$

Hence setting $\alpha(t)=N_{\bN}(\He_\iota,t)$, we have for any $a_\iota -\e/2 \le a \le a_\iota $,
\begin{equation*}
\sum_{x\in \bN(K)} \He_\iota(x)^{-a}=\int_0^\infty t^{-a}\,
d\alpha(t)
=a \int_0^\infty t^{-a-1}\alpha(t)\, dt \ll
\int_0^\infty t^{-(1+ \epsilon/2)}\, dt
< \infty.
\end{equation*}

Since $c_x\ll \He_\iota(x)^{-a_\iota}$ by (\ref{cx}) and
 $d_x= \He_\iota(x)^{-a_\iota} (\log
 \He_\iota(x))^{b_\iota-1}$,
 it follows  that
\begin{equation}\label{ch}
c_{\He_\iota}:=\sum_{x\in\bN(K)} c_x<\infty\quad\text{and}\quad
\sum_{x\in\bN(K)} d_x <\infty.
\end{equation}

Let $\delta_0>0$ be as in (\ref{eq_delta0}) and let $\beta>0$ be 
such that 
Proposition \ref{p:rate} holds
for all $T>\beta \cdot \He_\iota(x)$. Let $\delta>0$ be a constant
given in the same proposition.

Applying Lemma \ref{l:prod} for $\bM$ and $\bN$
with $\kappa$ therein, we have
\begin{align}\label{eq:est}
&\sum_{x\in\bN(K): \He_\iota(x)>
 \beta ^{-1}T} N_{\pi^{-1}(x)}(\He_\iota,T)
\\&=\#\{xy\in \bN(K)\bM(K):  \He_\iota(x)>\beta^{-1} T, \;\;
 \He_\iota(xy)<T\} \nonumber
\\
&\le  N_{\bM}(\He_\iota,\kappa \beta)
\cdot N_{\bN}(\He_\iota,\kappa T\delta_0^{-1})  \nonumber
\\
& =O(T^{a_\iota-\epsilon}).  \nonumber
\end{align}
Now applying Proposition \ref{p:rate},
since $\sum_{x\in\bN(K)} d_x<\infty$,
\begin{align*}
&\sum_{x\in\bN(K):
 \He_\iota(x)\le \beta ^{-1}T} N_{\pi^{-1}(x)}(\He_\iota,T)\\
&=\left(\sum_{x\in\bN(K): \He_\iota(x)\le \beta^{-1} T} c_x\right)T^{a_\iota}(\log T)^{b_\iota-1}+ O(T^{a_\iota}(\log T)^{b_\iota-1-\delta})
.\end{align*}
Therefore as $T\to \infty$,
\begin{align*} N(\He_\iota, T)&=
\sum_{x\in\bN(K):
 \He_\iota(x)\le \beta ^{-1}T} N_{\pi^{-1}(x)}(\He_\iota,T) + O(T^{a_\iota-\e})  \\
&= \left(\sum_{x\in\bN(K): \He_\iota(x)\le \beta^{-1} T} c_x\right)
 T^{a_\iota}(\log T)^{b_\iota-1}(1+ O((\log T)^{-\delta})) . \end{align*}

Since $\sum_{x\in\bN(K): \He_\iota(x)\le \beta^{-1} T} c_x= C(\He_\iota)
+O(T^{ -\e})$, we have
$$ N(\He_\iota, T)=
 C(\He_\iota)\cdot T^{a_\iota}(\log T)^{b_\iota-1}(1+ O((\log T)^{-\delta}))$$
finishing the proof. \end{proof}

We now construct the probability measure on $X_\iota$ in order to
prove Theorem \ref{e} in the introduction in a general case.
We consider each $\tilde\mu_{x,\iota}$ as a measure on $X_\iota$ since $X_\bM$
is a closed subspace of $X_\iota$. For each $x\in \bN(K)$,
 denote by $x.\tilde \mu_{x, \iota}$ the measure defined by
$$(x.\tilde\mu_{x, \iota})(\psi):=\tilde\mu_{x, \iota}(\psi_x)$$
where $\psi_x(g)=\psi(gx)$.
Noting that each $x.\tilde\mu_{x, \iota}$ supported on $X_\bM x$,
defines a probability measure  $\mu_\iota$ on $X_{\iota}$ by
\begin{equation}\label{arith-measure}
\mu_\iota=\sum_{x\in\bN(K)} \frac{c_x}{c_{\He_\iota}} (x.\tilde\mu_{x,\iota}).
\end{equation}

Note that $\mu_\iota=\tilde\mu_\iota$
 in the case when $\iota$ is saturated (see Theorem \ref{smuiota-g} for the definition of
$\tilde\mu_\iota$).
\begin{Thm}\label{t:eqgeneral}
For any $\psi\in C(X_\iota)$,
$$
\lim_{T\to \infty}\frac{1}{N(\He_\iota,T)}
 \sum_{g\in\bG(K): \He_\iota(g)<T} \psi(g) = \int_{X_\iota} \psi\,
 d\mu_{\iota} .$$
\end{Thm}
\begin{proof}
Let $\psi\in C(X_\iota)$. We write
\begin{align*}
\mu_T (\psi) &=\frac{1}{N(\He_\iota,T)} \sum_{g\in \bG(K):\, \He_\iota(g)<T}
\psi(g),\\
\mu_{x,T}(\psi) &=\frac{1}{N_{\pi^{-1}(x)}(\He_\iota,T)} \sum_{g\in \bM(K):\, \He_\iota(gx)<T}
\psi({g})\quad\hbox{for each}\quad  x\in \bN(K) .
\end{align*}

Note that
\begin{equation}\label{eq:est2}
\mu_T (\psi)=\sum_{x\in\bN(K)} \frac{N_{\pi^{-1}(x)}(\He_\iota,T)}
{N(\He_\iota,T)}\mu_{x,T}(\psi).
\end{equation}
Since $\{\mu_T: T\gg 1\}$ is a sequence of probability measures on a compact space $X_\iota$,
it suffices to prove that any convergent subsequence of $\mu_T$, in the weak$^*$ topology,
has the same limit $\mu_\iota$. Hence without loss of generality, we may assume
that $\mu_T$ is convergent.

Let $\delta>0$ be as in Proposition \ref{p:rate}
and let $\beta>0$ be such that the same proposition holds for
all $T>\beta\cdot \He_\iota(x)$.
By (\ref{eq:est}) and Theorem \ref{uns},
there exists $\e>0$ such that $$
\sum_{x\in\bN(K):\,\He_\iota(x)>
 \beta ^{-1}T} \frac{N_{\pi^{-1}(x)}(\He_\iota,T)}{N(\He_\iota,T)}
\le O(T^{-\epsilon}).
$$

Hence
\begin{equation}\label{major}
\lim_{T\to\infty}\mu_T (\psi)=\lim_{T\to\infty} \sum_{x\in\bN(K):\,\He_\iota(x)\le
 \beta ^{-1}T}  \frac{N_{\pi^{-1}(x)}(\He_\iota,T)}
{N(\He_\iota,T)}\mu_{x,T}(\psi).
\end{equation}

For $x\in\bN(K)$ such that $\He_\iota(x)\le \beta
^{-1}T$, we deduce from Proposition \ref{p:rate} and Theorem \ref{uns}
that 
$$
\frac{N_{\pi^{-1}(x)}(\He_\iota,T)}{N(\He_\iota,T)}\le \frac{c_x}{c_{\He_\iota}} +\theta_x
(\log T)^{-\delta},
$$
where 
$\theta_x$ comes from the error terms in (\ref{eq:naive})
and (\ref{eq:errr}), and satisfies $\sum_{x\in \bN(K)}
\theta_x<\infty$. 
Since the sum in \eqref{major} is majorized by 
$$
\sum_{x\in \bN(K)} \left(\frac{c_x}{c_{\He_\iota}}
 +\theta_x (\log T)^{-\delta}\right)<\infty,
$$
we can apply the dominated convergence theorem to obtain
$$
\lim_{T\to\infty}\mu_T(\psi)=
\sum_{x\in \bN(K)}\frac{c_x}{c_{\He_\iota}}\left(\lim_{T\to\infty}\mu_{x,T}(\psi)\right).
$$
By Theorem \ref{e}, applied to $\bM$ and the height function
 $g\mapsto \He_\iota(gx)$, we have
 for every $\psi\in C(X_\iota)$, 
 $$
\lim_{T\to\infty} \mu_{x,T}(\psi)=\frac{1}{N_{\pi^{-1}(x)}(\He_\iota,T)}
 \sum_{g\in\bM(K): \He_\iota(gx)<T} \psi_x (g)=
 \int_{X_\iota} \psi_x \,d\tilde\mu_{x,\iota} 
$$
where
$\psi_x(g)=\psi(gx)$.

Therefore 
$$
\lim_{T\to\infty} \mu_T(\psi)=\sum_{x\in\bN(K)} \frac{c_x}{c_{\He_\iota}}\cdot \tilde \mu_{x,
  \iota}(\psi_x).
$$
\end{proof}
\section{Manin's and Peyre's conjectures}\label{manin}


In this section, we now explain our main results in the context of
Manin's conjecture on the asymptotic number of rational points of bounded height
for Fano varieties. Let $X$ be a smooth projective variety
defined over $K$. For every line bundle class $[L]$ 
on $X$ defined
over $K$, there exists an
associated height function $\He_{\Cal L}$
 on $X(K)$, unique up to the multiplication by
 bounded functions, via Weil's height machine (cf. \cite[Theorem B. 3.2]{Si}).
 For example, if $L$ is a very ample line bundle of $X$ with a $K$-embedding
 $\psi_L:X\to\mathbb P^N$, then a height function $\He_\Cal L$ on $X(K)$
is defined as $$\He_\Cal L:=\He \circ \psi_L $$ 
 for some height function $\He$ on $\mathbb P^N(K)$. 
We call a pair $\Cal L=(L, \He_\Cal L)$ a metrized line bundle. 
Due to
the freedom of choosing a height function $\He$ on $\mathbb P^N(K)$,
$\He_\Cal L$ is not uniquely determined and this is why we use the
subscript $\Cal L$ rather than $L$.

 For a metrized
ample line bundle $\Cal L=(L, \He_\Cal L)$ on $X$ and a subset $U$ of $X$, set
$$ N_{U}(\Cal L ,T):=\# \{g\in U\cap X(K)\;:\; \He_{\Cal L}(g)\,< \,T\}. $$
The goal of Manin's conjecture is to
obtain the asymptotic (as $T\to\infty $) of
 $N_{U}(\Cal L ,T)$ for some Zariski open subset $U$ of $X$, possibly by
passing to a finite extension field of $K$.

Two important geometric invariants here are:
\begin{align*}
a_L&:=\inf \{a\in \q^+: a[L]+[K_X]\in \Lambda_{\op{eff}}(X)\}
\hbox{ --- the Nevanlinna invariant of $L$},\\
b_L&:= \,\hbox{the codimension of the face of
$\Lambda_{\op{eff}}(X)$ containing $a_L[L]+[K_X]$ in its interior}
\end{align*}
where $[K_X]$ denotes the canonical line bundle class and
$\Lambda_{\op{eff}}(X)$ denotes the cone of classes of effective
 line bundles on $X$.

Now let $\bG$ be
a connected semisimple adjoint algebraic group defined over $K$.
 Let $X$ denote the projective $K$-variety, which
is the wonderful compactification of $\bG$ constructed by
 De Concini and Procesi \cite{DP} and by De Concini and Springer \cite{DS}.
It is shown in \cite{DP} that $X$ is a Fano variety.

One way of constructing $X$ explicitly
is to take the Zariski closure of the
image of $\bG$ in $\mathbb P(\op{M}_N)$ under an irreducible
faithful representation $\bG \to \GL_N$ whose highest weight is
regular.  A dominant weight
$\chi$ is called regular if $\chi=\sum_{\alpha\in \Delta} m_\alpha \omega_\alpha$ with all
$m_\alpha>0$ where $\{\omega_{\alpha}:\alpha\in \Delta\}$ is the set
of fundamental weights.

The Picard group $\op{Pic}(X)_{\bar K}$ is isomorphic to the
weight lattice of $\bG$. Under this isomorphism, the simple roots $\alpha$
correspond to the boundary divisors $D_\alpha$ such that
$X-\bG=\cup_\alpha D_\alpha$, and the Galois action on
$\op{Pic}(X)_{\bar K}$ corresponds to the twisted Galois action
(also called the $*$-action, see \cite[2.3]{Ti1}) on the
weight lattice. Hence, the Picard group
$\op{Pic}(X)$ is freely generated by
the line bundles corresponding to
the orbits of the fundamental weights under the twisted Galois
action. The closed
cone $\Lambda_{\op{eff}}(X)$ of the effective line bundles
is the positive cone spanned by $D_{\Gamma_K.\alpha}$,  $\alpha\in \Delta$,
i.e., $$\Lambda_{\op{eff}}(X)=\oplus\, \br_{\ge 0}
\,[ D_{\Gamma_K.\alpha}] $$
where the sum is  taken over the $\G_K$-orbits $\G_K.\alpha$
 in
the set $\{\alpha\in \Delta\}$ of simple roots and
$D_{\Gamma_K .\alpha}:=\sum_{{\beta}\in \Gamma_K.\alpha} D_{\beta}$,
and the anticanonical class $[-K_X]$ corresponds to
$2\rho+\sum_{\alpha\in \Delta}\alpha$. Moreover any ample line
bundle class $[L]$ of $X$ over $K$ corresponds to a regular
 dominant weight in such a way that
if $[L']:=m[L]$ corresponds to
$\chi\in X^*( \mathbf T)$ for $m\in \n$,
 the restriction of
 $\He_{\Cal L'}$ to $\bG(K)$ coincides with
 a height function $\He_\iota$ with
 respect to the irreducible representation
 $\iota$ defined over $K$ with the highest weight $\lambda_\iota$
and $a_L=a_\iota$ and $b_L=b_\iota$ (cf. \cite[Proposition 6.3]{STT2}).
In particular, $\He_{\Cal L'}$
has a natural extension to $\bga$.
 We refer to \cite[Ch.6]{BK} for a more detailed
account on the wonderful compactification,and
\cite[4.1]{T} and \cite[section 6]{STT2}
 on metrized line bundles.

Therefore the following theorem, conjectued by
Manin, follows from Theorem \ref{mt}.
\begin{Thm}\label{manp} Let $X$ be the wonderful compactification of a connected adjoint semisimple $K$-group
$\bG$, and $\Cal
L=(L, \He_\Cal L)$ a metrized ample line bundle on $X$.
 Then  there exist $c_{\Cal L}>0$ and $\delta>0$ such that
$$ N_\bG (\Cal L ,T)
= c_{\Cal L} \cdot T^{a_L} (\log T)^{b_L-1}(1+O((\log T) ^{-\delta})).$$
\end{Thm}



\vs
In order to describe the distribution of rational points,
we construct a finite measure $\tau_{\mathcal L}$ on $X(\A)$
first for each saturated ample line bundle $\mathcal L$ and then for
any ample line bundle $\mathcal L$.

\begin{Def}  We call an ample line bundle
$L$ on $X$ {\it saturated} if the representation defined by the
 corresponding dominant weight is saturated.
\end{Def}

We note that
 if $\bG$ is $K$-simple, every ample line bundle is saturated,
and that the
 anticanonical
line bundle $-K_X$ is always saturated for any $\bG$.
 Batyrev and Tschinkel introduced the notion of a strongly
saturated line bundle in \cite{BT3}: A line bundle
$\mathcal{L}$ is
called {\it strongly saturated} if 
for any Zariski open dense subset $U$
of $X$, 
\begin{equation}\label{sst}
\lim_{T\to\infty} \frac{ N_{U}(\Cal L ,T)}{ N_{\bG}(\Cal L ,T)}=1.
\end{equation}

\begin{Lem}
A strongly saturated (see \eqref{sst}) ample line bundle $L$ is saturated.
\end{Lem}
\begin{proof}
Suppose not. Then by Theorem
\ref{vol}, there exists a connected normal $K$-subgroup
$\bM$ of $\bG$ such that $\iota|_{\bM}$ is saturated and
the volume of $B_T\cap \bM(\A)$ is
of order $T^{a_L}(\log T)^{b_L-1}$.
By Theorem \ref{e}, $\# B_T\cap \bM(K)$ has the order
of $T^{a_L}(\log T)^{b_L-1}$. This contradicts to the assumption
that $L$ is strongly saturated.
\end{proof}




\begin{Lem}\label{p:tam} Let $\mathcal{L}$ and $\mathcal{L}'$ be
 metrizations of a
saturated line bundle $L$ and $(B_T, W_f)$, and $(B_T',W_f')$ be defined as
above with respect to $\mathcal L$ and $\mathcal L'$ respectively. Then
$$
\lim_{T\to\infty}\frac{\tau_{W_f}(B_T\cap
  G_{W_f})}{\tau_{W_f'}(B_T'\cap
  G_{W_f'})}=
\frac{\tau_\mathcal{L}(\bG(\A))}{\tau_{\mathcal{L}'}(\bG(\A))}.
$$
\end{Lem}

\begin{proof}
Let $V_f=W_f\cap W_f'$. By Proposition \ref{tauf}, it suffices to
show that
\begin{equation}\label{eq:scale}
\lim_{T\to\infty} \frac{\tau(B_T\cap
  G_{V_f})}{\tau(B_T'\cap
  G_{V_f})}=
\frac{\tau_\mathcal{L}(\bG(\A))}{\tau_{\mathcal{L}'}(\bG(\A))}.
\end{equation}

Let $S$ be a finite set such that
$\He_{\mathcal{L}}$ and $\He_{\mathcal{L}'}$ are equal on $\bG^S$. If we
set $\He_{\mathcal{L},S}=\He_{\mathcal{L}}|_{\bG_S}$, then it
follows from (\ref{eq:tauL}) and Theorem \ref{ab} that
\begin{equation}\label{eq:tam}
\frac{\tau_{\mathcal{L}}(\bG(\A))}{\tau_{\mathcal{L}'}(\bG(\A))}=
 \frac{\int_{\bG_S} \He_{\mathcal{L},S}(g)^{-a_L}\gamma_S(g)\,d\tau_S}
{\int_{\bG_S} \He_{\mathcal{L}',S}(g)^{-a_L}\gamma_S(g)\,d\tau_S}.
\end{equation}
Theorem \ref{vol} with $S=\emptyset$ and (\ref{eq:tam}) imply that the
both sides of (\ref{eq:scale})
stay the same when $\He_{\mathcal{L},S}$ and $\He_{\mathcal{L}',S}$ are
replaced by constant multiples. Hence, we can assume that
\begin{equation}\label{eq:e}
\He_{\mathcal{L},S}(e)=\He_{\mathcal{L}',S}(e)=1.
\end{equation}
As in the proof of Lemma \ref{measure}, we obtain
\begin{align*}
\tau(B_T\cap  G_{V_f}) &= \int_{g\in \bG_S}
\tau^S(B_{T\He_{\mathcal{L},S}(g)^{-1}}\cap
s_g G_{V_f}\cap
\bG^S) \; d\tau_S(g)\\
& \sim \left(\int_{g\in \bG_S}
\He_{\mathcal{L},S}(g)^{-a_L}\gamma_S(g)d\tau_S(g)\right)\cdot \tau^S(B_T\cap
G_{V_f}\cap \bG^S).
\end{align*}
 Similarly,
\begin{align*}
\tau(B_T'\cap  G_{V_f})
 \sim \left(\int_{g\in \bG_S} \He_{\mathcal{L}',S}(g)^{-a_L}\gamma_S(g)d\tau_S(g)\right)\cdot \tau^S(B_T'\cap G_{V_f}\cap
\bG^S).
\end{align*}
 Since by (\ref{eq:e}),
$$
B_T\cap \bG^S=B_T'\cap
 \bG^S,
$$
this finishes the proof.
\end{proof}

Now we define a finite measure $\tau_\mathcal L$ on $X$ which describes
the asymptotic distribution of rational points
with respect to a metrized ample line bundle $\mathcal L=
(L,\He_{\mathcal L})$.
If $L$ is saturated, we set $\tau_{\mathcal L}$ to be
$\tilde \tau_{\mathcal L}$
 defined in the following
Proposition. Let $W_{\mathcal L}$ denote
the maximal compact open subgroup of $\bG(\A_f)$ under which
$\He_{\mathcal L}$ is bi-invariant.
\begin{Prop}\label{taul}
For any metrized ample line bundle $\mathcal{L}=(L,\He_\mathcal{L})$,
there exists a unique finite measure $\tilde \tau_{\mathcal L}$
on $X(\A)$ such that for all $\psi\in C(X(\A))$ invariant
under a co-finite subgroup of $W_{\mathcal L}$,
\begin{equation}\label{eq:tauL}
\tilde \tau_{\Cal L}(\psi)=d_K ^{-\op{dim}(X)/2}\cdot \sum_{\chi\in \Lambda}
\left(\lim_{s\to a_L^+} (s-a_L)^{b_L} 
\int_{\bga} \He_\mathcal{L}(g)^{-s}\chi(g) \;\psi(g)\; d\tau(g)\right)
\end{equation}
\end{Prop}
\begin{proof}
Let $\iota$ be the representation
with highest weight given by the regular dominant weight corresponding
to $L$.
Then $X(\A)=X_\iota$.
By Theorem \ref{smuiota-g}
it suffices to set
 $$\tilde \tau_{\mathcal L}=d_K^{-\operatorname{dim}(X)/2}\cdot \gamma_{W_{\mathcal L}}(e)
\cdot
 \tilde \mu_{\iota} .$$
\end{proof}

For a general ample line bundle $L$, the variety $X$ has an {\it asymptotic
arithmetic $\mathcal{L}$-fibration} in the sense of \cite{BT3}.
By the results in section \ref{arith}, there exist
a connected semisimple $K$-group $\bN$ and
a surjective $K$-homomorphism $\pi:\bG\to \bN$
such that
for each $x\in \bN(K)$, there exist a finite measure
$\tilde \tau_{x, \Cal L}$ on $X(\A)$ supported on $\pi^{-1}(x)(\A)$ satisfying the following:
\begin{enumerate}
 \item 
for any $\psi\in C(X(\A))$ invariant under a compact open subgroup of $\bG(\A_f)$,
\begin{equation*}
\tilde \tau_{x,\Cal L}(\psi)=d_K ^{-\op{dim}(X_M)/2}\cdot \sum_{\chi\in \Lambda_{\bM}}
\left(\lim_{s\to a_L^+} (s-a_L)^{b_L} 
\int_{\bM(\A)} \He_\mathcal{L}(gx)^{-s}\chi(g) \;\psi(gx)\; d\tau_{\bM}(g)\right)
\end{equation*}
where $\bM=\pi^{-1}(e)$, $\tau_{\bM}$ is the Haar measure on $\bM(\A)$ with
$\tau_{\bM}(\bM(K)\ba \bM(\A))=1$, $\Lambda_{\bM}$ is defined in the same way as $\Lambda$ for $\bM$
and $X_M$ is the closure of $\bM$ in $X$;
\item there exists $c_x>0$ such that as $T\to \infty$, $$N_{\pi^{-1}(x)}(\Cal L,T)
\sim  c_x
 \cdot T^{a_L} (\log T)^{b_L-1}.$$
\end{enumerate}

By Lemma \ref{p:tam}, there exists $c_L>0$ (independent of metrization)
such that
$c_x=c_L\cdot \tilde \tau_{x,\Cal L}(X(\A))$ for each $x\in \bN(K)$. 
Theorem \ref{uns} implies that
$$\sum_{x\in \bN(K)}\tilde \tau_{x, \Cal L}(X(\A))<\infty $$
and that the following defines a finite measure on $X$
satisfying Theorem \ref{equigeneral}:
\begin{equation*}\label{generaltauL}\tau_{\Cal L}:=\sum_{x\in \bN(K)} \tau_{x, \Cal L}.\end{equation*}

\begin{Thm}\label{equigeneral} 
For any metrized ample line bundle $\mathcal L=(L,\He_{\mathcal L})$
on $X$, and for any $\psi\in C(X(\A))$,
$$\lim_{T\to \infty}\frac{1}{N_\bG(\Cal L, T)} \sum_{g\in \bG(K): \He_\Cal L(g)<T} \psi(g)
= \frac{1}{\tau_\Cal L (X(\A))} \int _{X(\A)} \psi\, d\tau_\Cal L .$$
Moreover, if
 $L$ is saturated and $c_\mathcal L$ is as in
Theorem \ref{manp},
the ratio $\frac{c_{\mathcal L}}{\tau_\mathcal L (X(\A))}$ is independent
of the metrization $\He_{\mathcal L}$. \end{Thm}

\noindent{\bf Remark: } Peyre \cite{Pe1} defined the Tamagawa measure $\tau_{-\Cal K_X}$ on
$X(\A)$ associated with the anti-canonical
metrized line bundle $-\Cal K_X=(-K_X, \He_{-\Cal K_X})$:
$$\tau_{-\Cal K_X}:= c_0\cdot
\lim_{s\to 1^+}
(s-1)^{\op{rank } (\op{Pic}(X))}\left(\prod_{v\in R-S} L_v (s, \op{Pic}(X))\right)
\cdot    \He_{-\Cal K_X}(g)^{-1} d\tau(g) $$ where
$S\subset R$ is a finite subset of places with bad reduction,
and $c_0= d_K ^{-\frac{\op{dim}(X)}{2}}\cdot
 \prod_{v\in R-S}
L_v(1, \op{Pic} (X))^{-1}$
with $d_K$ the discriminant of $K$.

Note that $a_{-K_X}=1$, $b_{-K_X}=\hbox{rank} (\hbox{Pic}(X))$, and
$$
\lim_{s\to 1^+} (s-1)^{\hbox{\tiny rank} (\hbox{\tiny Pic}(X))} \int_{\bga} \He_{-\mathcal{K}_X}(g)^{-s}\chi(g) d\tau(g)=0
$$
for all $\chi \ne 1$. Hence 
 for $\Cal L=-\Cal K_X$,
(\ref{eq:tauL}) gives Peyre's measure $\tau_{-\mathcal{K}_X}$.
An analog of Peyre's measure for general line bundles was also
introduced in \cite{BT3}, but the measure $\tau_\mathcal{L}$ seems to be different, in
general, from the measure defined in \cite{BT3}.

\end{document}